\documentclass[12pt]{amsart}
\usepackage[utf8]{inputenc}

\usepackage{amsmath}
\usepackage{amsxtra}
\usepackage{amscd}
\usepackage{amsthm}
\usepackage{amsfonts}
\usepackage{amssymb}

\usepackage[all]{xy}
\usepackage{bbm}
\usepackage{hyperref}
\usepackage[mathscr]{eucal}

\numberwithin{equation}{section}
\newtheorem{theorem}{Theorem}
\newtheorem{lemme}[equation]{Lemma}
\newtheorem{proposition}[equation]{Proposition}
\newtheorem{definition}[equation]{Definition}
\newtheorem{corollary}[equation]{Corollary}
\newtheorem{conjecture}[equation]{Conjecture}
\newtheorem*{claim}{Claim}

\newtheorem*{corollary*}{Corollary}

\newcommand\calE{{\mathcal E}}

\newcommand\build[3]{\mathrel{\mathop{\kern
0pt#1}\limits_{\textstyle #2}^{\textstyle #3}}}

\newcommand\bbA{{\mathbbm A}}

\newcommand\CC{{\mathbbm C}}

\newcommand\FF{{\mathbbm F}}
\newcommand\GG{{\mathbbm G}}

\newcommand\NN{{\mathbbm N}}
\newcommand\PP{{\mathbbm P}}
\newcommand\QQ{{\mathbbm Q}}

\newcommand\XX{{\mathbbm X}}
\newcommand\ZZ{{\mathbbm Z}}
\newcommand\bbx{{\mathbbm x}}

\newcommand\bbG{\GG}

\newcommand\GL{{\rm GL}}
\newcommand\SO{{\rm SO}}
\newcommand\Tr{{\rm Tr}}
\newcommand\Gal{{\rm Gal}}

\newcommand\pr{{\rm pr}}

\newcommand\Frob{{\rm Frob}}
\newcommand\id{{\rm id}}
\newcommand\Spec{{\rm Spec}}
\newcommand\Mod{{\rm Mod}}
\newcommand\Rep{{\rm Rep}}
\newcommand\Gm{{\mathbb{G}}_{m}}
\newcommand\SL{{\rm SL}}
\newcommand\PGL{{\rm PGL}}
\newcommand\Sp{{\rm Sp}}
\newcommand\SU{\textrm{SU}}

\newcommand\Aut{{\rm Aut}}

\newcommand\Sym{{\rm Sym}}
\newcommand\Lie{{\rm Lie}}
\newcommand\Bun{{\rm Bun}}

\newcommand\ev{{\rm ev}}
\newcommand\ad{{\rm ad}}
\newcommand\rk{{\rm rk}}

\newcommand\Out{{\rm Out}}

\newcommand\hookr\hookrightarrow
\newcommand\hookl\hookleftarrow

\newcommand\inv{{\rm inv}}
\newcommand\isom{\,\smash{\mathop{\hbox to 5mm{\rightarrowfill}}
    \limits^\sim}\,}

\newcommand\geom{{\rm geo}}

\newcommand\Ql{{\overline\QQ_\ell}}

\newcommand\der{\on{der}}

\def\lc{{\it loc. cit}}

\newcommand\calO{{\mathcal O}}

\newcommand\cG{{\mathscr G}}
\newcommand\cZ{{\mathscr Z}}

\newcommand\tilw{{\widetilde{w}}}
\newcommand\tilc{\widetilde{c}}

\newcommand\xcoch{\XX_{*}}
\newcommand\xch{\XX^{*}}

\newcommand\Fl{\textup{Fl}}
\newcommand{\tBun}{\widetilde{\Bun}}

\newcommand\Kum{\textup{Kum}}
\newcommand\Kl{\textup{Kl}}
\newcommand{\Swan}{\textup{Swan}}
\newcommand\triv{\textup{triv}}
\newcommand{\Span}{\textup{Span}}
\newcommand{\mult}{\textup{mult}}
\newcommand{\add}{\textup{add}}
\newcommand{\opp}{{\textup{opp}}}
\newcommand{\coker}{\textup{coker}}

\newcommand\Ad{\textup{Ad}}
\newcommand\St{\textup{Std}}
\newcommand\sep{\textup{sep}}
\newcommand{\GO}{\textup{GO}}
\newcommand{\so}{\mathfrak{so}}
\newcommand{\Cox}{\textup{Cox}}
\newcommand{\Vect}{\textup{Vec}}
\newcommand{\Loc}{\textup{Loc}}
\newcommand{\canon}{\textup{can}}
\newcommand{\gr}{\textup{gr}}

\newcommand{\dR}{\textup{dR}}
\newcommand\Conn{\textup{Conn}}

\newcommand\barB{B^\opp}
\newcommand\barN{U^\opp}

\newcommand\bark{{\overline{k}}}

\newcommand\dualG{{\check{G}}}
\newcommand\dualB{{\check{B}}}
\newcommand\dualT{{\check{T}}}

\newcommand\dualX{{\check{X}}}
\newcommand{\dualg}{\check{\mathfrak{g}}}
\newcommand{\dualt}{\check{\mathfrak{t}}}

\newcommand{\dsigma}{{\check{\sigma}}}

\newcommand\dG{\leftexp{\phi}{\check{G}}}
\newcommand\dB{\leftexp{\phi}{\check{B}}}
\newcommand\dT{\leftexp{\phi}{\check{T}}}

\newcommand\pline{\PP^1_{\backslash\{0,\infty\}}}

\newcommand{\plineC}{\PP^1_{\backslash\{0,\infty\},\CC}}
\newcommand{\Px}{\PP^1\backslash\{x\}}
\newcommand\geompline{\pline\otimes_k\bark}
\newcommand{\tpline}{\widetilde{\PP}^1_{\backslash\{0,\infty\}}}

\newcommand{\una}{\underline{a}}
\newcommand{\unb}{\underline{b}}

\newcommand{\jiao}[1]{\langle{#1}\rangle}
\newcommand{\quash}[1]{}  
\newcommand{\leftexp}[2]{{\vphantom{#2}}^{#1}{#2}}
\newcommand{\twtimes}[1]{\stackrel{#1}{\times}}

\newcommand\frg{\mathfrak{g}}
\newcommand{\frsl}{\mathfrak{sl}}
\newcommand{\frS}{\mathfrak{S}}
\newcommand{\tfrS}{\widetilde{\mathfrak{S}}}

\newcommand\thv{{\theta^\vee}}
\newcommand\gav{{\gamma^\vee}}

\newcommand\Pnth{P_{-\theta}}
\newcommand\gnth{\frg_{-\theta}}

\newcommand{\Png}{P_{-\gamma}}
\newcommand{\Vng}{V_{-\gamma}}
\newcommand{\Qng}{Q_{-\gamma}}
\newcommand{\qng}{q_{-\gamma}}

\newcommand\inert{\mathscr{I}}
\newcommand{\iinf}{\inert_\infty}

\newcommand{\grot}{\bbG^{\textup{rot}}_m}

\newcommand{\Srot}{S^{\textup{rot}}}

\newcommand{\Grsubreg}{\Gr_{\textup{subr}}}
\newcommand\const[1]{\overline{\QQ}_{\ell,#1}}
\newcommand\conv[1]{\stackrel{#1}{\ast}}

\newcommand{\Qlt}{\overline{\QQ}_\ell^\times}

\newcommand{\uAut}{\underline{\textup{Aut}}}

\newcommand{\generic}{{\overline{\eta}}}

\newcommand{\Sat}{\textup{Sat}}

\theoremstyle{definition}
\newtheorem*{example}{Example}
\newtheorem{remark}[equation]{Remark}

\DeclareMathOperator{\Res}{Res}
\DeclareMathOperator{\cind}{c-Ind}
\DeclareMathOperator{\Hecke}{Hecke}
\DeclareMathOperator{\GR}{GR}
\DeclareMathOperator{\Gr}{Gr}
\DeclareMathOperator{\AS}{AS}
\DeclareMathOperator{\Perv}{Perv}
\DeclareMathOperator{\IC}{IC}
\DeclareMathOperator{\Hk}{Hk}
\DeclareMathOperator{\G}{G}

\newcommand{\tensor}{\otimes}

\newcommand{\map}[1]{\stackrel{#1}{\longrightarrow}}
\newcommand{\incl}[1]{\stackrel{#1}{\hookrightarrow}}

\newcommand{\mat}[4]{
 \left(  \begin{matrix} #1 & #2 \\ #3 & #4 \end{matrix} \right)}

\def\Qbar{{\overline{\bQ}}}

\def\cB{\mathcal{B}}
\def\cC{\mathcal{C}}
\def\cD{\mathcal{D}}
\def\cE{\mathcal{E}}

\def\cG{\mathcal{G}}

\def\cI{\mathcal{I}}

\def\cL{\ensuremath{\mathcal{L}}}

\def\cN{\mathcal{N}}
\def\cO{\mathcal{O}}
\def\cP{\mathcal{P}}

\def\cS{\mathcal{S}}
\def\cT{\mathcal{T}}
\def\cU{\mathcal{U}}
\def\cV{\mathcal{V}}

\def\cZ{\mathcal{Z}}

\def\bA{{\bbA}}

\def\bF{{\FF}}
\def\bG{{\GG}}

\def\bN{{\NN}}

\def\bP{{\PP}}
\def\bQ{{\QQ}}

\def\bZ{{\ZZ}}

\def\P{\textsf{P}}

\def\P1min2{\bP^1_{\backslash\{0,\infty\}}}
\def\PSU{\textrm{PSU}}
\def\LG{{{}^L\!\cG}}
\def\der{\textrm{der}}

\begin{document}
	\title{Kloosterman sheaves for reductive groups}
	\author{Jochen Heinloth}
	\address{University of Amsterdam, Korteweg-de Vries Institute for Mathematics, Science Park 904, 1098 XH Amsterdam,
The Netherlands}
           \email{J.Heinloth@uva.nl}
	\author{Bao-Ch\^au Ng\^o}
	\address{School of Mathematics, Institute for Advanced Study, Princeton, NJ 08540, USA}
	\author{Zhiwei Yun}
	\address{School of Mathematics, Institute for Advanced Study, Princeton, NJ 08540, USA}
        \email{zhiweiyun@gmail.com}
	\date{\today}

\begin{abstract}
Deligne constructed a remarkable local system on $\bP^1-\{0,\infty\}$ attached to a family of Kloosterman sums. Katz calculated its monodromy and asked whether there are Kloosterman sheaves for general reductive groups and which automorphic forms should be attached to these local systems under the Langlands correspondence. 

Motivated by work of Gross and Frenkel-Gross we find an explicit family of such automorphic forms and even a simple family of automorphic sheaves in the framework of the geometric Langlands program. We use these automorphic sheaves to construct $\ell$-adic Kloosterman sheaves for any reductive group in a uniform way, and describe the local and global monodromy of these Kloosterman sheaves. In particular, they give motivic Galois representations with exceptional monodromy groups $G_2,F_4,E_7$ and $E_8$. This also gives an example of the geometric Langlands correspondence with wild ramifications for any reductive group.
\end{abstract}

\maketitle

\section*{Introduction}	








\subsection{Review of classical Kloosterman sums and Kloosterman sheaves}\label{ss:review}

Let $n$ be a positive integer and $p$ be a prime number.  For every finite extension $\FF_q$ of $\FF_p$ and $a\in\FF_q^\times$, the Kloosterman sum in $n$-variables in defined as the exponential sum
\begin{equation*}
\Kl_n(a;q):=(-1)^{n-1}\sum_{x_1x_2\cdots x_n=a,x_i\in\FF_q}\exp(\frac{2\pi i}{p}\Tr_{\FF_q/\FF_p}(x_1+x_2+\cdots+x_n)).
\end{equation*}
Kloosterman sums occur in the Fourier coefficients of modular forms. They satisfy the Weil bound
\begin{equation}\label{Weilbound}
|\Kl_n(a;q)|\leq nq^{(n-1)/2}.
\end{equation}
When $n=2$, for any $a\in\overline{\FF}^\times_p$, define the angle $\theta(a)\in[0,\pi]$ to be such that $$2p^{\deg(a)/2}\cos(\theta(a))=\Kl_2(a;p^{\deg(a)}),$$ where $\deg(a)$ is the degree of $a$ over $\FF_p$. Then the angles $\{\theta(a)|a\in\overline{\FF}^\times_p\}$ are equidistributed according to the Sato-Tate measure $\frac{2}{\pi}\sin^2\theta d\theta$ on $[0,\pi]$, as $\deg(a)$ tends to $\infty$. In general, a similar equidistribution theorem for $\Kl_n(a;q)$ was proved by Katz \cite{KatzKloosterman}, using Deligne's results in\cite{WeilII}. 

The above properties of Kloosterman sums were proved using a sheaf-theoretic incarnation. 
Let us recall this construction. In \cite{DeligneSommes}, Deligne considered the diagram
\begin{equation*}
\xymatrix{& \Gm^n\ar[dl]_{\sigma}\ar[dr]^{\pi}\\ \GG_a & & \Gm}
\end{equation*}
Here $\Gm$ is the multiplicative group, $\GG_a\cong\bbA^1$ is the additive group, and $\sigma$ (resp. $\pi$) is the map of taking the sum (resp. product) of the $n$-coordinates of $\Gm^n$. Let $\psi:\FF_p\to\QQ_\ell(\mu_p)^\times$ be a nontrivial  character (here $\mu_p$ is the set of $p$-th roots of unity), and let $\AS_\psi$ be the associated Artin-Schreier local system on the additive group $\GG_a$ over $\FF_p$. Deligne then defined the Kloosterman sheaf as the following complex of sheaves with $\QQ_\ell(\mu_p)$-coefficients on $\Gm$ over $\FF_p$:
\begin{equation*}
\Kl_n:=R\pi_!\sigma^*\AS_\psi[n-1].
\end{equation*}
Fix an embedding $\iota:\QQ_\ell(\mu_p)\hookrightarrow\CC$ such that $\iota\psi(x)=\exp(2\pi ix/p)$ for $x\in\FF_p$. For any $a\in\Gm(\FF_q)=\FF_q^\times$, denote by $\Frob_a$ the geometric Frobenius at $a$, acting on the geometric stalk $(\Kl_n)_{\overline{a}}$. Then, by the Grothendieck-Lefschetz trace formula we have:
\begin{equation*}
\Kl_n(a;q) =\iota\Tr(\Frob_a,(\Kl_n)_{\overline{a}}).
\end{equation*}
In this sense, $\Kl_n$ is a sheaf-theoretic incarnation of the Kloosterman sums $\{\Kl_n(a;q)\}_{a\in\FF_q^\times}$.

In \cite[Th\'eor\`eme 7.4, 7.8]{DeligneSommes}, Deligne proved:
\begin{enumerate}
\item\label{ur} $\Kl_n$ is concentrated in degree $0$, and is a local system of rank $n$.
\item\label{ram0} $\Kl_n$ is tamely ramified around $\{0\}$, and the monodromy is unipotent with a single Jordan block.
\item\label{raminf} $\Kl_n$ is totally wildly ramified around $\{\infty\}$ (i.e., the wild inertia at $\infty$ has no nonzero fixed vector on the stalk of $\Kl_n$), and the Swan conductor $\Swan_\infty(\Kl_n)=1$. 
\item\label{pure} $\Kl_n$ is pure of weight $n-1$ (which implies the estimate \eqref{Weilbound}).
\end{enumerate}

In \cite[\S8.7]{KatzKloosterman}, Katz proved the unicity of Kloosterman sheaves: if a rank $n$ local system $L$ on $\Gm$ satisfies properties (\ref{ram0}) and (\ref{raminf}) listed above, then $L$ is isomorphic to $\Kl_n$ on $\Gm\otimes\overline{\FF}_p$ up to a translation action on $\Gm$. Moreover, in \cite[\S11]{KatzKloosterman}, Katz further studied the global geometric monodromy of $\Kl_n$. Fix a geometric point $\generic$ of $\Gm$ and denote by 
\begin{equation*}
\varphi:\pi_1(\Gm,\generic)\to\GL_n(\QQ_\ell(\mu_p))
\end{equation*}
the monodromy representation associated to the Kloosterman sheaf $\Kl_n$. 
Let $\varphi^\geom$ be the restriction of $\varphi$ to $\pi_1(\Gm\otimes_{\FF_q}\overline{\FF}_q,\generic)$. Katz determined the Zariski closure $\dualG_{\geom}$ 
of the image of $\varphi^\geom$ to be:
\begin{equation}\label{Klnimage}
\dualG_{\geom}=\begin{cases}\Sp_n & n \textup{ even}\\ 
\SL_n & n \textup{ odd}, p\textup{ odd}\\
\SO_n & n \textup{ odd}, n\neq3, p=2\\
G_2 & n=7, p=2.\\                               \end{cases}
\end{equation}
Using Deligne's theorem in \cite{WeilII} this allowed him to deduce that the conjugacy classes for the Frobenius elements $\Frob_a$ are equidistributed according to a "Sato-Tate" measure.

\subsection{Motivation and goal of the paper}\label{ss:goal}
In view of the mysterious appearance of the exceptional group $G_2$ as the global geometric monodromy, Katz asked  (\cite{KatzDE} p. I-5) whether all semisimple groups appear as geometric monodromy of local systems on $\Gm$. Alternatively, are there exponential sums whose equidistribution laws are governed by arbitrary simple groups, and especially by exceptional groups?

In this paper we find a uniform construction of such local systems. For any split reductive group $\dualG$, we will construct a $\dualG$-local system $\Kl_\dualG$ on $\Gm=\pline$ with similar local ramifications as $\Kl_n$, and when $\dualG=\GL_n$, we recover the Kloosterman sheaf $\Kl_n$ of Deligne. We will determine the Zariski closure of its global geometric monodromy (which turns out to be ``large''), prove purity of the sheaf and deduce equidistribution laws. Finally, we give a conjecture about the unicity of such Galois representations (or local systems).

For the purpose of the introduction, let us restrict to the following cases. Assume $G$ is either a split, almost simple and simply-connected group over $k=\FF_q$, or $G=\GL_n$ over $k$. Let $\dualG$ be its Langlands dual group over $\QQ_\ell(\mu_p)$, which is either a split, simple and adjoint group, or $\GL_n$.

The motivation of our construction comes from the Langlands correspondence for the rational function field $K=k(t)$. 
Already in his study of Kloosterman sheaves, Katz \cite{KatzKloosterman} suggested the following: 

``$\cdots${\em It would be interesting to compare this result with the conjectural description of such sheaves, provided by the Langlands philosophy, in terms of automorphic forms}.''

In a series of work \cite{GrossReeder},\cite{Gross},\cite{GrossLetter},\cite{FrenkelGross}, Gross, partly joint with Reeder and with Frenkel, proposed a candidate automorphic representation $\pi$ of $G(\bbA_K)$ which, in the case of $G=\GL_n$, should give the Kloosterman sheaf $\Kl_n$. Let us briefly review their work.

In \cite{GrossReeder}, Gross and Reeder gave the following construction of a representation $V_\phi$ of $G(k((s)))$. Fix a Borel subgroup $B\in G$ and denote by $U$ the unipotent radical of $B$. Denote by $I(0):=\{ g\in G(k[[s]]) | g \mod s \in B\}$ the Iwahori subgroup of $G(k[[s]])$ and by $I(1):=\{g\in G(k[[s]]) | g\mod s \in U\}$ the unipotent radical of $I(0)$. For any affine generic character
$\phi\colon I(1) \to \QQ_\ell(\mu_p)$ (see \S\ref{GenericCharacter}), let $V_\phi:= \cind_{I(1)\times Z(G)}^{G(k((s)))}(\phi\otimes 1)$ be the compactly induced representation of $G(k((s)))$. Gross and Reeder show that this representation is irreducible and supercuspidal.

For any global field $F$, Gross \cite{Gross} managed to use the trace formula to obtain an expression for the multiplicities of automorphic representations $\pi$ of $G(\bbA_F)$ whose ramified local components are either the Steinberg representation or the representation $V_\phi$. In particular,  when the global field is $K=k(t)$, this formula implies that, for any semisimple simply connected group $G$, there is unique cuspidal automorphic representation $\pi=\pi(\phi)$ of $G(\bbA_K)$ such that
\begin{enumerate}
\item[\eqref{ur}] $\pi$ is unramified outside $\{0,\infty\}$;
\item[\eqref{ram0}] $\pi_0$ is the Steinberg representation of $G(k((t)))$;
\item[\eqref{raminf}] $\pi_\infty$ is the simple supercuspidal representation $V_\phi$ of $G(k((s)))$, $s=t^{-1}$.
\end{enumerate}

Motivated by the Langlands philosophy, Gross \cite{GrossLetter} raised the following conjecture: the automorphic representation $\pi=\pi(\phi)$ should correspond to a $\dualG$-local system $\Kl_\dualG(\phi)$ on $\Gm=\pline$, which we will call the {\em Kloosterman sheaf} associated to $\dualG$ and $\phi$, with the following properties parallel to that of $\pi$:
\begin{enumerate}
\item[\eqref{ur}] $\Kl_\dualG(\phi)$ is a $\dualG$-local system on $\Gm=\pline$.
\item[\eqref{ram0}] $\Kl_\dualG(\phi)$ is tamely ramified around $\{0\}$, and the monodromy is a regular unipotent element in $\dualG$.
\item[\eqref{raminf}] The local system $\Kl^\Ad_\dualG(\phi)$ associated to the adjoint representation of $\dualG$ (which is a local system on $\pline$ of rank $\dim\dualG$) is totally ramified around $\{\infty\}$ (i.e., the inert group at $\infty$ has no nonzero fixed vector), and $\Swan_\infty(\Kl^\Ad_\dualG(\phi))=r(\dualG)$, the rank of $\dualG$. 


\item[\eqref{pure}] For any irreducible representation $V$ of $\dualG$, the associated local system $\Kl_\dualG^V(\phi)$ is pure. 
\end{enumerate}

For $\dualG=\GL_n$, the sheaf $\Kl_n$ has the above properties.  For $\dualG$ of type $A,B,C$ and $G_2$, $\dualG$-local systems with the above properties were constructed earlier by Katz \cite{KatzDE}, using a case-by-case construction, and using $\Kl_n$ as building blocks.

In \cite{FrenkelGross}, Gross and Frenkel constructed an analog of such local systems over the complex numbers. Namely, they  defined a $\dualG_{\CC}$-connection $\nabla_\dualG$ on $\PP^1_{\CC}\backslash\{0,\infty\}$ (depending on the choice of nonzero vectors in the affine simple root spaces of $\dualg$). This connection has regular singularity at $\{0\}$ and irregular singularities at $\infty$, parallel to the properties \eqref{ram0} and \eqref{raminf} above. They furthermore show that $\nabla_\dualG(\dualX)$ is cohomologically rigid \cite[Theorem 1]{FrenkelGross}, and compute the differential Galois group of these connections. For $\dualG$ of type $A,B,C$ and $G_2$, they verify that $\nabla_\dualG$ coincide with connections constructed by Katz \cite{KatzDE} which are the analogs of Katz's $\ell$-adic $\dualG$-local systems mentioned above. All these give strong evidence that $\nabla_\dualG(\dualX)$ should be the correct de Rham analog of the conjectural local system $\Kl_\dualG(\phi)$. 

The predictions about the conjectural local system $\Kl_\dualG(\phi)$ made in \cite{FrenkelGross} served as a guideline for our work.



\subsection{Method of construction}
Our construction of the Kloosterman sheaves can be summarized as follows.

We start with the automorphic representation $\pi$ mentioned above. The key observation is, that $\pi$ contains a Hecke eigenfunction $f_\phi$, which can be written down explicitly using the combinatorics of the  double coset $G(K)\backslash G(\bbA_K)/I_0\times I_\infty(2)\times\prod_{x\neq0,\infty}G(\calO_x)$ (\S\ref{ss:eigenfn}).


The points of this double coset are the rational points of a moduli stack $\Bun_{G(0,2)}$ of $G$-bundles  on $\bP^1$ with a particular level structure at $0$ and $\infty$ and we can upgrade the function $f_\phi$ to a sheaf $A_\phi$ on this stack. We then prove that this sheaf is indeed a Hecke eigensheaf. In particular we show that the eigenvalues of the geometric Hecke operators applied to $A_\phi$ define a $\dualG$-local system on $\pline$ -- our Kloosterman sheaf $\Kl_\dualG(\phi)$. Technically, the proof of this result relies on the fact that the sheaf $A_\phi$ is a clean extension of a local system with affine support.

In the case $G=\GL_n$ we compute the local system $\Kl_\dualG(\phi)$ explicitly. For a particular $\phi$, the rank $n$ local system associated to $\Kl_\dualG(\phi)$ and the standard representation of $\dualG=\GL_n$ turns out to be the Kloosterman sheaf $\Kl_n$ of Deligne. Here Deligne's diagram will reappear as a part of a geometric Hecke transformation.

Our construction also works when we replace $G$ by a certain class of quasi-split group schemes $\cG$ over $\pline$ and there is also a variant depending on an additional Kummer character. These generalization seems to be the natural setup needed in order to compare Kloosterman sheaves for different groups. Since the generalization does not require additional arguments, we will give the construction in this more general setup. 



\subsection{Properties of Kloosterman sheaves}
The longest part of the paper (\S\ref{ProofThm}--\S\ref{s:mono}) is then devoted to the study of local and global properties of Kloosterman sheaves $\Kl_{\dualG}(\phi)$ for simple adjoint groups $\dualG$.

First, using \cite{WeilII}, the expected purity property \eqref{pure} is a corollary to our construction.
When $\dualG$ is adjoint we can moreover normalize $\Kl^V_\dualG(\phi)$ to be pure of weight $0$. 
Thus, fixing an embedding $\iota:\QQ_\ell(\mu_p)\hookrightarrow\CC$ as before, the exponential sums
\begin{equation*}
\Kl^V_\dualG(\phi;a;q)=\iota\Tr(\Frob_a, (\Kl^V_\dualG(\phi))_{\overline{a}}),  \textup{ for }a\in\FF_q^\times, V\in\Rep(\dualG) 
\end{equation*}
satisfy the Weil bound
\begin{equation*}
|\Kl^V_\dualG(\phi;a;q)|\leq\dim V.
\end{equation*}

The property \eqref{ram0} about the monodromy at $0$ (see Theorem \ref{Heckeeigensheaf}(2)) will be prove in \S\ref{ss:unipmono}.

The property \eqref{raminf} about the monodromy at $\infty$ (see Theorem \ref{th:Swan-inv}) will be proved in \S\ref{s:coho}. The calculation of the Swan conductors involves a detailed study of the geometry of certain Schubert varieties in the affine Grassmannian, which occupies a large part of \S\ref{s:coho}. Using a result of Gross-Reeder \cite{GrossReeder}, one can give an explicit description of the monodromy at $\infty$ (Corollary \ref{c:simple-wild}).
Together these results show that the local system $\Kl_\dualG(\phi)$ satisfies the properties (1)-(4) expected by Gross.

Again, let
\begin{equation*}
\varphi:\pi_1(\pline,\generic)\to\dualG(\QQ_\ell(\mu_p))
\end{equation*}
be the monodromy representation associated to the Kloosterman sheaf $\Kl_\dualG(\phi)$, and let $\varphi^\geom$ be its restriction to $\pi_1(\geompline,\generic)$. 
We find (Theorem \ref{th:glob}) that for $p>2$ (and $p>3$ in the case $\dualG=B_3$) the Zariski closure $\dualG_{\geom}$ of the image of $\varphi^\geom$
coincides with the differential Galois group $\dualG_{\nabla}$ of $\nabla_\dualG(\dualX)$ calculated by Frenkel-Gross in \cite[Cor.9,10]{FrenkelGross}, which we list in Table \ref{geomMon}.
\begin{table}[h]
\caption{Global geometric monodromy of $\Kl_\dualG(\phi)$}\label{geomMon}
\begin{tabular}{l|l}
$\dualG$ & $\dualG_{\geom}$\\\hline
$A_{2n}$ & $A_{2n}$ \\
$A_{2n-1}, C_n$ &  $C_n$ \\ 
$B_n, D_{n+1}$ $(n\geq4)$ & $B_n$ \\
$E_7$ & $E_7$\\
$E_8$ & $E_8$\\
$E_6, F_4$ & $F_4$\\
$B_3,D_4, G_2$ & $G_2$.
\end{tabular}
\label{t:globalmono}
\end{table}


Using Katz's argument (\cite[\S3,\S13]{KatzKloosterman}) this implies the following 
equidistribution law for the conjugacy classes $\varphi(\Frob_a)$.

\begin{corollary*}\label{th:Equidist}
Suppose $\dualG$ is simple and adjoint (in this case the whole image of $\varphi$ lies in $\dualG_\geom$). Let $\dualG_{\geom,c}\subset\dualG_\geom(\CC)$ be a compact real form (using the embedding $\iota:\QQ_\ell(\mu_p)\hookrightarrow\CC$ to make sense of $\dualG_\geom(\CC)$). Let $\dualG^\natural_{\geom,c}$ be the set of conjugacy classes of $\dualG_{\geom,c}$.

For $a\in\overline{k}^\times$, the conjugacy class $\varphi(\Frob_a)$ of $\dualG_{\geom}(\CC)$ in fact belongs to  $\dualG^\natural_{\geom,c}$. As $\deg(a)$ (the degree of the field extension $k(a)/k$) tends to $\infty$, the conjugacy classes 
$$\{\varphi(\Frob_a)|a\in\overline{k}^\times\}\subset\dualG^\natural_{\geom,c}$$ 
become equidistributed according to the push-forward of the Haar measure from $\dualG_{\geom,c}$ to $\dualG_{\geom,c}^\natural$.
\end{corollary*}

\subsection{Open problems}
Compared to the known results about $\Kl_n$, the only missing piece for $\Kl_\dualG(\phi)$ is the unicity. 
We prove that $\Kl_{\dualG}(\phi)$ is {\em cohomologically rigid}, i.e.,
\begin{equation*}
H^1(\PP^1,j_{!*}\Kl^\Ad_\dualG(\phi))=0.
\end{equation*}
This gives evidence for the {\em physical rigidity} of $\Kl_\dualG(\phi)$: any other local system $L$ on $\pline$ satisfying properties \eqref{ram0}  and \eqref{raminf} should be isomorphic to $\Kl_\dualG(\phi)$ over $\geompline$. We will state the precise unicity conjecture in Conjecture \ref{conj:weak} and \ref{conj:strong}.

Table \ref{t:globalmono} also suggests that the Kloosterman sheaves for different groups $\dualG$ appearing in the same line of the table should be essentially the same (see Conjecture \ref{conj:Hasse-Davenport}). This can be viewed as a functoriality statement for Kloosterman sheaves.

Also, the monodromy of Kloosterman sheaves is studied in detail only for split groups. 
In \S\ref{ss:conjqs}, we state our predictions on the local and global monodromy of $\Kl_\LG(\phi)$ for $\cG$ a quasi-split, simple and simply-connected group scheme over $\PP^1$ with good reduction at all places outside $\{0, \infty\}$, which is split by a tame Kummer cover $[N]:\pline\to\pline$ ($z\mapsto z^N$).

\subsection{Organization of the paper}
In \S\ref{GroupSchemes}, we define the various group schemes $\cG(m_0,m_\infty)$ over $\PP^1$ encoding the level structures of the moduli stacks of $G$-bundles on which the automorphic sheaves will be defined. We also give a description of the geometry of these moduli stacks. 
In \S\ref{s:eigensheaf}, we construct the automorphic sheaf $A_\phi$ 
and state 
the main results of the paper. 
The proofs are given in \S\ref{ProofThm}--\S\ref{s:mono}. \S\ref{s:exam} is devoted to the case $G=\GL_n$, where we recover classical Kloosterman sheaves of Deligne. In \S\ref{s:fonc}, we state our conjectures about rigidity of Kloosterman sheaves, and our expectations about Kloosterman sheaves for quasi-split groups.
 
The paper contains four appendices. In \S\ref{a:moduli}, we prove the structure theorem for the moduli spaces $\Bun_\cG$ in the generality of quasi-split groups. 
In \S\ref{a:Sat}, we prove compatibility between outer automorphisms and the geometric Satake equivalence. This is needed to describe the global monodromy image of Kloosterman sheaves and we couldn't find a reference for the result. In \S\ref{a:qmcom}, we collect some facts about quasi-minuscule representations of $\dualG$. In \S\ref{a:g2}, we analyze the geometry of the adjoint Schubert variety for $G_2$, which does not fit into the uniform treatment for groups of other types.

\noindent{\bf Acknowledgement} The authors met and started the project during their stay at the Institute for Advanced Study in 2009-2010. They would like to thank the IAS for excellent working condition.

At various stages of the project, one of the authors (Z.Y.) had numerous discussions with Dick Gross, who provided us with invaluable insight and advice. The authors also benefitted from conversations with P.\ Deligne, D.\ Gaitsgory, N.\ Katz, M.\ Reeder and P.\ Sarnak.

The work of J.H.\ is supported by NSF (under the agreement No DMS-0635607) and the GQT-cluster of NWO.
The work of B-C.N.\ is supported by NSF and the Simonyi foundation.
The work of Z.Y.\ is supported by NSF (under the agreement No.DMS-0635607) and Zurich Financial Services.

\section{Structural groups}\label{GroupSchemes}
%
We will work over a fixed finite field $k$ of characteristic $p$. 
We fix a coordinate $t$ of our base curve $\bP^1$, so that $\bA^1=\Spec(k[t])\subset \bP^1$. We write $s:=t^{-1}$ for the coordinate around $\infty\in \bP^1$. 

For any closed point $x\in \bP^1$ we will denote by $\cO_x$ the completed local ring at $x$ and by $K_x$ the fraction field of $\cO_x$. 

Since we are interested in the geometric Langlands correspondence for (wildly) ramified local systems we will need to consider principal bundles with various level structures. It will be convenient to view these as torsors under group schemes $\cG$ over $\bP^1$. Moreover, in order to formulate the conjectured functoriality, it will be useful to allow quasi-split group schemes. We will introduce these group schemes in several steps.

\subsection{\texorpdfstring{Quasi-split group schemes over $\pline$}{Quasi-split group schemes over P1-(0,infinity)}}\label{ss:qs}
We will assume that $\cG|\pline$ is a quasi-split reductive group. Moreover we will assume that there is a finite extension $k^\prime/k$ and an integer $N$ with $(N,p)=1$, such that $\cG$ splits over the tame extension $[N]\colon \bG_{m,k^\prime}\to \bG_{m,k^\prime}$ defined by $t\mapsto t^N$. We may and will assume that $k^\prime$ contains all $N-$th roots of unity. We write $\mu_N$ for the group of $N$-th roots of unity.

We will fix subgroups $S \subset \cT \subset \cB \subset \cG$, where $S$ is a maximal split torus, $\cT$ is a maximal torus and $\cB$ is a Borel subgroup. We also fix a quasi-pinning of $\cG$.

In order to describe our group schemes, we will for simplicity assume that $k=k^\prime$. In general, the construction we give will be invariant under $\Gal(k^\prime/k)$, so that Galois-descent will then give the general case.

By \cite{SGA3III} (Exp. XXIV, Thm 3.11) these groups can be described explicitly as follows. 
By assumption there is a split reductive group $G$ over $k$ such that $[N]^*\cG|\bG_{m}=G\times \bG_{m}$. The automorphism group of the covering $[N]$ is $\mu_N$. We fix $T\subset B \subset G$ a split maximal torus and a Borel subgroup as well as a pinning $\dagger$ of $G$.

This data defines a morphism $\sigma\colon\mu_N \to \Aut^\dagger(G)$, where $\Aut^\dagger(G)$ is the automorphisms of $G$ respecting the pinning. Let us denote by $[N]_*(G\times \bG_m)$ the Weil restriction of $G\times \bG_m$ for the covering $[N]$. Then $\mu_N$ acts on this Weil restriction by the action of $\sigma$ on $G$ and by the action on the covering $[N]$. We will denote the diagonal action again by $\sigma$. Then descent implies:
$$\cG|\pline\cong ([N]_*(G\times\pline))^{\sigma} = (G\times\bG_m)/(\mu_N).$$

\subsection{Level structures at $0$ and $\infty$}

Next we need to describe several extensions of $\cG|\pline$ to $\bP^1$. We will again denote by $[N]:\bP^1\to \bP^1$ the map given by $t\mapsto t^N$.

The group $\cG|\pline$ will be extended as the (special) Bruhat-Tits group $\cG$ defined by the convex function, which is $0$ on all roots. Since we assumed that $\cG|\pline$ splits over a tame extension this group can be described as the connected component of the group $([N]_*(G\times \bP^1))^{\sigma}$. The subgroups $S\subset \cT|\pline \subset \cB|\pline \subset \cG|\pline$ define closed subgroups $S\subset \cT \subset \cB\subset \cG$. 

We will need to introduce level structures at $0$ and $\infty$. These will correspond to the first steps of the Moy-Prasad filtration of $\cG$. A self-contained exposition of the construction of these group schemes can be found in the preprint of Yu \cite{Yu}.

We first consider tori. For a split torus $T=\bG_m^r \times \bP^1$ we will denote by $T(m_0,m_\infty)$ the smooth group scheme over $\bP^1$ such that 
\begin{align*} 
T(m_0,m_\infty)|{\pline}&=\Gm^r\times \pline\\
T(m_0,m_\infty)(\cO_{\bP^1,0})&=\{ g\in T(\cO_{\bP^1,0}) | g \equiv 1 \mod t^{m_0} \}\\
T(m_0,m_\infty)(\cO_{\bP^1,\infty})&=\{ g\in T(\cO_{\bP^1,\infty}) | g \equiv 1 \mod s^{m_\infty} \}.
\end{align*}
This also defines a filtration for induced tori. For an arbitrary torus $\cT$, pick an embedding into an induced torus $\cT\incl{} \cI$ and define the filtration by pulling back the filtration on $\cI$. 
By \cite{Yu} Section 4, this definition is independent of the chosen embedding and since we assumed that $\cT$ splits over a tame extension the subgroups define connected groups.

For reductive groups $\cG|\pline$ we want to define models $\cG(m_0,m_\infty)$ for $m_0,m_\infty\in\{0,1,2\}$, corresponding to the $m_0$-th (resp. $m_{\infty}$-th) step in the Moy-Prasad flitration of the Iwahori subgroup at $0$ (resp. $\infty$).

First assume that $\cG=G\times \bP^1$ is a split, semisimple group. As before, we fix $S=T\subset B \subset G$ a split maximal torus and a Borel subgroup and a pinning of $G$. We denote by $U\subset B$ the unipotent subgroup.
Denote by $\Phi=\Phi(G,S)$ the set of roots and $\Phi^\pm\subset \Phi$ the set of positive and negative roots with respect to the chosen Borel subgroup. Let $\alpha_1,\dots \alpha_r$ denote the positive, simple roots and by $\{\theta_j\}$ the highest roots of $G$. Finally for each root $\alpha$ let $U_\alpha$ denote the corresponding root subgroup of $G$. 

We will consider the following bounded subgroups of $G(k[[s]])$:
\begin{align*}
I(0)&:=\{g\in G(k[[s]]) | g \mod s \in B\} \textrm{ is the Iwahori subgroup.}\\
I(1)&:=\{g\in G(k[[s]]) | g \mod s \in U\} \textrm{ is the unipotent radical of } I(0).
\end{align*}

To describe $I(2)$ let $f\colon \Phi\cup\{0\} \to \bN$ be the concave function defined by $f(0)=1$ and
$$ f(\alpha) = \left\{ \begin{array}{cl} 
0 & \textrm{ if } \alpha \in \Phi^+\setminus\{\alpha_i\} \textrm{ is positive, but not simple}\\
1 & \textrm{ if } \alpha=\alpha_i \textrm{ is a positive simple root}\\
1 & \textrm{ if } \alpha\in \Phi^- \setminus\{ - \theta_j\}\\
2 & \textrm{ if } \alpha=-\theta_j \textrm{ is the negative of a highest root.}\\
\end{array}\right.$$

$ I(2)\subset G(k[[s]])$ denotes the bounded subgroup defined by the concave function $f$,
i.e., the subgroup generated by $\{u \in U_\alpha | u \equiv 1 \mod s^{f(i)} \}$ and $\{ g\in T(k[[s]]) | g \equiv 1\mod s \}$. 

Note that by definition $I(1)/I(2)\cong \oplus_{\alpha \textrm{ simple affine}} \bG_a$ is isomorphic to the sum of the root subgroups $U_\alpha\subset G(k[[s]])$ for which $\alpha$ is a simple affine root.

Similarly we define $I(i)^\opp\subset G(k[[t]])$ to be the analogous groups obtained by using the opposite Borel subgroup $B^\opp$ in the above definition. 

\begin{example} For $G=\SL_n$ we choose $S=T$ to be the diagonal matrices and $B$ the upper triangular matrices. Then the subgroup $I(0)\subset SL_n(k[[s]])$ is the subgroup of matrices, such that the lower diagonal entries are divisible by $s$. $I(1)\subset I(0)$ is the subgroup such that the diagonal entries are elements of $1+sk[[s]]$

The root subgroups of $\SL_n$ for the simple roots are given by the above-diagonal entries $a_{i,i+1}$ and the $U_{-\tilde{\alpha}}$ of the negative of the longest root is given by the entry in the lower left corner. So $I(2)$ consists of matrices of the form 
{\small $\left(\begin{array}{cccc} 
1+sk[[s]] & sk[[s]] & k[[s]] & k[[s]] \\
sk[[s]]   & 1+sk[[s]] & sk[[s]] & k[[s]] \\
sk[[s]] &   sk[[s]]  & 1+sk[[s]] & sk[[s]] \\ 
s^2k[[s]] & sk[[s]] & sk[[s]] & 1+k[[s]]
\end{array}\right)$}.
We obtain an isomorphism $I(1)/I(2)\cong \bA^n$ by mapping a matrix $(a_{ij})$ to the leading coefficients of the entries $a_{i,i+1}$ and $a_{n,1}$. 
\end{example}

For a general, split reductive group $\cG=G\times \bP^1$ we consider the derived group $G_{\der}$ and the connected component of the center $Z(G)^\circ$, which is a torus. We will temporarily denote by $I_{G_\der}(i),I_{G_\der}^\opp(i)$ the groups defined above for the semisimple group $G_\der$ and define
$$I(i):= Z(G)^\circ(i)(k[[s]]) \cdot I_{G^{\der}}(i)\subset \cG(k[[s]])$$ and similarly $I(i)^\opp.=Z(G)^\circ(i,m_\infty)(k[[t]]) I^\opp_{G_\der}(i)$.  

For split, reductive groups $\cG=G\times \bP^1$ the group $G(m_0,m_\infty)$ denotes the Bruhat Tits group scheme such that
\begin{align*}
G(m_0,m_\infty)|\pline &= G \times \bG_m\\
G(m_0,m_\infty)(\cO_{\infty})&=I(m_\infty) \textrm{ and }\\
G(m_0,m_\infty)(\cO_{0})&=I(m_0)^\opp.
\end{align*}

Finally for a general quasi-split group we define 
$$\cG(m_0,m_\infty):= \big(\left([N]_*G(m_0,m_\infty)\right)^{\sigma}\big)^\circ.$$ 
Note that for tori this does not give a new definition. 

We will abbreviate:
\begin{align*} 
I(m_\infty)&:= \cG(m_0,m_\infty)(\cO_\infty).\\
I^-(m_0)&:=\cG(m_0,m_\infty)(\bP^1-\{\infty\}).
\end{align*}

Recall from \cite{PappasRapoport} that the groups $I(m)$ have a natural structure as (infinite dimensional) group schemes over $k$.

\subsection{Affine generic characters}\label{GenericCharacter}\label{h-generic}
As indicated before our construction depends on the choice of a character of $I(1)$.
We call a linear function $\phi\colon I(1)/I(2) \to \bA^1$ {\em generic} if for any simple affine root $\alpha$, the restriction of $\phi$ to $U_\alpha$ is non-trivial. 
Throughout we will fix such a generic $\phi$. 

We will fix a non-trivial additive character $\psi\colon \bF_p \to \Qbar_\ell$ and denote the character $k \to \Qbar_\ell$ defined as $\psi\circ \Tr_{k/\bF_p}$ again by $\psi$.
With this notation, the character $\psi\circ \phi$ is called an {\em affine generic character} of $I(1)$.

\subsection{Principal bundles}\label{ss:moduli}

Having defined our groups, we need to collect some basic results on the geometry of the moduli stacks of $\cG$-bundles $\Bun_\cG$. All of these are well-known for constant groups (see e.g., \cite{Faltings} for a recent account). In order to generalize these results to our setup we rely on \cite{PappasRapoport} and \cite{HainesRapoport}, where the corresponding results on twisted loop groups are explained. Let us point out that, except for the computation of the connected components of $\Bun_\cG$, all results are particular to group schemes over $\bP^1$, that split over a tamely ramified covering $[N]:\bP^1\to \bP^1$.

First, we recall some results and notations from \cite{HainesRapoport}. We denote by $\cN$ the normalizer of $\cT$ and by $\widetilde{W}:= \cN(k((s)))/\cT(k[[s]])$ the Iwahori-Weyl group. Furthermore denote by $W_0:=\cN(k((s)))/\cT(k((s)))$ the relative Weyl group. This is isomorphic to the Weyl group of the reductive quotient of the special fiber $\cG_{\infty,red}$ (loc.cit.\ Proposition 13) and $\widetilde{W}\cong \xcoch(\cT)_{\pi_1(\pline)}\ltimes W_0$. 

Denote by $W_a$ the affine Weyl group of the root system of $\cG(k((s)))$, which can be identified with the Iwahori-Weyl group of the simply connected cover of the derived group of $\cG$ (\cite{HainesRapoport} p.196). This is a Coxeter group.
Then $\widetilde{W}\cong W_a \ltimes \Omega$, where $\Omega \cong \xch(Z(\check{G})^{\pi_1(\bG_m)})$ is the stabilizer of an alcove.

Finally, for $x\in \bP^1-\{\infty\}$ we denote by $\Gr_{\cG,x}$ the affine Gra\ss mannian (see \cite{PappasRapoport}) so that $\Gr_{\cG,x}(k)=\cG(K_x)/\cG(\cO_x)$. It can also be defined as the ind-scheme parametrizing $\cG$-bundles $\cP$ on $\bP^1$ together with a trivialization $\varphi\colon\cP|_{\bP^1\backslash\{x\}}\map{\cong}\cG_{\bP^1\backslash\{x\}}$ (e.g.,\cite{Uniformization}).

\begin{proposition}\label{P1-uniformization} Assume that the ground field $k$ is either finite, or algebraically closed.
Let $\cG$ over $\bP^1$ be a quasi-split group scheme as defined in \S\ref{GroupSchemes}. Then the following holds:
\begin{enumerate}
\item For any $x\in \bP^1$, the canonical map $\Gr_{\cG,x} \to \Bun_\cG$ has sections, locally in the smooth topology on $\Bun_\cG$. Moreover, this map is essentially surjective on $k$ points, i.e., for any $\cG$-bundle over $\bP^1$ its restriction to $\bP^1\backslash\{x\}$ is trivial.
\item $\pi_0(\Bun_\cG)\cong \pi_1([N]^*\cG|\bG_m)_{\pi_1(\bG_m)}\cong \Omega$.
\item Every $\cG$-bundle on $\bP^1$ admits a reduction to $\cT$.
\item (Birkhoff-Grothendieck decomposition) 
$$\cG(0,0)(k((s)))= \coprod_{w\in\widetilde{W}} I^-(0) \cdot \widetilde{W} \cdot I(0).$$
\end{enumerate}
\end{proposition}
If $\cG=G\times \bP^1$ is a split group, a proof of this result can be found in \cite{Faltings}. Since the case of split groups was our starting point, we will postpone the proof of the general case to the appendix. For $\gamma\in \Omega$ we will denote by $\Bun_{\cG(m,n)}^\gamma$ the corresponding connected component of $\Bun_{\cG(m,n)}$.

Let us collect some consequences of this result:
\begin{corollary}\label{HeckeOmega} Let $n\in\{0,1,2\}, m\in\{0,1\}$. 
Any $\gamma\in\Omega$ defines an isomorphism $\Hk_\gamma\colon \Bun_{\cG(m,n)} \to \Bun_{\cG(m,n)}$. 
This induces an isomorphism of the connected components $\Bun_{\cG(m,n)}^0 \to \Bun_{\cG(m,n)}^\gamma$. 

In particular all connected components of $\Bun_{\cG(m,n)}$ are isomorphic.
\end{corollary}
\begin{proof}
By part (2) of the proposition, the connected components of $\Bun_{\cG}(m,n)$ are indexed by $\Omega\subset \widetilde{W}$. 

Since $\Omega$ is the stabilizer of the alcove defining $I(0)^\opp$ any $\gamma\in \Omega$ normalizes $I(0)^\opp$ and its unipotent radical $I(1)^\opp$. Thus $\Omega$ acts by right-multiplication on $\Gr_{\cG(m,n),0}=\cG(k((t)))/I(m)^{\opp}$. This defines an isomorphism $\Hk_{\gamma}\colon \Bun_{\cG(m,n)} \to \Bun_{\cG(m,n)}$. Moreover, \cite{HainesRapoport} Lemma 14 implies that the restriction of this map to $\Bun_{\cG(m,n)}^0$ induces an isomorphism $\Bun_{\cG(m,n)}^0 \to \Bun_{\cG(m,n)}^\gamma$.
\end{proof}

\begin{corollary}\label{affine-embedding}
\begin{enumerate}
\item Denote by $\cG(\bP^1)$ the automorphism group of the trivial $\cG$-bundle. 
Then the inclusion $B \cG(\bP^1) \to \Bun_{\cG}$ is an affine, open embedding.
\item Denote by $\cT(\bP^1)=H^0(\bP^1,\cT)=H^0(\bP^1,\cG(0,0))$ the automorphism group of the trivial $\cG(0,0)$-bundle. 
The inclusion of the trivial bundle defines an affine open embedding $B(\cT(\bP^1)) \hookrightarrow \Bun_{\cG(0,0)}$.
\item The trivial $\cG(0,1)$-bundle defines an affine, open embedding $\Spec(k) \hookrightarrow \Bun_{\cG(0,1)}^0$.
\item Applying the action of $I(1)/I(2)$ on $\Bun_{\cG(0,2)}$ to the trivial $\cG(0,2)$-bundle, we obtain a canonical map $j_0\colon I(1)/I(2) \hookrightarrow \Bun_{\cG(0,2)}$. This is an affine open embedding. 
\item For any $\gamma\in \Omega$ the map $$j_\gamma:=\Hk_\gamma\circ j_0\colon I(1)/I(2) \hookrightarrow \Bun_{\cG(0,2)}$$ is an affine open embedding, called the {\em big cell}.
\item Applying the action of $\cT_0^{\textrm{red}}\times I(1)/I(2)$ on $\Bun_{\cG(1,2)}$ to the trivial $\cG(1,2)$-bundle we obtain canonical affine embeddings: \begin{align*} \tilde{j}_0\colon \cT_0^{\textrm{red}}\times I(1)/I(2) &\hookrightarrow \Bun_{\cG(1,2)} \textrm{ and }\\ \tilde{j}_\gamma :=\Hk_\gamma\circ \tilde{j}_0\colon \cT_0^{\textrm{red}}\times I(1)/I(2) &\hookrightarrow \Bun_{\cG(1,2)}. \end{align*}
\end{enumerate}
\end{corollary}
\begin{proof}
Let us first recall, why this corollary holds for $\GL_n$. For (1) note that the only vector bundle on $\bP^1$ of rank $n$ with trivial cohomology is the bundle $\cO(-1)^n$. Therefore the inverse of the determinant of cohomology line bundle on $\Bun_{\GL_n}$ has a section vanishing precisely on the trivial bundle. This proves (1) in this case. Next, recall that a $\GL_n(0,0)$ bundle is a vector bundle together with full flags at $0$ and $\infty$. Let us denote by $\Mod_{i,0}:\Bun_{\GL_n(0,0)}^0 \to \Bun_{\GL_n(0,0)}^{i}$ the $i$-th upper modification along the flag at $0$. The inverse of this map is given by the $i$-th lower modification, which we will denote by $\Mod_{-i,0}$. To prove (2) denote by $\Mod_{i,\infty}$ the $i$-th modification at $\infty$. 
The trivial $\GL_n(0,0)$-bundle is the bundle given by $\cO^n$ and the canonical opposite flags $(V_{i,0})_{i=0,\dots,n}$ of the fiber $(\cO^n)_0$ at $0$ and $(V_{i,\infty})_{i=0,\dots,n}$ at $\infty$. 
This is the only $\GL_n(0,0)$-bundle $(\cE,V_{i,0},V_{i,\infty})$ of degree $0$ on $\bP^1$ such that the complex $\cE \to \cE_0/V_{i,0}\oplus \cE_{\infty}/V_{n-i,\infty}$ has trivial cohomology for all $i$. Thus again, the inclusion of the trivial $\GL_n(0,0)$ bundle in $\Bun_{\GL_n(0,0)}$ is defined by the non-vanishing of sections of line bundles. 

For general $\cG$ we pick a faithful representation $\rho\colon \cG \to \GL_{n}\times \bP^1$. This defines a map $\Bun_{\cG}\to \Bun_{\GL_n}$. The Birkhoff-Grothendieck decomposition implies that a $\cG$ bundle is trivial, if and only if the associated $\GL_n$ bundle is trivial. Moreover, in order to check that the reductions to $\cB$ at $0,\infty$ are opposite it is also sufficient to check this on the induced $\GL_n$-bundle. This proves (1) and (2).

(3) follows from (2) and the isomorphism $H^0(\bP^1,\cT)\cong \cT_0^{red}$ from the proof of Proposition \ref{P1-uniformization}. (4) follows from this, because the map $\Bun_{\cG(0,2)}\to \Bun_{\cG(0,1)}$ is an $I(1)/I(2)$-torsor. By corollary \ref{HeckeOmega} (5) follows from (4). Finally (6) follows from (5) since $\Bun_{\cG(1,2)}\to \Bun_{\cG(0,2)}$ is a $\cT_0^{\textrm{red}}$-torsor.
\end{proof}

\section{Eigensheaf and eigenvalues---Statement of main results}\label{s:eigensheaf}
%

In this section we will construct the automorphic sheaf $A_\phi$. 
We will also state our main results about the local and global monodromy of the Kloosterman sheaf $\Kl_\LG(\phi,\chi)$, which will be defined to be the eigenvalue of $A_\phi$ 

\subsection{The eigenfunction}\label{ss:eigenfn}
In this subsection we give a simple formula for an eigenfunction in Gross' automorphic representation $\pi$. Surprisingly this calculation turns out to be independent of Gross' result. 

We will use the 
generic affine character $\psi\circ\phi\colon I(1)/I(2) \to \Qbar_\ell$ chosen in Section \ref{GenericCharacter}.

By Proposition \ref{P1-uniformization} we know
\begin{align*}
\Bun_{\cG(0,2)}(k) &= \cG(k(t)) \backslash \cG(0,2)(\bA_{k(t)}) / \prod_{x} \cG(0,2)(\cO_x)\\
&  = I^-(0) \backslash \cG(k((s))) / I(2).
\end{align*}
and $$\cG(k((s)))=\coprod_{w\in \widetilde{W}} I^-(0) w I(1).$$

Suppose we were given, as in the introduction, an automorphic representation $\pi=\tensor^\prime \pi_\nu$ of $\cG(\bA_{k(t)})$ such that for $\nu\not\in \{0,\infty\}$ the local representation $\pi_\nu$ is unramified, $\pi_0$ is the Steinberg representation and such that $\pi_\infty$ occurs in $\cind_{Z(\cG(0,1))\times I(1)}^{\cG(k((s)))} \psi\circ\phi$. Then there exists a function $f$ on $\Bun_{\cG(0,2)}(k)$ such that $f(gk)=\psi(\phi(k))f(g)$ for all $k\in I(1)$. 

The following lemma characterizes such functions in an elementary way:

\begin{lemme}\label{eigenfunction}
Let $f\colon \Bun_{\cG(0,2)}(k) \to \Qbar_\ell$ be a function such that for all $k\in I(1), g\in  \cG(0,2)(k((s)))$ we have $f(gk)=\psi(\phi(k)) f(g)$. Then $f$ is uniquely determined by the values $f(\gamma)$ for $\gamma \in \Omega$. Moreover $f(w)=0$ for all $w\in \widetilde{W}-\Omega$.
\end{lemme}
\begin{proof}
For any $w\in \widetilde{W}$ with $l(w)>0$ there exists a simple affine root $\alpha_i$ such that $w.\alpha_i$ is negative. This implies that $wU_{\alpha_i} w^{-1}\subset I^-(0)$. This implies that for all $u\in U_{\alpha_i}$ we have $\psi(\phi(u))f(w)=f(wu)=f(w)$ so $f(w)=0$. For $w$ with $l(w)=0$, i.e. $w\in \Omega$ the value of $f(\gamma)$ can be arbitrary by Corollary \ref{affine-embedding}.
\end{proof}

\begin{remark}
Any function $f$ in the above Lemma is automatically cuspidal. In fact, for any parabolic $Q\subset G$ with unipotent radical $N_Q$, the constant term function
\begin{equation*}
f_{N_Q}(g)=\int_{N_Q(K)\backslash N_Q(\bbA_K)}f(ng)dn.
\end{equation*}
is left invariant under $N_Q(K_\infty)$ and right equivariant under $I_\infty(1)$ against the character $\psi\phi$. A similar argument as in Lemma \ref{eigenfunction} shows that such a function must be zero.
\end{remark}

\subsection{The eigensheaf}\label{eigensheaf}
Let us reformulate the preceding observation geometrically. Denote by $\AS_\psi$ the Artin-Schreier sheaf on $\bA^1$ defined by the character $\psi$. We set $\AS_\phi:=\phi^*(\AS_\psi)$, the pull back of the Artin-Schreier sheaf to $I(1)/I(2)\cong \bA^{d}$. 

Let us denote by $\Perv(\Bun_{\cG(0,2)})^{I(1),\AS_\phi}$ the category of perverse sheaves on $\Bun_{\cG(0,2)}$ that are $(I(1),\AS_\phi)$-equivariant. We will use this notation more generally for any stack with an action of $I(1)/I(2)$.

For any $\gamma\in \Omega$ denote by $j_\gamma\colon I(1)/I(2) \incl{} \Bun_{\cG(0,2)}^\gamma$ the embedding of the big cell and by $i_\gamma\colon \Spec(k) \to \Bun^\gamma_{\cG(0,2)}$ the map given by the $\cG(0,2)$ bundle defined by $\gamma$.

The following is the geometric analog of Lemma \ref{eigenfunction}.
\begin{lemme}\label{clean}
The sheaf $\AS_\phi$ satisfies $j_{\gamma,!}\AS_\phi= j_{\gamma,*}\AS_\phi$. 

Moreover, for any $\gamma\in \Omega$ the functor 
\begin{align*}
\Perv(\bG_m) &\to \Perv(\Bun_{\cG(0,2)}^\gamma \times \bG_m) \\
   F         &\mapsto j_{\gamma,!}\AS_\phi \boxtimes F [\dim(\Bun_{\cG(0,2)})]
\end{align*}
is an equivalence of categories. An inverse is given by 
\begin{align*}
 K &\mapsto (i_\gamma\times \id_{\pline})^*K[-\dim(\Bun_{\cG(0,2)})].
\end{align*}
\end{lemme}
\begin{proof}
For any $w\in \widetilde{W}-\Omega$ we pick a representative in $\cN(\cT((t)))$, again denoted by $w$. Consider the $\cG(0,2)$-bundle $P_w$ defined by $w$. Let $U_{\alpha}\subset I(1)$ be a root subgroup corresponding to a simple affine root $\alpha$ such that $w.\alpha$ is negative, i.e. such that $U_{w.\alpha}\subset I^-(0)$. This defines an inclusion $U_\alpha \incl{} \Aut(P_w)$.  Thus we get a commutative square:
$$\xymatrix{
U_\alpha \times \Spec(k) \ar[d]^{\id,P_w} \ar[r] & \Spec(k) \ar[d]^{P_w} \\
U_\alpha \times \Bun_{\cG(0,2)} \ar[r]^-{act} & \Bun_{\cG(0,2)}
}.$$
This implies that $(\AS_{\phi}|_{U_\alpha})\boxtimes K|_{P_w} \cong \Qbar_\ell|_{U_\alpha} \boxtimes K|_{P_w}$. Since we assumed that $\AS_\phi|_{U_\alpha}$ is defined by a nontrivial character of $U_\alpha$ it follows that the stalk of $K$ at $P_w$ vanishes. Dually the same result holds for the costalk of $K$ at $P_w$.

This proves our first claim. Also for our second claim, this implies that any $(I(1),\AS_\phi)$-equivariant perverse sheaf on $\Bun_{\cG(0,2)}^\gamma$ is its $!$-extension form the substack $j_\gamma(I(1)/I(2))$. 
On this substack tensoring with the local system $\AS_{\phi}$ gives an equivalence between $(I(1)/I(2))$-equivariant sheaves and $(I(1)/I(2),\AS_\phi)$-equivariant sheaves. This proves our claim. 
\end{proof}
Using this lemma we can now define our automorphic sheaf:
\begin{definition}\label{d:eigensheaf}
We define $A_\phi\in \Perv(\Bun_{\cG(0,2)})^{I(1),\AS_\phi}$ to be the perverse sheaf given on the component $\Bun_{\cG(0,2)}^\gamma$ by $j_{\gamma,!}\AS_\phi[\dim(\Bun_{\cG(0,2)})]$. We will denote by $A^\gamma_\phi$ the restriction of $A_\phi$ to the component $\Bun_{\cG(0,2)}^\gamma$.
\end{definition}
\begin{remark}[A variant with multiplicative characters]\label{r:multchar}
We can generalize the above construction of $A_\phi$ slightly. Recall that the open cell in $\Bun_{\cG(1,2)}^\gamma$ is canonically isomorphic to $\cT^{red}_0\times I(1)/I(0)$. Any character $\chi\colon \cT_0^{red}(k) \to \Qbar_\ell$ defines a character sheaf $\Kum_\chi$ on $\cT^{red}_0$.  Then the above Lemma \ref{clean} also holds for  $(\cT^{red}_0\times I(1)/I(0), \Kum_\chi\boxtimes \AS_\phi)$-equivariant sheaves on $\Bun_{\cG(1,2)}$. We will denote the corresponding perverse sheaf on $\Bun_{\cG(1,2)}$ by $A_{\phi,\chi}$.
\end{remark}

\subsection{The geometric Hecke operators}\label{ss:Hecke}

In order to state our main result, we need to recall the definition of the geometric Hecke operators. 

The stack of Hecke modifications is the stack

$$\Hecke^{\bA^1}_{\cG(m,n)}(S):=\left\langle (\cE_1,\cE_2,x,\varphi) \bigm| \begin{array}{l} \cE_i\in \Bun_{\cG(m,n)}(S), x\colon S\to \bA^1,\\ \varphi\colon\cE_1|_{\bP^1-x\times S} \map{\cong} \cE_2|_{\bP^1-x}\end{array}\right\rangle.$$
This stack has natural forgetful maps:
\begin{equation}\label{HeckeOp}
\xymatrix{ & \Hecke^{\bA^1}_{\cG(m,n)} \ar[dl]_{\pr_1}\ar[dr]^{\pr_2}\ar[r]^-{\pr_{\bA^1}} & \bA^1\\ \Bun_{\cG(m,n)} && \Bun_{\cG(m,n)}\times \bA^1.
}\end{equation}
We will denote by  $\Hecke^{\pline}_{\cG(m,n)}:=\pr_{\bA^1}^{-1}(\pline)$.

\begin{remark}\label{RemarkHeckeDiagram}
\begin{enumerate}
\item The fiber of $\pr_1$ over the trivial bundle $\cG(m,n)\in\Bun_\cG$ is called the Beilinson-Drinfeld Gra\ss mannian. It will be denoted by $\GR_{\cG(m,n)}$. The fibers of $\GR_{\cG} \to \bA^1$ over a point $x\in\pline$ are isomorphic to the affine Gra\ss mannian $\Gr_{\cG,x}$, the quotient $\cG(K_x)/\cG({\cO}_x)$.


\item The geometric fibers of $\pr_2$ over $\Bun_{\cG(m,n)}\times\pline$ are (non-canonically) isomorphic to the affine Gra\ss mannian $\Gr_{G}$. Locally in the smooth topology on $\Bun_{\cG}\times\pline$ this fibration is trivial (e.g., Remark \ref{fibration}).

\item\label{sym} There diagram \eqref{HeckeOp} has a large group of symmetries. The group $I(0)/I(2)=(I(1)/I(2))\rtimes\cT^{red}_\infty$ acts on $\Bun_{\cG(m,2)}$ by changing the $I(2)$-level structures at $\infty$, and this action extends to the diagram \eqref{HeckeOp} (i.e., it also acts on $\Hecke^{\bbA^1}_{\cG(m,2)}$, and the maps $\pr_i$ are equivariant under these actions). The one dimensional torus $\grot$ acts on the curve $\PP^1$ fixing the points $\{0\}$ and $\{\infty\}$, hence it also acts on \eqref{HeckeOp}. Finally, the pinned automorphisms $\Aut^\dagger(G)$ act on \eqref{HeckeOp}. So we see that the group $I(0)/I(2)\rtimes(\grot\times\Aut^\dagger(G))$ acts on the diagram \eqref{HeckeOp}.

\end{enumerate}
\end{remark}

Let us first recall the Hecke-operators for constant group schemes. For this, we collect some facts about the geometric Satake equivalence (see \cite{MV}, \cite{Ginz}).

Let $\Gr_G=G((\tau))/G[[\tau]]$ be the abstract affine Gra\ss mannian, without reference to any point on $\pline$. Let $\calO=k[[\tau]]$ and let $\Aut_\calO$ be the pro-algebraic $k$-group of continuous (under the $\tau$-adic topology) automorphisms of $\calO$. Then $G[[\tau]]\rtimes\Aut_\calO$ acts on $\Gr_G$ from the left.
The $G[[\tau]]$ orbits on $\Gr_G$ are indexed by dominant cocharacter $\mu\in\xcoch(T)^+$. The orbits are denoted by $\Gr_{G,\mu}$ and their closures (the Schubert varieties) are denoted $\Gr_{G,\leq\mu}$. We denote the intersection cohomology sheaf of on $\Gr_{G,\leq\mu}$ by $\IC_\mu$. We will normalize $\IC_\mu$ to be of weight 0 (see Remark \ref{r:norm} below).

The {\em Satake category} $\Sat=\Perv_{\Aut_\calO}(G[[\tau]]\backslash G((\tau))/G[[\tau]])$, is the category of $G[[\tau]]\rtimes\Aut_\calO$-equivariant perverse sheaves (with finite-type support) on $\Gr_G$.
Similarly, we define $\Sat^\geom$ by considering the base change of the situation to $\bark$. Finally we define the {\em normalized semisimple Satake category} $\cS$ to be the full subcategory of $\Sat$ consisting of direct sums of $\IC_\mu$'s.


%

In \cite{MV} and \cite{Ginz}, it was shown that
\begin{itemize}
\item $\Sat^\geom$ carries a natural tensor structure (which is also defined for $\Sat$), such that the global cohomology functor $h=H^*(\Gr_G,-):\Sat^\geom\to\Vect$ is a fiber functor.
\item $\uAut^\otimes(h)$ is a connected reductive group over $\Ql$ which is Langlands dual to $G$. Let $\dualG:=\uAut^\otimes(h)$, then the Tannakian formalism gives the {\em geometric Satake equivalence} of tensor categories
\begin{equation*}
\Sat^\geom\cong\Rep(\dualG).
\end{equation*}
\item By construction $\dualG$ is equipped with a maximal torus $\dualT$, and a natural isomorphism $\xch(\dualT)\cong\xcoch(T)$ (in fact, $\dualG$ is equipped with a canonical pinning; see Lemma \ref{l:outSatakeH}). The geometric Satake equivalence sends $\IC_\mu$ to the irreducible representation $V_\mu$ of extremal weight $\mu$. We denote the inverse of this equivalence by $V \mapsto \IC_V$.
\end{itemize}


\begin{remark}
In \cite[\S3.5]{ArkhipovBezrukavnikov}, it was argued that $\cS$ is closed under the tensor structure on $\Sat$. Therefore $\cS$ is naturally a tensor category.

The pull-back along $\Gr_G\otimes_k\bark\to\Gr_G$ gives a tensor functor $\cS\to\Sat^\geom$, which is easily seen to be an equivalence because both categories are semisimple with explicit simple objects. Therefore, the above results in \cite{MV} and \cite{Ginz} all apply to $\cS$. In particular, we have the (semisimplified $k$-version of) the geometric Satake equivalence
\begin{equation*}
\cS\cong\Sat^\geom\cong\Rep(\dualG).
\end{equation*}
\end{remark}

\begin{remark}[Normalization of weights]\label{r:norm}
We use the normalization making this complex pure of weight $0$, i.e., we choose a square root of $q$ in $\Qbar_\ell$ and denote by $\IC_\mu$ the intersection complex, tensored by $\Qbar_\ell(\frac{1}{2}\dim(\Gr_\mu))$.

As was pointed out in \cite{FrenkelGross} it is not necessary to make this rather unnatural choice, which is made to obtain the group $\dualG$ from the category $\cS$. Alternatively, we can enlarge the category $\cS$ by including all Tate-twists of the intersection cohomology sheaves $\IC_\mu(n)$. By the previous remark this is still a neutral tensor category, defining group $\dualG_1$, which is an extension of $\dualG$ by a central, one dimensional torus.
\end{remark}

The stack $\Hecke^{\pline}_{G(m,n)}$ is a locally trivial fibration over $\Bun_{G(m,n)}\times \pline$ with fiber $\Gr_G$ and the $\G[[\tau]]$-orbits $\Gr_\mu$ on $\Gr_G$ define substacks $\Hecke^{\pline}_\mu \subset \Hecke^{\pline}_{G(m,n)}$. By abuse of notation we will also denote bu $\IC_\mu$ the intersection cohomology complex of $\Hecke^{\pline}_\mu$, shifted in degree such that $\IC_\mu$ restricts to the intersection complex on every fiber.

One defines the geometric Hecke operators as a functor (see \cite{GaitsgoryDeJong}):
\begin{eqnarray*}
\Hk:\Rep(\dualG)\times D^b(\Bun_{G(m,n)}) &\to& D^b(\Bun_{G(m,n)} \times \P1min2)\\
(V,K)&\mapsto&\Hk_{V}(K):=\pr_{2,!} \pr_1^*(K\tensor \IC_V).
\end{eqnarray*}
In order to compose these operators one extends $\Hk_V$ to an operator $D^b(\Bun_{G(m,n)} \times \pline)\to D^b(\Bun_{G(m,n)} \times \pline)$ defined as $K \mapsto \pr_{2!}((\pr_1\times\pr_{\pline})^*K\tensor\IC_V)$. 

Let $E$ be a $\dualG$-local system on $\pline$, viewed a tensor functor
\begin{equation*}
E:\cS\cong\Rep(\dualG)\to\Loc(\pline).
\end{equation*}
sending $V\in\Rep(\dualG)$ to $E_V\in\Loc(\pline)$.

A {\em Hecke eigensheaf} with eigenvalue $E$ is a perverse sheaf $K\in\Perv(\Bun_{\cG(m,n)})$ together with an isomorphisms $\Hk_V(K) \isom K \boxtimes E^V$, that are compatible with the symmetric tensor structure on $\Rep(\dualG)$ and composition of Hecke correspondences (see \cite{GaitsgoryDeJong} for a detailed exposition of the compatibility axiom).

\begin{remark}\label{rem:monodromy}
Since local systems are usually introduced differently, let us briefly recall how, the tensor functor $E$ allows to reconstruct the monodromy representation of the local system. This will be useful to fix notations for the monodromy representation.
Choose a geometric point $\generic$ over $\Spec K_0\in\pline$. The restriction to $\generic$ defines a tensor functor
\begin{equation}\label{omphi}
\omega_E:\cS\xrightarrow{E}\Loc(\pline)\xrightarrow{j^*_\generic}\Vect,
\end{equation}
i.e. (see \cite[Theorem 3.2]{DM}), a $\dualG$ torsor with an action of $\pi_1(\pline,\generic)$. Choosing a point of the torsor this defines the monodromy representation
\begin{equation}
\varphi:\pi_1(\pline,\generic)\to\dualG(\Ql).
\end{equation}
We denote by $\varphi^\geom$ the restriction of $\varphi$ to $\pi^\geom_1(\pline,\generic)=\pi_1(\geompline,\generic)$, and call it the {\em global geometric monodromy representation} of $E$.
\end{remark}


The analog of this construction for twisted groups $\cG$ will produce $\LG$-local systems. Here we define the L-group of $\cG$ to be $\LG=\dualG\rtimes \langle \sigma\rangle$, where $\sigma$ is the automorphism of order $N$ of $G$ used in the definition of $\cG$. This automorphism induces an automorphism of $\dualG$ via the geometric Satake isomorphism. In Lemma \ref{l:outSatakeH} we will check that this automorphism indeed preserves the pinning of $\dualG$. 

Let us give a definition of $\LG$-local systems, that is sufficient for our purposes. The reason why we cannot immediately apply the geometric Satake isomorphism is that the fibers of $\pr_2$ are not constant along $\P1min2$, so only some of the Hecke operators $\Hk_V$ will define global Hecke operators over $\P1min2$. However, we can pull-back the convolution diagram by the map $[N]:\pline\to \pline$. We will denote this covering by $\tpline\to \pline$. After pull-back we can as before define Hecke operators on $\Hecke^{\tpline}_{G(m,n)}$.
\begin{equation*}
\Hk\colon \Rep(\dualG)\times D^b(\Bun_{G(m,n)})\cong\cS\times D^b(\Bun_{G(m,n)}) \to D^b(\Bun_{G(m,n)} \times \tpline)
\end{equation*}
Moreover the covering group $\mu_N$ acts on the convolution diagram over $\tpline$. This defines a $\mu_N$-equivariant structure on the functor $\Hk$: on the source $\mu_N$ acts on $\cS$ via $\sigma:\mu_N\to\Aut^\dagger(G)$, which can be identified with the action of $\Aut^\dagger(\dualG)$ on $\Rep(\dualG)$ (Lemma \ref{l:outSatakeH}); on the target $\mu_N$ acts on $\tpline$. 

Let $E$ be a $\dualG$-local system on $\tpline$ together with compatible isomorphisms $\zeta^*E \cong E\times^{\dualG,\sigma(\zeta)}\dualG$ for $\zeta\in\mu_N$. We view $E$ as a tensor functor
\begin{equation*}
E:\Rep(\dualG)\to\Loc(\tpline)
\end{equation*}
together with a $\mu_N$-equivariant structure. We can define a Hecke eigensheaf $K\in \Perv(\Bun_{\cG(m,n)})$ with eigenvalue $E$ as before, but now we have to specify an isomorphism of functors $\epsilon(V):\Hk_V(K)\isom K\boxtimes E^V$ compatible with the tensor structure on $\Rep(\dualG)$ which commutes with the $\mu_N$-equivariant structures of the functor $V\mapsto \Hk_V(K)$ and the functor $E$.

Note that, if $\sigma\colon \mu_N \to \Aut(G)$ is trivial, the isomorphisms $\zeta^*E \cong E\times^{\dualG,\sigma(\zeta)} \check{G}$ define a descent datum for $E$. So in this case the definition coincides with the definition for constant groups.

With these definitions we can state our first main result:
\begin{theorem}\label{Heckeeigensheaf}
\begin{enumerate}
\item The sheaves $A=A_\phi$ and $A_{\phi,\chi}$ are Hecke eigensheaves. We will denote the eigenvalue of $A_\phi$ (resp. $A_{\phi,\chi}$)  by $\Kl_{\LG}(\phi)$ (resp. $\Kl_{\LG}(\phi,\chi)$).
\item If $\cG=\bP^1 \times G$ is a constant group scheme, the local system $\Kl_{\dualG}(\phi)$ is tamely ramified at $0$.
The monodromy action at $0$ on $\Kl_{\dualG}(\phi)$, is given by a principal unipotent element in $\dualG$.
\item For any irreducible representation $V\in \Rep(\dualG)$ the sheaf $\Kl_{\LG}(\phi,\chi)^V$ is pure.
\end{enumerate}
\end{theorem}
Since we defined $\Kl_{\dualG}(\phi)$ using geometric Hecke operators, for any point $x\in \pline$ the $\dualG$-conjugacy class of $\Frob_x$ defined by the local system is given by the Satake parameter of Gross' automorphic form $\pi(\phi)$ at $x$.


\subsection{The monodromy representation}\label{ss:statemono}
In this section, we assume $\cG=G\times\PP^1$ where $G$ is an almost simple split group over $k$. As in Remark \ref{rem:monodromy}, we will denote by $\varphi$ the monodromy representation for $\Kl_\dualG(\phi)$ and by $\varphi^\geom$ its geometric monodromy representation.

Our next result is about the geometric monodromy of $\Kl_\dualG(\phi,\chi)$ at $\infty$. Let $\Kl^\Ad_\dualG(\phi,\chi)$ be the local system on $\pline$ induced from the adjoint representation $V=\dualg$. Fixing an embedding of the local Galois group $\Gal(K^{\sep}_\infty/K_\infty)$ into $\pi_1(\pline,\generic)$, we get an action of $\Gal(K^{\sep}_\infty/K_\infty)$ on the geometric stalk of $\Kl^\Ad_\dualG(\phi,\chi)$ at the formal punctured discs $\Spec K^{\sep}_\infty$. Choosing an identifications of this stalk with $\dualg$, we get an action of $\Gal(K^{\sep}_\infty/K_\infty)$ on $\dualg$ (well-defined up to $\dualG$-conjugacy).

Let $\iinf\subset\Gal(K^{\sep}_\infty/K_\infty)$ be the inertia group. Let $\iinf^+\subset\iinf$ be the wild inertia group, and $\iinf^t=\iinf/\iinf^+$ be the tame inertia group.

\begin{theorem}\label{th:Swan-inv} Suppose $p=\textup{char}(k)$ is good for $G$ if $G$ is not simply-laced (i.e., $p>2$ when $G$ is of type $B_n,C_n$, and $p>3$ when $G$ is of type $F_4,G_2$). Then
\begin{itemize}
\item $\Swan_{\infty}(\Kl^\Ad_\dualG(\phi,\chi))=r(\dualG)$, the rank of $\dualG$;
\item $\dualg^{\iinf}=0$.
\end{itemize}
\end{theorem}

The Swan equality will be proved in Corollary \ref{c:Swan}; the vanishing of $\iinf$-invariants will be proved in Prop. \ref{p:van}(2). 


\begin{corollary}\textup{(Gross-Reeder \cite[Proposition 5.6]{GrossReeder})}\label{c:simple-wild}
Suppose $p=\textup{char}(k)$ does not divide $\#W$, then the local geometric Galois representation $\varphi^\geom_\infty:\iinf\to\dualG$ is a {\em simple wild parameter} defined in \cite[\S 6]{GrossReeder}. More precisely, up to $\dualG$-conjugation, we have a commutative diagram of exact sequences
\begin{equation*}
\xymatrix{\iinf^+\ar[r]\ar[d]^{\varphi^{\geom,+}_\infty} &\iinf\ar[r]\ar[d]^{\varphi^\geom_\infty} & \iinf^t\ar[d]\\
\dualT\ar[r] & N(\dualT)\ar[r] & W}
\end{equation*}
where
\begin{itemize}
\item A topological generator of the tame inertia $\iinf^t$ maps to a Coxeter element $\Cox\in W$ (well-defined up to $W$-conjugacy, see \cite[Ch.V,\S6]{Bourbaki});
\item The wild inertia $\iinf^+$ maps onto a subgroup $\dualT(\zeta)\subset\dualT[p]$. Here $\zeta\in\overline{\FF}^\times_p$ is a primitive $h$-th root of unity ($h$ is the Coxeter number), and $\dualT(\zeta)\subset\dualT[p]$ is the unique $\FF_p[\Cox]$-submodule of $\dualT[p]$ isomorphic to $\FF_p[\zeta]$ (on which $\Cox$ acts by multiplication by $\zeta$).
\item The nonzero breaks of the $\iinf$-representation $\Kl^\Ad_\dualG(\phi,\chi)$ are equal to $1/h$.
\end{itemize}
\end{corollary}


Finally we can state our result on the global geometric monodromy of $\Kl_\dualG(\phi)$. Let $\dualG_\geom\subset\dualG$ be the Zariski closure of the image of the global geometric monodromy representation $\varphi^\geom$. 

\begin{theorem}\label{th:glob} Suppose $\textup{char}(k)>2$, then the geometric monodromy group $\dualG_\geom$ for $\Kl_\dualG(\phi)$ is connected, and
\begin{itemize}
\item $\dualG_\geom=\dualG^{\Aut^\dagger(\dualG),\circ}$, if $\dualG$ is not of type $A_{2n}$ ($n\geq2$) or $B_3$;
\item $\dualG_\geom=\dualG$ if $\dualG$ is of type $A_{2n}$;
\item $\dualG_\geom=G_2$ if $\dualG$ is of type $B_3$ and $\textup{char}(k)>3$.
\end{itemize}
\end{theorem}
The proof will be given in \S\ref{ss:glob}.

\subsection{Variant}
There is a variant of our construction using $D$-modules instead of $\ell$-adic sheaves. The base field is then taken to be $k=\CC$. The Artin-Schreier local system $\AS_\psi$ is replaced by the exponential $D$-module $\CC\jiao{x,\partial_x}/(\partial_x-1)$ on $\bbA^1_{\CC}=\Spec\CC[x]$, and all the rest of the construction carries through, and we get a tensor functor:
\begin{equation*}
\Kl_\LG(\phi)_{\dR}: \Rep(\LG_{\CC})\to\Conn(\plineC),
\end{equation*}
where $\Conn(\plineC)$ is the tensor category of vector bundles with connections on $\plineC$.

We conjecture that our construction should give the same connection as the Frenkel-Gross construction. To state this  precisely, let $G$ be almost simple, and $\cG$ be a quasi-split form of $G$ over $\pline$ given by $\sigma:\mu_N\to\Aut^\dagger(G)$. Recall from \cite[\S 5]{FrenkelGross} that the Frenkel-Gross connection on the trivial $\LG$-bundle on $\plineC$:
\begin{equation*}
\nabla_{\LG}(\dualX_{0},\cdots,\dualX_{r_\sigma})=d+\sum_{i=0}^{r_\sigma}\dualX_{i}\frac{d}{dz}
\end{equation*}
where $r_\sigma$ is the rank of $G^\sigma$, and $\dualX_{i}$ is a basis of the $(-\alpha_i)$-root space of the (twisted) affine Kac-Moody Lie algebra associate to $\dualg$ and $\sigma$.

\begin{conjecture}\label{conj:FG}
There is a bijection between the set of generic linear functions $\phi: I(1)/I(2)\to\GG_{a,\CC}$ and the set of bases $(\dualX_0,\cdots,\dualX_{r_\sigma})$, such that whenever $\phi$ corresponds to $(\dualX_0,\cdots,\dualX_{r_\sigma})$ under this bijection, there is a natural isomorphism between $\LG$-connections on $\plineC$: $$\Kl_\LG(\phi)_{\dR}\cong(\LG,\nabla_\LG(\dualX_0,\cdots,\dualX_{r_\sigma})).$$
\end{conjecture}




\section{Example: Kloosterman sheaf for $\GL_{n}$}\label{s:exam}


In this section, we calculate the Kloosterman sheaf $\Kl_{\GL_n}(\phi,\chi)$ for the constant group $G=\GL_n$ over $\PP^1$. Its Langlands dual is $\GL_{n,\Ql}$ and we will denote the standard representation by $\St$. To describe this $\GL_n$-local system over $\pline$ is then the same as describing the rank $n$ local system $\Kl^\St_{\GL_n}(\phi,\chi)$. We will see that this rank $n$ local system coincides with the classical Kloosterman sheaf defined by Deligne in \cite{DeligneSommes}.

\subsection{Another modular interpretation}
We want to interpret $\GL_n(1,2)$-bundles in terms of vector bundles. 
We first define a variant of $\Bun_{\GL_n(1,2)}$. Let $\Bun_{n,1,2}$ be the stack classifying the data $(\calE,F^*\calE,\{v^i\},F_*\calE,\{v_i\})$ where
\begin{enumerate}
\item $\calE$ is a vector bundle of rank $n$ on $\PP^1$;
\item a decreasing filtration $F^*\calE$ giving a complete flag of the fiber of $\calE$ at $0$:
\begin{equation*}
\calE=F^0\calE\supset F^1\calE\supset\cdots\supset F^n\calE=\calE(-\{0\});
\end{equation*}
\item a nonzero vector $v^i\in F^{i-1}\calE/F^{i}\calE$ for each $i=1,\cdots,n$;
\item an increasing filtration $F_*\calE$ giving a complete flag of the fiber of $\calE$ at $\infty$: 
\begin{equation*}
\calE(-\{\infty\})=F_0\calE\subset F_1\calE\subset\cdots\subset F_n\calE=\calE;
\end{equation*}
\item a vector $v_i\in F_i\calE/F_{i-2}\calE$ which does not lie in $F_{i-1}\calE/F_{i-2}\calE$, for $i=1,\cdots,n$ (we understand $F_{-1}\calE$ as $(F_{n-1}\calE)(-\{\infty\})$). 
\end{enumerate}

Note that $\Bun_{n,1,2}$ is the moduli stack of $\cG$-torsors over $\PP^1$, where $\cG$ is the Bruhat-Tits group scheme over $\PP^1$ such that
\begin{itemize}
\item $\cG|_{\pline}=\GL_n\times(\pline)$;
\item $\cG(\calO_0)=I_{\GL_n}(1)^{\opp}$;
\item $\cG(\calO_\infty)=Z_{\GL_n}(1)(k[[s]])\cdot I_{\SL_n}(2)\supset I_{\GL_n}(2)$.
\end{itemize}
The only difference between $\cG$ and $\GL_n(1,2)$ is that they take different level structures for the center $\Gm=Z_{\GL_n}$ at $\infty$. Therefore we have a natural morphism $\GL_n(1,2)\to\cG$, hence a natural morphism $\Bun_{\GL_n(1,2)}\to\Bun_{n,1,2}$, which is a $\GG_a$-torsor. 

Choosing a trivialization of the bundle $\cE$ over $\pline$, we can rewrite the moduli problem for $\Bun_{n,1,2}$ in the following way. Let $\Lambda$ be the free $k[t,t^{-1}]$-module with basis $\{e_1,\cdots,e_n\}$. Let $e_{i+jn}=t^je_i$ for $1\leq i\leq n$ and $j\in\ZZ$, then $\Lambda$ is a $k$-vector space with basis $\{e_i\}_{i\in\ZZ}$. For any $k$-algebra $R$, an {\em $R[t]$-lattice} in $R\otimes\Lambda$ is an $R[t]$-submodule $\Lambda'\subset R\otimes\Lambda$ such that there exists $M\in\ZZ_{>0}$ such that:
\begin{equation*}
\Span_R\{e_i|i>M\}\subset\Lambda'\subset\Span_R\{e_i|i\geq-M\}
\end{equation*}
and that both $\Lambda'/\Span_R\{e_i|i>M\}$ and $\Span_R\{e_i|i\geq-M\}/\Lambda'$ are projective $R$-modules. Similarly we can define the notion of $R[t^{-1}]$-lattices in $R\otimes\Lambda$.

Let $\tBun_{n,1,2}(R)$ classify the data $(\Lambda^*,\{v^i\}_{0\leq i\leq n-1},\Lambda_*,\{v_i\}_{1\leq i\leq n})$:
\begin{enumerate}
\item $R\otimes\Lambda\supset\Lambda^0\supset\Lambda^1\supset\cdots\supset\Lambda^n=t\Lambda^0$ is a chain of $R[t]$-lattice such that $\Lambda^i/\Lambda^{i+1}$ is a projective $R$-module of rank 1. We let $\Lambda^{i+jn}=t^j\Lambda^i$ for any $1\leq i\leq n$ and $j\in\ZZ$.
\item $v^i\in\Lambda^i/\Lambda^{i+1}$ is an $R$-basis for $i=0,\cdots,n-1$;
\item $R\otimes\Lambda\supset\Lambda_n\supset\Lambda_{n-1}\supset\cdots\supset\Lambda_0=t^{-1}\Lambda_n$ is a chain of $R[t^{-1}]$-lattices such that $\Lambda_i/\Lambda_{i+1}$ is a projective $R$-module of rank 1. We let $\Lambda_{i+jn}=t^j\Lambda_i$ for any $1\leq i\leq n$ and $j\in\ZZ$.
\item $v_i\in\Lambda_i/\Lambda_{i-2}$ whose image in $\Lambda_i/\Lambda_{i-1}$ is an $R$-basis, for $i=1,\cdots,n$.
\end{enumerate}
The group $\GL_n(k[t,t^{-1}])$ (viewed as an ind-group scheme over $k$) acts on $\Lambda$, and hence it acts on the stack $\tBun_{n,1,2}$. Moreover, $\Bun_{n,1,2}$ is naturally isomorphic to the quotient stack $\tBun_{n,1,2}/\GL_n(k[t,t^{-1}])$.

Associated with the chains of lattices $(\Lambda^*,\Lambda_*)$ is the locally constant integer-valued function on $\Spec R$:
\begin{equation*}
\deg(\Lambda^*,\Lambda_*):=\chi_R(\Lambda^0\oplus\Lambda_0\xrightarrow{(\iota^0,\iota_0)}R\otimes\Lambda)
\end{equation*}
Here $\iota^0,\iota_0$ are inclusion maps and the Euler characteristic on the RHS is defined as $\rk_R\ker(\iota^0,\iota_0)-\rk_R\coker_R(\iota^0,\iota_0)$.

For $d\in\ZZ$, let $\Bun^d_{n,1,2}\subset\Bun_{n,1,2}$ be the substack classifying $\{\Lambda^*,v^i,\Lambda_*,v_i\}$ with $\deg(\Lambda^*,\Lambda_*)=d$. 

The embedding of the big cell $j^d:T\times\GG_a^n\hookrightarrow\Bun^d_{n,1,2}$ can be fixed as follows. For $\una=(a_1,\cdots,a_n)\in T(R)=(R^\times)^n$ and $\unb=(b_1,\cdots,b_n)\in R^n$, the point $j^d(\una,\unb)\in\Bun^0_{n,1,2}(R)$ is given by the $\GL_n(R[t,t^{-1}])$-orbit of the following data: 
\begin{enumerate}
\item $\Lambda^i=\Span_R\{e_j|j>i\}\subset R\otimes\Lambda$;
\item $v^i=a_{i+1}e_{i+1}\in\Lambda^i/\Lambda^{i+1}$;
\item $\Lambda_i=\Span_R\{e_j|j\leq i+d\}\subset R\otimes\Lambda$;
\item $v_i=e_{i+d}+b_ie_{i+d-1}\in\Lambda_{i}/\Lambda_{i-2}$.
\end{enumerate}
In particular, we have a base point $\star_d=j^d(\underline{1},\underline{0})\in\Bun^d_{n,1,2}$. We denote the $k[t]$-chain (resp. $k[t^{-1}]$-chain) of $\star_d$ by $\Lambda^*(\star_d)$ (resp. $\Lambda_*(\star_d)$). Notice that the underlying vector bundle for any point $j^d(\una,\unb)$ is the bundle
\begin{equation*}
\calE_d=\calO(m+1)^{\oplus r}\oplus\calO(m)^{\oplus n-r}
\end{equation*}
where $m\in\ZZ$ and $0\leq r\leq n-1$ is uniquely determined by $d=mn+r$.

\subsection{The Kloosterman sheaf associated with the standard representation of $\GL_n$}
Let $\chi:T(k)\to\Qlt$ be a character, which defines a Kummer local system $\Kum_\chi$ on $T$. The perverse sheaves
\begin{equation*}
A^d_{\phi,\chi}=j^d_!(\Kum_\chi[n]\boxtimes\phi^*\AS_\psi[n])
\end{equation*}
on $\Bun^d_{n,1,2}$ form a Hecke eigensheaf on $\Bun_{n,1,2}$ by Theorem \ref{Heckeeigensheaf}.


Consider the Hecke correspondence $\Hecke_{\omega_1}$ given by the coweight $\omega_1=(1,0,\cdots,0)$ ($\omega_1$ defines the representation $\St$). We restrict it to the components of $\Bun^{0}$ and $\Bun^{1}$ and the curve $\pline$:
\begin{equation}\label{HkGL}
\xymatrix{& \Hecke_{\omega_1}\ar[dl]_{pr_1}\ar[dr]^{pr_2} &\\ \Bun^{0}_{n,1,2} & & \Bun^{1}_{n,1,2}\times(\pline)}
\end{equation}
We would like to evaluate the eigenvalue local system $\Kl^\St_{\GL_n}(\phi,\chi)$ on $\pline$ characterized by:
\begin{equation*}
\pr_{2,!}\pr^*_1A^0_{\phi,\chi}[n-1]\cong A^1_{\phi,\chi}\boxtimes\Kl^\St_{\GL_n}(\phi,\chi).
\end{equation*}
The intersection cohomology sheaf $\IC_{\omega_1}$ on $\Hk_{\omega_1}$ is simply $\Ql[n-1]$ because $\pr_1:\Hk_{\omega_1}\to\Bun_{n,1,2}$ is a $\PP^{n-1}$-bundle \footnote{In this calculation, we will {\em not} normalized $\IC_{\omega_1}$ to be of weight 0, as we did in Remark \ref{r:norm}. So strictly speaking, the sheaf $\Kl^\St_{\GL_n}(\phi,\chi)$ differs from the one defined in \S\ref{ss:Hecke} by a Tate twist $(\frac{n-1}{2})$. }. We restrict the diagram \eqref{HkGL} to the fiber of $\star_1\times(\pline)\subset\Bun^1_{n,1,2}\times(\pline)$ under $\pr_2$, and the fiber $\GR_{\omega_1}^{\circ}$ over the big cell in $\Bun^0_{n,1,2}$ under $\pr_1$, we get
\begin{equation}\label{diaU}
\xymatrix{& \GR_{\omega_1}^{\circ}\ar[r]\ar[dl]_{\pr_1} & \GR_{\omega_1}\ar[dl]_{\pr_1}\ar[dr]^{\pi} &\\ T\times\GG_a^n\ar[r]^{j^0} & \Bun^0_{n,1,2} & & \pline}
\end{equation}
The $R$-points of the open subscheme $\GR_{\omega_1}^{\circ}\subset\GR_{\omega_1}$ classifies $R[t,t^{-1}]$-homomorphisms $M:R\otimes\Lambda\to R\otimes\Lambda$  sending the chains $(\Lambda^*(\star_0),\Lambda_*(\star_0))$ to $(\Lambda^*(\star_1),\Lambda_*(\star_1))$ up to pre-composing with automorphisms of $(R\otimes\Lambda,\Lambda^*(\star_0),\Lambda_*(\star_0))$. With respect to the $R[t,t^{-1}]$-basis $\{e_1,\cdots,e_n\}$ of $R\otimes\Lambda$, any such $M$ takes the form
\begin{equation*}
\left( \begin{array}{ccccc}
a_1 &     &   & &tb_n \\
b_1 & a_2     \\
    & b_2 & a_3 \\
& & \ddots & \ddots\\
& & & b_{n-1}& a_n \end{array} \right)
\end{equation*}
up to right multiplication by diagonal matrices (here $a_i,b_i\in R^\times$). Therefore, we can normalize the matrix of $M$ to be
\begin{equation}\label{mtau}
M=\left( \begin{array}{cccc}
a_1 &     &    &t \\
1 & a_2     \\
    & \ddots & \ddots \\
& &  1 & a_n \end{array} \right)
\end{equation}
In other words, we get an identification $U\cong\GG^n_m$ by sending $M$ to $(-a_1,\cdots,-a_n)$.

We now describe the projections $\pr_1$ and $\pi$ using this identification. Using the isomorphism $M$ we can pull back the vectors $\{v^i=e_{i+1}\},\{v_i=e_{i+1}\}$ in the data of $\star_1$ to get the corresponding vectors for the point $\pr_1(M)=(\Lambda,\Lambda^*(\star_0),\{v^i\},\Lambda_*(\star_0),\{v_i\})\in\Bun^0_{n,1,2}$:
\begin{eqnarray*}
v^i=a_{i+1}^{-1}e_{i+1} & i=0,\cdots,n-1;\\
v_i=e_i-a_ie_{i-1} & i=1,\cdots,n.
\end{eqnarray*}
This means that $\pr_1(M)$ for $M$ as in \eqref{mtau} has coordinates:
\begin{equation*}
\pr_1(M)=j^0(a_1^{-1},\cdots,a_n^{-1},-a_1,\cdots,-a_n).
\end{equation*}
On the other hand, the point $\pi(M)\in\pline$ is the value of $t$ such that $\det(M)=0$, i.e.,
\begin{equation*}
\pi(M)=(-1)^na_1\cdots a_n.
\end{equation*}
In summary, the maps $\pr_1$ and $\pi$ from $\GR_{\omega_1}^{\circ}$ can be identified with the following maps
\begin{equation*}
\xymatrix{& \Gm^n\ar[dl]_{(-\inv,\iota)}\ar[dr]^{\mult}\\ \Gm^n\times\GG_a^n & & \Gm=\pline}
\end{equation*}
where $\inv:\Gm^n\to\Gm^n$ is the coordinate-wise inverse and $\iota:\Gm^n\hookrightarrow\GG_a^n$ is the natural inclusion.

We thus get
\begin{proposition}
The Kloosterman sheaf $\Kl^{\St}_{\GL_n}(\phi,\chi)$ associated with the Hecke eigensheaf $A_{\phi,\chi}$ and the standard representation of the dual group $\GL_{n,\Ql}=\widehat{\GL}_n$ takes the form
\begin{equation*}
\Kl^\St_{\GL_n}(\phi,\chi)\cong\mult_!((-\inv)^*\Kum_\chi\otimes\phi^*\AS_\psi)[n-1]
\end{equation*}
Here we denote the restriction of $\phi:\GG_a^n\to\GG_a$ to $\Gm^n$ still by $\phi$.
\end{proposition}
\begin{remark}
When $\phi=\add:\GG_a^n\to\GG_a$ is the addition of the coordinates, and $\chi$ is written as $n$-multiplicative characters $(\chi_1,\cdots,\chi_n)$, the Kloosterman sheaf $\Kl^\St_{\GL_n}(\phi,\chi)$ is the same as the Kloosterman sheaf $\Kl(\psi;\chi_1,\cdots,\chi_n;1,\cdots,1)$ defined by Katz in \cite[\S 4.1]{KatzKloosterman}, which is a generalization of Deligne's Kloosterman sheaves \cite[\S 7]{DeligneSommes}.
\end{remark}


\section{Proof of Theorem \ref{Heckeeigensheaf}}\label{ProofThm}



In this section we prove Theorem \ref{Heckeeigensheaf}.

\subsection{\texorpdfstring{First Step: $\Hk_V(A)$ is perverse.}{First Step: Hk(A) is perverse.}}
We want to show that for every $V\in \Rep(\check{G})$ the complex $\Hk_V(A_{\phi,\chi})[-1]$ is a perverse sheaf.
In order to simplify notations, we will only consider the sheaf $A=A_\phi$. The proof for $A_{\phi,\chi}$ is identical, one only needs to replace $\cG(0,2)$ by $\cG(1,2)$ everywhere in the argument.

Let us recall the convolution diagram from section \ref{eigensheaf}: 
$$\xymatrix{ & \Hecke^{\bA^1}_{\cG(0,2)} \ar[dl]_{\pr_1}\ar[dr]^{\pr_2}\ar[r]^-{\pr_{\bA^1}} & \bA^1\\ \Bun_{\cG(0,2)} && \Bun_{\cG(0,2)}\times \bA^1
}$$
Our proof is based on a few simple geometric observations. First we need to recall that the maps $\pr_i$ in the above diagram are locally trivial fibrations:
\begin{remark}\label{fibration}
The map $\pr_1$ is a locally trivial fibration, i.e., there exists a smooth atlas $U\to \Bun_{\cG(0,2)}$, such that $$U\times_{\Bun_{\cG(0,2)}} \Hecke_{\cG(0,2)}^{\bA^1} \cong U \times \GR_{\cG(0,2)}.$$ 

Furthermore, the map $\pr_2$ is also locally trivial on this atlas, i.e.,
$$(U\times \bA^1)\times_{\Bun_{\cG(0,2)}\times \bA^1} \Hecke_{\cG(0,2)}^{\bA^1}  \cong U \times \GR_{\cG(0,2)}.$$
\end{remark}
\begin{proof}
In order to find an atlas $p\colon U\to \Bun_{\cG(0,2)}$ satisfying these conditions, we only need that the family of $\cG(0,2)$-bundles corresponding to $p$ on $U\times \bP^1$ is trivial over $U\times \bA^1$. By Proposition \ref{P1-uniformization} for any $U$ this condition is satisfied locally in the \'etale topology on $U$.
\end{proof}
For $\gamma\in \Omega$ we denoted by $j_\gamma\colon I(1)/I(2) \incl{} \Bun_{\cG(0,2)}$ the canonical embedding. Recall that $A^\gamma=j_{\gamma,!}\AS_\phi[\dim(\Bun_{\cG(0,2)})]$.  Denote by $j^\prime_\gamma\colon \pr_1^{-1}(I(1)/I(2)) \incl{} \Hecke_{\bA^1}^{\cG(0,2)}$ the inverse image of this embedding into the Hecke stack.
\begin{remark}\label{affine}
The restriction of $\pr_2$ to $\pr_1^{-1}(I(1)/I(2))$:  
$$\pr_2\colon \pr_1^{-1}(I(1)/I(2)) \to \Bun_{\cG(0,2)} \times \bA^1$$
is affine.
\end{remark}
\begin{proof}
By Corollary \ref{affine-embedding} the open subset $I(1)/I(2) \hookrightarrow \Bun_{\cG(0,2)}$ is defined by the non-vanishing of sections of line bundles $\cL_i$, which can be defined as the pull-back of the corresponding sections for $\GL_n$ under a faithful representation $\cG\to \GL_n$. If $\cG=\GL_n$ then for any $x\in\bA^1$ the pull-backs of these bundles generate of the Picard-group $\Gr_{\GL_n(0,2),x}$, which is an ind-projective scheme (e.g., \cite[Theorem  7 and 8]{Faltings}).  Thus the preimage of $I(1)/I(2)$ in $\GR_{\GL_n(0,2)}$ is affine over $\bA^1$.

For general $\cG$ the ind-scheme $\Gr_{\cG(0,2),x}$ is usually constructed as a closed sub-scheme of $\Gr_{\GL_n}$ for a suitable faithful representation (\cite{PappasRapoport}). So again, the claim follows from the case $\cG=\GL_n$. \end{proof}

Now we can prove the first step. Fix any $V\in \Rep(\check{G})$ and the corresponding perverse sheaf  $\IC_V$ on $\Hecke^{\bA^1}_{\cG(0,2)}$. We claim that Remark \ref{fibration} implies: 
$$j^\prime_!(\pr_1^*\AS_\phi \tensor \IC_V) \cong j^\prime_*(\pr_1^*\AS_\phi \tensor \IC_V).$$

To see this, recall that the formation of $j_!$ and $j_*$ commutes with smooth base-change. Thus, to prove the claimed isomorphism, we may choose a smooth atlas $p\colon U\to \Bun_{\cG(0,2)}$ such that the pull-back of $\pr_1$ is isomorphic to the projection $U\times \GR_{\cG} \to U$. But in this case the claim follows from the isomorphism $j_!\AS_\phi = j_* \AS_\phi$ (Lemma \ref{clean}).

In particular, we find that $j^\prime_!(\pr_1^*A \tensor \IC_V)[1]$ is a perverse sheaf and
\begin{align*}
\Hk_V(A) & = (\pr_2 \circ j^\prime)_! (\pr_1^*(A)\tensor \IC_V)\\
&=  (\pr_2 \circ j^\prime)_* (\pr_1^*(A)\tensor \IC_V).
\end{align*}
By Remark \ref{affine}, the map $(\pr_2 \circ j^\prime)$ is affine. Therefore $(\pr_2 \circ j^\prime)_*$ is right-exact for the perverse t-structure (\cite[\S 4.1.1]{BeilinsonBernsteinDeligne}). Thus $\Hk_V(A)[1]$ must be perverse. 

\subsection{Second step: $A$ is an eigensheaf.}

We already noted (Remark \ref{RemarkHeckeDiagram}) that the action of $I(1)/I(2)$ on $\Bun_{\cG(0,2)}$ extends to an action on the convolution diagram. Thus $\Hk_V(A)$ is $(I(1)/I(2),\phi)$-equivariant. By Lemma \ref{clean} this implies that for any $\gamma\in \pi_0(\Bun_{\cG(0,2)})$ we have $$\Hk_V(A)|_{\Bun_{\cG(0,2)}^\gamma\times X}= A^\gamma \boxtimes E^\gamma_V,$$ where $E^\gamma_V[1]$ is a perverse sheaf on $X$.

We claim that the sheaf $E^\gamma_{V}$ does not depend on $\gamma$. To show this, we may assume that $\IC_V$ is supported only on one connected component $\GR_\cG^\gamma \subset \GR_{\cG}$, because the functor $\Hk_V$ is isomorphic to the direct sum of the functors defined by the restriction of $\IC_V$ to the connected components. 

Note that for any $\gamma^\prime\in \Omega$ the Hecke operator $\Hk_{\gamma^\prime}$ commutes with $\Hk_{V}$, because Hecke operators supported at different points of $\bP^1$ commute. By definition $A$ is an eigensheaf for the operators $\Hk_{\gamma^\prime}$. Thus 
\begin{align*} 
A^{0} \boxtimes E_V^{0} &= \Hk_V(A^{-\gamma}) = \Hk_V(\Hk_{-\gamma}(A^0)) \\
&= \Hk_{-\gamma}(\Hk_V(A^0)) =\Hk_{-\gamma}(A^{\gamma}\boxtimes E_V^{\gamma}) \\
& = A^{0} \boxtimes E_{V}^{\gamma}. 
\end{align*}
This implies that the sheaves $E_V^\gamma$ are canonically isomorphic to $E_V^0$, so that we may drop the index $\gamma$.

The $\Hk_V$ are compatible with the tensor product of representations, 
 so in particular we have 
$$ \Hk_{V\tensor V} (A) = \Hk_{V} \circ \Hk_V (A) = A \boxtimes (E_V\tensor E_V).$$

This implies that $E_V\tensor E_V[-1]$ is again a perverse sheaf. Therefore, $E_V$ must be a perverse sheaf concentrated in cohomological degree $0$. So for any $V$ the complex $E_V$ is a sheaf, such that $E_V=j_{\eta,*} E_{V,\eta}$, where we denoted by $j_\eta\colon\eta \hookrightarrow \bP^1$ the inclusion of the generic point.

Also, we obtain a tensor functor $V \mapsto (E_V)_\eta$ with values in the category of local systems on $\eta$. This is a rigid tensor category, so by \cite[Theorem  3.2]{DM} this defines a $\dualG$-local system over $\eta$. We need to show that $E_{V,\eta}$ extends to a local system on $\pline$.

For the trivial Hecke operator $\Hk_1$ we have canonical isomorphisms $E_1=\Qbar_\ell$ together with maps $E_{1} \to E_V \tensor E_{V^*} \to E_{1}$ such that the composition is equal to multiplication by $\dim(V)$. We already know that $E_V=j_{\eta,*}E_{V,\eta}$. Thus for any geometric point $\overline{x}$ of $\bG_m$ the fiber $E_{V,\overline{x}}$ is a subspace of the geometric generic fiber $E_{V,\overline{\eta}}$ and the canonical map $\id_{\overline{\eta}}:\Qbar_{\ell,\overline{\eta}} \to (E_V \tensor E_{V^*})_{\overline{\eta}}$ factors through $E_{V,\overline{x}}\tensor E_{V^*,\overline{x}}$. This implies that the sheaves $E_V$ are locally constant, because $\id_{\overline{\eta}}$ corresponds to the identity of $E_{V,\overline{\eta}}$. 

Thus the $E_V$ define a tensor functor from $\Rep(\dualG)$ to the category $\textrm{Loc}(\pline)$ of local systems on $\pline$. Since again this is a rigid tensor category this defines a $\dualG$-local system on $\pline$.

\subsection{Third step: The monodromy at the tame point.}\label{ss:unipmono}

As in the statement of Theorem \ref{Heckeeigensheaf} (2) we now assume $\cG=G\times \bP^1$ is a constant split group. To compute the monodromy of $\Kl_{\dualG}(\phi)$ at $\{0\}$ we rephrase an argument of Bezrukavnikov \cite{Bezrukavnikov}. His argument relies on results on central sheaves of Gaitsgory \cite{Gaitsgory} and on Gabber's result that the monodromy filtration on nearby cycles coincides with the weight filtration \cite{BeilinsonBernstein}. 

In order to explain the argument we need to recall Gaitsgory's construction (\cite{Gaitsgory}). He considered the diagram 
$$\xymatrix{
\GR^{\pline}_{\cG}\ar@{^(->}[r]^-{j_{GR}}\ar[d] & \GR_{\cG}\ar[d] & \GR_{\cG,0}=\Fl_G\ar@{_(->}[l]\ar[d]\\
\pline \ar@{^(->}[r]^-{j} & \bA^1 & \{0\}\ar@{_(->}[l]
}$$ 
and the induced nearby cycles functor $\Psi: \Perv(\GR^{\pline}_{\cG}) \to \Perv_{I_0}(\Fl_G)$. 
He showed that the monodromy action on the sheaves $\Psi(\IC_V)$ is unipotent.

Since the map $\pr_1: \Hecke^{\bA^1}_{\cG(0,2)}\to \Bun_{\cG(0,2)}$ is locally isomorphic to the product with fibers isomorphic to $\GR_{\cG(0,2)}$, we know that $\Psi(\pr_1^*(A)\tensor \IC_V)=\pr_1^*(A)\tensor\Psi(\IC_V)$. In particular, by Gaitsgory's result \cite[Theorem 2]{Gaitsgory} we find that the monodromy action on this sheaf is unipotent. Therefore (\cite[Lemma 5.6]{GoertzHaines}) the monodromy action on $\pr_{2,!}(\Psi(\pr_1^*(A)\tensor \IC_V))$ is again unipotent. Since taking nearby-cycles commutes with proper push-forward, the monodromy action on 
\begin{equation}\label{prPsi}
\pr_{2,!}(\pr_1^*(A)\tensor\Psi(\IC_V))=\pr_{2,!}(\Psi(\pr_1^*(A)\tensor \IC_V))=\Psi(A\boxtimes E_V)=A\boxtimes\Psi(E_{V})
\end{equation}
is also unipotent. By definition, $\Psi(E_V)$ is the stalk of $E_V$ at the geometric point $\Spec K^{\sep}_0$ over the punctured formal neighborhood of $0$, carrying the $\Gal(K^{\sep}_0/K_0)$-action as monodromy. Restricting \eqref{prPsi} to the trivial $\cG(0,2)$-bundle, we therefore get
\begin{equation}\label{conv}
\Kl^V_{\dualG}|_{\Spec K^\sep_0}=R\Gamma_c(\pr^*_1A\otimes\Psi(\IC_V))
\end{equation}
with the inertia group $\inert_0\subset\Gal(K^\sep_0/K_0)$ acting tamely and unipotently.

We have to show that the monodromy action is given by a principal unipotent element.

Recall that the $\cG(0,2)(k[[t]])=I(0)^\opp$-orbits on the affine flag manifold $\Fl_G=\GR_{\cG,0}$ are parametrized by the Iwahori-Weyl group $\widetilde{W}$ (\cite[Prop.8.1]{PappasRapoport}). The intersection cohomology sheaves of the closures of these orbits will be denoted by $\IC_{\tilw}$. The convolution with these sheaves defines Hecke operators $\Hk_{\tilw}(A):=\pr_{2,!}(\pr_1^*(A)\tensor \IC_{\tilw})$.

By a result of G\"ortz and Haines \cite[Corollary 1.2]{GoertzHaines}, the sheaf $\Psi(\IC_V)$ has a filtration such that the associated graded sheaves are isomorphic to $\IC_w(i)$ for some $w\in \widetilde{W},i\in\bZ$ and the multiplicity of the $\IC_1(i)$ is equal to $\dim H^{2i}(\Gr_G, \IC_V)$. 

Now, for any $\tilw\in \widetilde{W}$ of length $l(\tilw)>0$ there exists a simple reflection $s$ such that $\tilw=\tilw^\prime s$ and $l(\tilw^\prime)<l(\tilw)$. Write $P_s:= I(0)^\opp\cup I(0)^\opp sI(0)^\opp$ for the parahoric subgroup generated by $I(0)^\opp$ and $s$, so that the projection $\pr_s:\Fl \to G((t))/P_s$ is a $\bP^1$-bundle. From this we see that  the sheaf $\IC_\tilw$ is of the form $\pr_s^*(IC_{\tilw^\prime})$. 
Therefore the complex $\Hk_{\tilw}(A)$ 
is invariant under the larger parahoric subgroup $P_s\supset I(0)^\opp$. Since $\Hk_{\tilw}(A)$ is also $(I(1)/I(2),\AS_\phi)$-equivariant, this implies that $\Hk_{\tilw}(A)=0$.

Thus only the Hecke operators of length $0$ act non-trivially on $A$. Therefore, by the result of G\"ortz and Haines the dimensions of the weight filtration of $\Psi(\Kl^V_\dualG)=\Kl^V_{\dualG,\generic}$ are given by the dimensions of $H^{2i}(\Gr_G,\IC_V)$. Since the monodromy filtration agrees with the weight filtration \cite{BeilinsonBernstein} the monodromy must act as a principal nilpotent element. This proves our claim. 

\begin{remark} A key ingredient needed in the article of G\"ortz and Haines is the formula for the trace function of the sheaves $\Psi(\IC_\mu)$. This trace is given by Bernstein's formula for the central elements of the Iwahori-Hecke algebra (see \cite[\S 2.7]{GoertzHaines}). One can also deduce our proposition directly from this formula.
\end{remark} 




\section{Cohomological properties of Kloosterman sheaves}\label{s:coho}


In this section, let $G$ be a split almost simple group over $k$, viewed as a constant group scheme over $\PP^1$. Recall that for each root $\alpha\in\Phi$, $U_\alpha$ is the root group in $G$. We fix an isomorphism $u_\alpha:\GG_a\isom U_\alpha$ for each root $\alpha$.

\noindent{\bf Notation.} For a scheme $X$ defined over $k$ and a $\Ql$-complex $F$ of sheaves on $X$, $\chi_c(X,F)$ and $H^*(X,F)$ mean the Euler characteristic and the cohomology of the {\em pull back} of $F$ to $X\otimes_k\bark$.


\subsection{The Euler characteristics, Swan conductors and cohomological rigidity}\label{ss:Euler}

There are two uniform choices of $V\in\Rep(\dualG)$ which have small dimensions. One is $V=\Ad=\dualg$, the adjoint representation; the other is $V=V_\thv$, the representation whose nonzero weights consist of short roots of $\dualg$. We call $V_\thv$ the {\em quasi-minuscule} representation of $\dualG$ (which coincides with the adjoint representation when $\dualG$ is simply-laced). Basic facts about $V_\thv$ are summarized in Appendix \S\ref{a:qmcom}. We will denote the number of long (resp.\ short) simple roots of $G$ by $r_\ell(G)$ (resp.\ $r_s(G)$) and call it the {\em long rank}  (resp. the {\em short rank} of $G$). See Lemma \ref{l:samenumber}) for equivalent descriptions of these numbers.

Let $\Kl^\thv_{\dualG}(\phi,\chi)$ (resp. $\Kl^\Ad_{\dualG}(\phi,\chi)$) be the local system associated to $\Kl_\dualG(\phi,\chi)$ and the quasi-minuscule representation (resp. the adjoint representation) of $\dualG$.
\begin{theorem}\label{th:Eulerchar}
\begin{enumerate}
\item []
\item $-\chi_c(\pline,\Kl^{\thv}_{\dualG}(\phi,\chi))$ equals the number of long simple roots of $G$. 
\item Suppose $\textup{char}(k)$ is good for $G$ when $G$ is not simply-laced, then $-\chi_c(\pline,\Kl^{\Ad}_{\dualG}(\phi,\chi))$ equals the rank of $G$.
\end{enumerate}
\end{theorem}

By the Grothendieck-Ogg-Shafarevich formula, we get
\begin{corollary}\label{c:Swan}
\begin{enumerate}
\item []
\item $\Swan_\infty(\Kl^\thv_\dualG(\phi,\chi))=r_s(\dualG)$.
\item $\Swan_\infty(\Kl^\Ad_\dualG(\phi,\chi))=r(\dualG)$, under the same assumption on $\textup{char}(k)$ as in Theorem \ref{th:Eulerchar}(2).
\end{enumerate}
\end{corollary}

The following subsections will be devoted to the proof of Theorem \ref{th:Eulerchar}. We first draw some consequences.

\begin{lemme}\label{l:noglobinv} Suppose the principal nilpotent element acts on $V$ without trivial Jordan block, then $H^0(\pline,\Kl^V_\dualG(\phi))=0$.
\end{lemme}
\begin{proof}
Suppose the contrary, then $\Kl^V_\dualG(\phi)$ contains the constant sheaf as a sub-local-system. By construction, $\Kl^V_\dualG(\phi)$ is pure. By \cite[Th.5.3.8]{BeilinsonBernsteinDeligne}, over $\geompline$, $\Kl^V_\dualG(\phi)$ is a direct sum of simple perverse sheaves. Hence the constant sub-sheaf must be a direct summand. 

On the other hand, by assumption, $V$ does not contain any direct summand under which the tame monodromy at $\{0\}\in\PP^1$ (i.e., a principal unipotent element, by Theorem \ref{Heckeeigensheaf}) acts trivially. This yields a contradiction.
\end{proof}

Let $i_0:\{0\}\hookrightarrow\PP^1;i_\infty:\{\infty\}\hookrightarrow\PP^1$ and $j:\pline\hookrightarrow\PP^1$ be the inclusions. For any local system $L$ on $\pline$, we abuse the notation $j_{!*}L$ to mean $j_{!*}(L[1])[-1]$.

As in the discussion before Theorem \ref{th:Swan-inv}, the local Galois groups $\Gal(K^{\sep}_\infty/K_\infty)$ and $\Gal(K^{\sep}_0/K_0)$, hence the inertia groups $\inert_0$ and $\iinf$, act on the corresponding geometric stalks of $\Kl^V_\dualG(\phi)$, defining representations on $V$ up to conjugacy.

\begin{proposition}\label{p:van} Under the same assumption on $\textup{char}(k)$ as in Theorem \ref{th:Eulerchar}(2), we have
\begin{enumerate}
\item \textup{(Cohomological rigidity)} $H^*(\PP^1_\bark,j_{!*}\Kl^\Ad_\dualG(\phi))=0$.
\item $\dualg^{\inert_\infty}=0$ (Note that $\dualg$ is the space of $V=\Ad$).
\end{enumerate}
The same statements hold for $\Kl^\thv_\dualG(\phi)$ in place of $\Kl^\Ad_\dualG(\phi)$, with no restriction on $\textup{char}(k)$.
\end{proposition}
\begin{proof} The statements being geometric, we will ignore Tate twists in this proof.
Since the adjoint representation $\dualg$ is self-dual, the local system $\Kl^\Ad_\dualG(\phi)$ is also self-dual. Fixing such an isomorphism $\Kl^\Ad_\dualG(\phi)\isom(\Kl^\Ad_\dualG(\phi))^\vee$, 
we get isomorphisms $H^0(\PP^1,j_{!*}\Kl^\Ad_\dualG(\phi))\cong H^2(\PP^1,j_{!*}\Kl^\Ad_\dualG(\phi))^\vee$ and $H^0(\pline,\Kl^\Ad_\dualG(\phi))\cong H^2_c(\pline,\Kl^\Ad_\dualG(\phi))^\vee$.

We have a distinguished triangle in $D^b_c(\PP^1,\Ql)$
\begin{equation*}
j_!\Kl^\Ad_\dualG(\phi)\to j_{!*}\Kl^\Ad_\dualG(\phi)\to H^0i_0^*j_*\Kl^\Ad_\dualG(\phi)\oplus H^0i_\infty^*j_*\Kl^\Ad_\dualG(\phi)\to
\end{equation*}
which induces a long exact sequence
\begin{eqnarray}
\label{long1}0=H^0_c(\pline,\Kl^\Ad_\dualG(\phi))\to &H^0(\PP^1,j_{!*}\Kl^\Ad_\dualG(\phi))&\to \dualg^{\inert_0}\oplus \dualg^{\inert_\infty}\xrightarrow{d}\\
\label{long2}\xrightarrow{d}H^1_c(\pline,\Kl^\Ad_\dualG(\phi))\to &H^1_c(\PP^1,j_{!*}\Kl^\Ad_\dualG(\phi))&\to0\to\\
\label{long3}\to H^2_c(\pline,\Kl^\Ad_\dualG(\phi))\to &H^2_c(\PP^1,j_{!*}\Kl^\Ad_\dualG(\phi))&\to0
\end{eqnarray}
By \eqref{long3} we first conclude
\[
H^2(\PP^1,j_{!*}\Kl^\Ad_\dualG(\phi))\cong H^2_c(\pline,\Kl^\Ad_\dualG(\phi))\cong H^0(\pline,\Kl^\Ad_\dualG(\phi))^\vee=0
\]
where the vanishing of the last term follows by applying Lemma \ref{l:noglobinv} to $V=\dualg$. By duality, $H^0(\PP^1,j_{!*}\Kl^\Ad_\dualG(\phi))=0$. By \eqref{long1} and the vanishing of $H^0$'s, the connecting homomorphism $d$ is injective.

On the one hand, by Theorem \ref{Heckeeigensheaf}(2), $\dim\dualg^{\inert_0}=r(\dualG)$ because $\inert_0$ acts on $\dualg$ through a principal unipotent element.
On the other hand, by Theorem \ref{th:Eulerchar},  $\dim H^1_c(\pline,\Kl^\Ad_\dualG(\phi))=-\chi_c(\pline,\Kl^\Ad_\dualG(\phi))=r(\dualG)$  (here we again used the vanishing of $H^i_c(\pline,\Kl^\Ad_\dualG(\phi))$ for $i=0,2$). Therefore $d$ must be an isomorphism. This implies $\dualg^{\inert_\infty}=0$. By \eqref{long2}, we conclude that $H^1(\PP^1,j_{!*}\Kl^\Ad_\dualG(\phi))=0$.

The statement for $\Kl^\thv_\dualG(\phi)$ is proved in the same way. We may apply Lemma \ref{l:noglobinv} to $V_\thv$, because of equation \eqref{V0V1} from the proof of Lemma \ref{l:samenumber}.
\end{proof}

\begin{remark}
Following Katz \cite[\S 5.0]{KatzRigid}, we call a $\dualG$-local system $L$ on $\geompline$ {\em cohomologically rigid} if $H^1(\PP^1,j_{!*}L^\Ad)=0$. We think of $H^1(\PP^1,j_{!*}L^\Ad)$ as the space of infinitesimal deformations of the $\dualG$-local system $L$ with fixed isomorphism type on the formal punctured discs around $0$ and $\infty$, although the notion of such deformations has not been defined. Proposition \ref{p:van} implies that $\Kl_\dualG(\phi)$ is cohomologically rigid, which provides evidence for its {\em physical rigidity}. For more precise conjectures, see \S\ref{s:fonc}.
\end{remark}

\subsection{General method of calculation}\label{gen}
In the following calculation, we will only consider the neutral component of $\Bun_{G(1,2)}$. Hence we may assume $G$ is simply-connected.

We denote by $\star\in \Bun_{G(1,2)}$ the point corresponding to the trivial bundle. In order to compute the sheaf $\Kl^V_\dualG(\phi,\chi)$ we restrict the convolution diagram (\ref{HeckeOp}) to $\star \times \P1min2\subset \Bun_{G(1,2)} \times \P1min2$:
$$\xymatrix{
 & \GR\ar[dl]_-{\pr_1}\ar[dr]^-{\pr_2} &  \\
\Bun_{G(1,2)}&  & \star \times \P1min2.\\
}$$
By proper base change we have $\Kl^V_\dualG(\phi,\chi)=\pr_{2,!} (\pr_1^* A \tensor \IC_V)$.
Furthermore, $A$ is supported only on the big cell $j\colon T \times I(1)/I(2) \subset \Bun_{(G(1,2))}$. Let us denote by $\Gr^\circ\subset \Gr$ the inverse image of the big cell and write $I(1)/I(2) \cong U/[U,U] \times U_{-\theta}$. We write $\phi=(\phi_+,\phi_0)$ according to this decomposition.

We obtain a diagram
$$\xymatrix{
                                & \GR^\circ\ar@{^(->}[r]\ar@/^3pc/[drr]^-{\pr^\circ_2}\ar[dl]_-{(f_T,f_0,f_+)} & \GR\ar[dl]_-{\pr_1}\ar[dr]^-{\pr_2} &  \\
 T\times U_{-\theta} \times U/[U,U]\ar@{^(->}[r] & \Bun_{G(1,2)} & & \P1min2
}$$
By proper base change,
\begin{equation}\label{preKlV}
\Kl^V_\dualG(\phi,\chi)=\pr^\circ_{2,!} (f_T^*\Kum_\chi f_0^*\AS_{\phi_0}\otimes f_+^*\AS_{\phi_+} \tensor \IC_V).
\end{equation}
So we need to describe the subspace $\GR^\circ$ and compute the maps $f_T$, $f_0$, $f_+$.

First, let us denote by $\GR^\triv\subset \GR$ the preimage of the trivial $G$-bundle under the forgetful map $\GR\to\Bun_G$. Let us fix a point $x\in \P1min2$ and chose $t_x:=1-\frac{t}{x}$ as a local parameter at $x$. So $t_x=1$ corresponds to $t=0$ and $t_x=\infty$ corresponds to $t=\infty$.

The ind-scheme $\Gr^\triv_x$ classifies isomorphism classes of pairs $(\calE,\varphi)$ such that $\calE$ is a trivial $G$-bundle on $\PP^1$ and $\varphi\colon\calE|_{\Px}\isom G\times(\Px)$ is an isomorphism. We can rigidify this moduli problem so that $\Gr^\triv_x$ classifies an automorphism of the trivial $G$-bundle on $\PP^1\backslash\{x\}$ which is identity at $\infty$. Hence we can identify
\begin{equation}\label{bcellGr1}
\Gr_x^\triv\isom G[t_x^{-1}]_1:=\ker(G[t_x^{-1}]\xrightarrow{\ev(t_x=\infty)}G)
\end{equation}
where $G[t_x^{-1}]$ is the ind-scheme whose $R$-points are $G(R[t_x^{-1}])$.

\begin{remark} There is a dilation action of $\grot$ on $\pline$, which extends to $\GR$: $\lambda\in\grot$ sends $(x,\cE,\varphi)\mapsto(\lambda x,\lambda^{-1,*}\cE,\lambda^{-1,*}\varphi)$, and stabilizes $\GR^\triv$. Note that the local coordinate $t_x$ is invariant under the simultaneous dilation on $t$ and $x$. 

Let $\Gr^\triv\subset\Gr=G((\tau))/G[[\tau]]$ be an abstract copy of the affine Grassmannian, not referring to any point on $\pline$. The dilation action together with \eqref{bcellGr1} gives a trivialization of the family $\GR^\triv$ over $\pline$:
\begin{eqnarray}\label{GRprod}
\GR^\triv&\cong & (\pline)\times G[\tau^{-1}]_1\cong(\pline)\times\Gr^\triv\\
\notag(x,g(t_x^{-1}))&\mapsto& (x,g(\tau^{-1}))
\end{eqnarray}
\end{remark}

By definition the big cell in $\Bun_{G(1,2)}$ is obtained from the action of $T\times I(1)/I(0)$ on the trivial bundle. Let $\widehat{\ev}_0\colon G[t_x^{-1}] \to G[[t]]$ and $\widehat{\ev}_\infty\colon G[t_x^{-1}] \to G[[t^{-1}]]$ denote the expansions around $0$ and $\infty$ and $\ev_0,\ev_\infty\colon G[t_x^{-1}] \to G$ the evaluation at $t=0$ and $t=\infty$.

The $G(1,2)$-level structure on the $G$-bundle $(\calE,\varphi)\in \Gr_{x}$ is obtained from the trivialization $\varphi$. For $g\in G[t_x^{-1}]$, the composition $\cE|_{\Px}\map{\varphi} G \times \Px \map{g} G \times \Px$ extends to an isomorphism $\cE \to G \times \bP^1$. Thus $g$ defines the level structure on the trivial bundle defined by $\ev_0(g),\ev_\infty(g)$. We have a commutative diagram:
$$\xymatrix{
G[t_x^{-1}] \ar[rr]^-{\widehat{\ev}_0,\widehat{\ev}_\infty} \ar[dr] & & G\backslash (G[[t]] \times G[[t^{-1}]])\ar[dl]\\
& \Bun_{G(1,2)} &. 
}$$
Here $G=\Aut(G\times \bP^1)$ acts diagonally on $G[[t]] \times G[[t^{-1}]]$.

Thus we find that $g\in G[t_x^{-1}]$ lies in $\Gr^\circ_x$ if and only if $(\ev_0(g),\ev_\infty(g))\in G\backslash G(B^\opp\times U)$. Using the identifications \eqref{bcellGr1} and \eqref{GRprod}, we get
\begin{eqnarray}\notag
\Gr^\circ_x&\isom&\{g(t_x^{-1})\in G[t_x^{-1}]_1|\ev_0(g)\in UB^\opp\};\\
\label{bcellFl1}
\GR^\circ&\isom&(\pline)\times\{g(\tau^{-1})\in G[\tau^{-1}]_1|g(1)\in UB^\opp\}\\
\notag&\subset&\GR^\triv=(\pline)\times G[\tau^{-1}]_1.
\end{eqnarray}

\begin{lemme}\label{l:ff1}
For $(x,g(\tau^{-1}))\in\GR^\circ$ under the parametrization in \eqref{bcellFl1}, write $g(1)=ub^\opp$ for $u\in U$ and $b^\opp\in B^\opp$, then we have
\begin{eqnarray*}
f_T(x,g)&=& b^\opp\mod U^\opp \in T;\\
f_+(x,g)&=&u^{-1}\mod[U,U]\in U/[U,U];\\
f_0(x,g)&=&xa_{-\theta}(g)\in U_{-\theta}\cong\gnth,
\end{eqnarray*}
where $a_{-\theta}:G[\tau^{-1}]\to\gnth$ sends $g$ to the $\gnth$-part of the tangent vector $\dfrac{dg(\tau^{-1})}{d(\tau^{-1})}\bigg|_{\tau^{-1}=0}\in\frg$.
\end{lemme}
\begin{proof}
The formulas for $f_T,f_+$ follow from our description in \eqref{bcellFl1}. By definition $f_0(x,g)$ is obtained by expanding $g(t_x^{-1})$ at $t=\infty$ using the local parameter $t^{-1}$ and taking the $\gnth$-part of the coefficient of $t^{-1}$. Note that
\begin{equation*}
\frac{dg(t_x^{-1})}{d(t^{-1})}\bigg|_{t^{-1}=0}=x\frac{dg(t_x^{-1})}{d(t_x^{-1})}\bigg|_{t_x^{-1}=0}\in\frg.
\end{equation*}
Moreover, under the identification \eqref{GRprod}, the parameter $t_x$ corresponds to $\tau$, therefore $f_0(x,g)=xa_{-\theta}(g)$. This proves the lemma.
\end{proof}

%

Let $\AS_{U,-\phi_+}$ be the pull-back of $\AS_\psi$ via $U\to U/[U,U]\xrightarrow{-\phi_+}\GG_a$. Let $\Kum_{\barB,\chi}$ be the pull-back of the Kummer local system $\Kum_\chi$ via $B^\opp\to T$. 
Let $j_{UB^\opp}:U\barB\hookrightarrow G$ be the inclusion, and denote
\begin{equation*}
J=J_{-\phi_+,\chi}:=j_{U\barB,!}(\AS_{U,-\phi_+}\boxtimes\Kum_{\barB,\chi})\in D^b_c(G,\Ql).
\end{equation*}

\begin{remark}
According to \cite{BBM}, or the argument of Lemma \ref{clean} (in this paper), we have the cleanness property of $J$: 
\begin{equation*}
j_{U\barB,!}(\AS_{U,-\phi_+}\boxtimes\Kum_{\barB,\chi})\isom j_{U\barB,*}(\AS_{U,-\phi_+}\boxtimes\Kum_{\barB,\chi})
\end{equation*}
We can view $J=J_{-\phi_+,\chi}\in D^b_c(G,\Ql)$ as a finite-field analog of the automorphic sheaf $A_{\phi,\chi}$ in Definition \ref{eigensheaf} and Remark \ref{r:multchar}. 
\end{remark}

With this notation, and using the identification \eqref{GRprod}, the formula \eqref{preKlV} becomes
\begin{equation}\label{KlV1}
\Kl^V_{\dualG}(\phi,\chi)\isom\pi^\triv_!(\IC^\triv_V\otimes f_0^*\AS_{\phi_0}\otimes\ev_{\tau=1}^*J).
\end{equation}
where $\pi^\triv:\GR^\triv=\pline\times G[\tau^{-1}]_1\to\pline$ is the projection. In the sequel, we often write $\ev_{\tau=1}$ simply as $\ev$.


The ultimate goal of this section is to calculate the Euler characteristics of $\Kl^V_\dualG(\phi,\chi)$ for $V=V_\thv$ and $\dualg$. Here we make a few reduction steps for general $V$. By \eqref{KlV1}, we need to calculate
\begin{equation}\label{chiKlV}
\chi_c(\pline,\Kl^V_\dualG(\phi,\chi))=\chi_c(\pline\times\Gr^{\triv},\IC^\triv_V\otimes f_0^*\AS_{\phi_0}\otimes\ev^*J).
\end{equation}

According to whether $a_{-\theta}$ vanishes or not, we decompose $\Gr^\triv\cong G[\tau^{-1}]_1$ into $\Gr^{\triv,a\neq0}$ and $\Gr^{\triv,a=0}$. Over $\pline\times \Gr^{\triv,a=0}$, the complex in \eqref{chiKlV} is constant along $\pline$, hence the Euler characteristic is 0. On the other hand, we have a change of variable isomorphism
\begin{equation*}
\pline\times\Gr^{\triv,a\neq0}\ni(x,g)\mapsto(xa_{-\theta}(g),g)\in\GG_m\times \Gr^{\triv,a\neq0}.
\end{equation*}
Under this isomorphism, $f_0=xa_{-\theta}$ becomes the projection to the $\GG_m$-factor, and we can apply the K\"unneth formula. Summarizing these steps, we get
\begin{eqnarray}\label{chiKlV2}
&&\chi_c(\pline,\Kl^V_\dualG(\phi,\chi))\\
\notag&=&\chi_c(\pline\times\Gr^{\triv,a\neq0},\IC^\triv_V\otimes f_0^*\AS_{\phi_0}\otimes\ev^*J)\\
\notag&=&\chi_c(\GG_m,\AS_{\phi_0})\chi_c(\Gr^{\triv,a\neq0},\IC^\triv_V\otimes\ev^*J)\\
\notag&=&-\chi_c(\Gr^{\triv,a\neq0},\IC^\triv_V\otimes\ev^*J)\\
\notag&=&-\chi_c(\Gr^\triv\backslash\star,\IC^\triv_V\otimes\ev^*J)+\chi_c(\Gr^{\triv,a=0}\backslash\star,\IC^\triv_V\otimes\ev^*J).
\end{eqnarray}
The last equality is because the base point $\star\in\Gr$ belongs to $\Gr^{\triv,a=0}$.

\begin{lemme}\label{l:chivan}
$\chi_c(\Gr^\triv\backslash\star,\IC^\triv_V\otimes\ev^*J)=0$.
\end{lemme}
\begin{proof}
Since $\IC_V$ is constant along the strata $\Gr_\lambda$,  it suffices to show that $\chi_c(\Gr^\triv_\lambda,\ev^*J)=0$ for dominant coweights $\lambda\neq0$.

Denote by $\conv{G}$ the convolution product on $D^b_c(G,\Ql)$: for $K_1,K_2\in D^b_c(G,\Ql)$,
\begin{equation*}
K_1\conv{G}K_2:=m_!(K_1\boxtimes K_2)
\end{equation*}
where $m:G\times G\to G$ is the multiplication map. Let $K_\lambda:=\ev_{!}\const{\Gr^\triv_\lambda}\in D^b_c(G)$. Then $R\Gamma_c(\Gr^\triv_\lambda,\ev^*J)$ is the stalk at $e\in G$ of the convolution $\AS_{U,\phi_+}\conv{G}K_\lambda\conv{G}\Kum_{\barB,\chi^{-1}}$. Since $\ev=\ev_{\tau=1}$ is $G$-equivariant (under conjugation), $K_\lambda$ carries a natural $G$-equivariant structure. Hence $$K_\lambda\conv{G}\Kum_{\barB,\chi^{-1}}\cong\Kum_{\barB,\chi^{-1}}\conv{G}K_\lambda.$$ In particular, over each Bruhat stratum $\barB w\barB$, $K_\lambda\conv{G}\Kum_{\barB,\chi^{-1}}$ has the form $S_w\otimes L_w$, where $L_w$ is a local system on $\barN w\barB$, and $S_w$ (a complex of $\Ql$-vector spaces) is the stalk of $K_\lambda\conv{G}\Kum_{\barB,\chi^{-1}}$ at $\dot{w}\in N_G(T)$ (a representative of $w\in W$).

Therefore $\AS_{U,\phi_+}\conv{G}K_\lambda\conv{G}\Kum_{\barB,\chi^{-1}}$ is a successive extension of $S_w\otimes(\AS_{U,\phi_+}\conv{G}L_w)$ for various $w\in W$. To prove the vanishing of the Euler characteristic of the stalks of $\AS_{U,\phi_+}\conv{G}K_\lambda\conv{G}\Kum_{\barB,\chi^{-1}}$, it suffices to show that $\chi(S_w)=0$ for all $w$.

By definition, 
\begin{eqnarray}
\notag S_w&=&(K_\lambda\conv{G}\Kum_{\barB,\chi^{-1}})_{\dot{w}}\\
\notag&=&R\Gamma_c(w\barB,K_\lambda\otimes\Kum_{\barB,\chi})\\
\label{lastline}&=&R\Gamma_c(\Gr^\triv_\lambda\cap\ev^{-1}(w\barB),\ev^*\Kum_{\barB,\chi}).
\end{eqnarray}
The $T$-action on $\Gr^\triv=G[\tau^{-1}]_1$ by conjugation preserves $\Gr^\triv_\lambda\cap\ev^{-1}(w\barB)$. The only $T$-fixed points on $\Gr$ are $\tau^\mu$ for $\mu\in\xcoch(T)$, which do not belong to any $\Gr^\triv_\lambda$ ($\lambda\neq0$). Moreover, the local system $\ev^*\Kum_{\barB,\chi}$ is monodromic under $T$-conjugation: there exists $m\geq1$ (prime to $p$) such that $\ev^*\Kum_{\barB,\chi}$ is equivariant under the $m$-th power of the $T$-conjugation. Since the $m$-th power of $T$-conjugation still has no fixed point on $\Gr^\triv_\lambda\cap\ev^{-1}(w\barB)$, the Euler characteristic in \eqref{lastline} is zero. This proves the lemma.
\end{proof}

\begin{corollary}\label{c:chiKlV} Let $\xcoch(T)^+\subset\xcoch(T)$ be the dominant coweights, and $V(\lambda)\subset V$ be the weight space for $\lambda\in\xcoch(T)$. Then
\begin{equation*}
\chi_c(\pline,\Kl^V_\dualG(\phi,\chi))=\sum_{\lambda\in\xcoch(T)^+,\lambda\neq0}\dim V(\lambda)\chi_c(\Gr^{\triv,a=0}_\lambda,\ev^*J).
\end{equation*}
\end{corollary}
\begin{proof}
By \cite[Theorem 6.1]{Lusztig}, the Euler characteristic of the stalks of $\IC_V$ along $\Gr_\lambda$ is $\dim V(\lambda)$. Therefore
\begin{eqnarray*}
\chi_c(\Gr^{\triv,a=0}\backslash\star,\IC^\triv_V\otimes\ev^*J)=\sum_{\lambda\in\xcoch(T)^+,\lambda\neq0}\dim V(\lambda)\chi_c(\Gr^{\triv,a=0}_\lambda,\ev^*J).
\end{eqnarray*}
Combining this with \eqref{chiKlV2} and Lemma \ref{l:chivan} yields the desired identity.
\end{proof}

\subsection{Quasi-minuscule Schubert variety}\label{qmSch}
For a coweight $\lambda\in\xcoch(T)$, let $P_{\lambda}\subset G$ be the parabolic generated by $T$ and the root spaces $U_\alpha$ for $\jiao{\alpha,\lambda}\geq0$.

For each root $\alpha$ and $i\in\ZZ$, let $U_{\alpha,\geq i}\subset U_\alpha((\tau))$ be the subgroup whose $R$-points are $U_\alpha(\tau^iR[[\tau]])$. Let $U_{\alpha,i}=U_{\alpha,\geq i}/U_{\alpha,\geq i+1}$.

Now let $\lambda\in\xcoch(T)^+$ and consider the (open) Schubert variety $\Gr_\lambda\subset\Gr$. By \cite[Lemme 2.3 and discussions preceding it]{NgoPolo}, we have (fixing a total ordering of $\Phi^+$)
\begin{eqnarray}\label{GPU}
G\twtimes{P_{-\lambda}}\left(\prod_{\jiao{\alpha,\lambda}\geq2}U_{\alpha,\geq1}/U_{\alpha,\geq\jiao{\alpha,\lambda}}\right)&\isom&\Gr_\lambda\\
\notag(g, u)&\mapsto&gu\tau^\lambda.
\end{eqnarray}
The action of $P_{-\lambda}$ on the product in \eqref{GPU} is given by the adjoint action of $P_{-\lambda}$ on $\prod_{\alpha}U_{\alpha,\geq1}$ followed by the projection onto factors $U_{\alpha,\geq1}$ for $\jiao{\alpha,\lambda}\geq2$.


Now consider the special case $\lambda=\thv$, the dominant short coroot. We need to recall the description of the quasi-minuscule Schubert variety $\Gr_{\thv}$ from \cite{NgoPolo}. We write $\Pnth$ for $P_{-\thv}$.


\begin{lemme}\label{l:qm}
The open locus $\Gr^\triv_\thv\subset\Gr_\thv$ consists of a single $G$-orbit. More precisely, the morphism
\begin{eqnarray*}
G\twtimes{\Pnth}U^\times_{-\theta,-1}&\to& G[\tau^{-1}]_1=\Gr^\triv\\
(g,u_{-\theta}(c\tau^{-1}))&\mapsto& \Ad(g)u_{-\theta}(c\tau^{-1})
\end{eqnarray*}
gives an isomorphism onto $\Gr^\triv_\thv$ (here $\Pnth$ acts on $U_{-\theta,-1}$ through adjoint action and $U^\times_{\alpha,i}=U_{\alpha,i}\backslash\{1\}$).
\end{lemme}
\begin{proof}
This is essentially \cite[Lemme 7.2]{NgoPolo}. Applying \eqref{GPU} to $\lambda=\thv$, we find the product in \eqref{GPU} consists of only one term $U_{\theta,1}$ (since $\jiao{\alpha,\thv}\leq1$ for any root $\alpha\neq\theta$). Therefore $\Gr^\triv_\thv\cong G\twtimes{\Pnth}U^\triv_{\theta,1}\cdot\tau^\thv$ for some proper open subset $U^\triv_{\theta,1}\subset U_{\theta,1}$ stable under $\Pnth$ ($\Gr^\triv_\thv$ cannot be equal to $\Gr_\thv$ because $\tau^\thv\notin\Gr^\triv$). Since $\Pnth$ acts on $U_{\theta,1}$ via dilation, $U^\triv_{\theta,1}$ must be $U^\times_{\theta,1}$. Hence
\begin{equation*}
\Gr^\triv_\thv\cong G\twtimes{\Pnth}U^\times_{\theta,1}\cdot\tau^\thv.
\end{equation*}

On the other hand, the following calculation in the $\SL_2$-subgroup defined by $\theta$:
\begin{eqnarray}\label{sl2}
\mat{1}{c\tau}{}{1}\mat{\tau}{}{}{\tau^{-1}} = \mat{\tau}{c}{}{\tau^{-1}} = \mat{1}{}{c^{-1}\tau^{-1}}{1}\mat{\tau}{c}{-c^{-1}}{}
\end{eqnarray}
shows that $u_{\theta}(c\tau)\tau^\thv\in\Gr$ equals $u_{-\theta}(c^{-1}\tau^{-1})\in U^\times_{-\theta,-1}\subset G[\tau^{-1}]_1$ for $c$ invertible. This proves the lemma.
\end{proof}

The Bruhat decomposition for $G=\sqcup_{w\in W/W_\theta}Uw\Pnth$ gives a decomposition
\begin{equation*}
G\twtimes{\Pnth}U^\times_{-\theta,-1}=\bigsqcup_{w\in W/W_\theta}Uw\Pnth\twtimes{\Pnth}U^\times_{-\theta,-1}.
\end{equation*}
The stratum $Uw\Pnth\twtimes{\Pnth}U^\times_{-\theta,-1}$ has image $\Ad(U)U^\times_{-w\theta,-1}$ in $G[\tau^{-1}]_1$. Since $w\mapsto-w\theta$ sets up a bijection between $W/W_\theta$ and the set of long roots of $G$, we can rewrite the above decomposition as
\begin{equation*}
\Gr^\triv_\thv=\bigsqcup_{\beta\textup{ long root}}\Ad(U)U^\times_{\beta,-1}
\end{equation*}
For each $w\in W/W_\theta$, we have an isomorphism (fixing total ordering on $\Phi^+$)
\begin{equation*}
\prod_{\alpha\in\Phi^+,\jiao{w^{-1}\alpha,\thv}>0}U_\alpha\cong Uw\Pnth/\Pnth.
\end{equation*}
Therefore, for $\beta=-w\theta$, the stratum $\Ad(U)U^\times_{\beta,-1}$ can be written as
\begin{eqnarray}\label{qmstratum}
\Ad(U)U^\times_{\beta,-1}&\cong&\prod_{\alpha\in\Phi^+,\jiao{\alpha,\beta^\vee}<0}U_\alpha\times U^\times_{\beta,-1}\\
\notag\Ad(\prod u_\alpha)u_{\beta}(c\tau^{-1})&\leftrightarrow&(\prod u_\alpha,u_{\beta}(c\tau^{-1}))
\end{eqnarray}

\begin{lemme}\label{l:qmf01} The function $a_{-\theta}:\Gr^\triv_\thv\to\gnth$ in Lemma \ref{l:ff1} restricted on each stratum $\Ad(U)U^\times_{\beta,-1}\subset\Gr^\triv_\thv$ is given by
\begin{equation*}
a_{-\theta}(\Ad(u)u_{\beta}(c\tau^{-1}))=\begin{cases} c\bbx_{-\theta} &  \textup{if }\beta=-\theta,\\ 0 & \textup{otherwise}.\end{cases}
\end{equation*}
Here $\bbx_{-\theta}\in\gnth$ corresponds to $u_{-\theta}(1)\in U_{-\theta}$. In particular,
\begin{equation}\label{qma=0}
\Gr^{\triv,a=0}_\thv=\bigsqcup_{\beta\textup{ long},\beta\neq-\theta}\Ad(U)U^\times_{\beta,-1}
\end{equation}
\end{lemme}
\begin{proof}
By definition, the value of $a_{-\theta}$ on $\Ad(u)u_{\beta}(c\tau^{-1})$ is $$\frac{d}{d(\tau^{-1})}\Ad(u)u_{\beta}(c\tau^{-1})|_{\tau^{-1}=0}\in\Ad(U)\frg_{\beta}.$$ If $\beta\neq-\theta$, $\Ad(u)\frg_{\beta}$ only involves roots $\geq\beta$. If $\beta=-\theta$, then the above derivative equals $\Ad(u)c\bbx_{-\theta}\in\frg$, whose $\gnth$-part is $c\bbx_{-\theta}\in\gnth$.
\end{proof}

%
%

\subsection{Proof of Theorem \ref{th:Eulerchar}(1)} Applying Corollary \ref{c:chiKlV} to $V=V_{\thv}$, which only has one nonzero dominant weight $\thv$, we get
\begin{equation*}
\chi_c(\pline,\Kl^\thv_\dualG(\phi,\chi))=\chi_c(\Gr^{\triv,a=0}_\thv,\ev^*J).
\end{equation*}
By the decomposition \eqref{qma=0} in Lemma \ref{l:qmf01}, we only need to calculate $\chi_c(\Ad(U)U^\times_{\beta,-1},\ev^*J)$ for long roots $\beta\neq-\theta$. Theorem \ref{th:Eulerchar}(1) thus follows from

\begin{claim}
Suppose $\beta\neq-\theta$ is a long root, then
\begin{equation*}
\chi_c(\Ad(U)U^\times_{\beta,-1},\ev^*J)=\begin{cases}-1 & \beta \textup{ is a simple long root}\\ 0 & \textup{otherwise}\end{cases}
\end{equation*}
\end{claim}

To prove the claim, we distinguish three cases:

\noindent{\bf Case I}: $\beta$ is positive but not simple.

Since $\beta>0$, we have $\ev(\Ad(U)U_{\beta,-1}^\times)\subset U$. Since $\beta$ is not simple, the image of $\ev(\Ad(U)U_{\beta,-1}^\times)$ in $U/[U,U]$ is trivial, hence $\ev^*J$ is the constant sheaf on $\Ad(U)U_{\beta,-1}^\times$. By \eqref{qmstratum}, $\Ad(U)U^\times_{\beta,-1}$ has a factor $U^\times_{\beta,-1}\cong\GG_m$, hence
\begin{equation*}
\chi_c(\Ad(U)U_{\beta,-1}^\times,\ev^*J)=\chi_c(\Ad(U)U_{\beta,-1}^\times)=0.
\end{equation*}

\noindent{\bf Case II}: $\beta=\alpha_i$ is a simple long root.

For $u\in U$, $\ev(\Ad(u)u_{\beta}(c\tau^{-1}))$ has image $u_{\alpha_i}(c)\in U/[U,U]$. In terms of the coordinates in \eqref{qmstratum}, $\ev^*J$ is the pull-back of $\AS_{\phi_i}$ from the $U^\times_{\beta,-1}\cong U^\times_{\alpha_i}$-factor. Hence 
\begin{equation*}
\chi_c(\Ad(U)U^\times_{\beta,-1},\ev^*J)=\chi_c(U^\times_{\alpha_i},\AS_{\phi_i})=-1.
\end{equation*}

\noindent{\bf Case III}: $\beta$ is negative and $\beta\neq-\theta$.

Since $\beta\neq-\theta$, there exists a simple root $\alpha_i$ such that $\beta-\alpha_i$ is still a root. Then $\alpha=-\beta+\alpha_i$ is a positive root. Moreover, since $\jiao{\alpha_i,\beta^\vee}\leq1$ (because $\alpha_i\neq\beta$), we have $\jiao{\alpha,\beta^\vee}=\jiao{-\beta+\alpha_i,\beta^\vee}=-2+\jiao{\alpha_i,\beta^\vee}<0$. Hence $U_\alpha$ appears in the decomposition \eqref{qmstratum}.

Using \eqref{qmstratum},  we write an element in $\Ad(U)U^\times_{\beta,-1}$ as $\Ad(u)u_{\beta}(c_\beta\tau^{-1})$, where $u=u^\alpha u_\alpha(c_\alpha)$ and $u^\alpha\in\bbA:=\prod_{\alpha'>0,\jiao{\alpha',\beta^\vee}<0,\alpha'\neq\alpha}U_{\alpha'}$. Note that
\begin{eqnarray}
\label{comm}\ev(\Ad(u)u_{\beta}(c_\beta\tau^{-1}))=uu_\beta(c_\beta)u^{-1}=u^\alpha[u_{\alpha}(c_\beta),u_{\beta}(c_\beta)]u_{\beta}(c_\beta)u^{\alpha,-1}.
\end{eqnarray}
Since $\alpha\neq\pm\beta$, we can apply Chevalley's commutator relation \cite[p.36,(4)]{Ch} to conclude that $[u_{\alpha}(c_\beta),u_{\beta}(c_\beta)]$ is a product of elements in the root groups $U_{i\alpha+j\beta}$ for $i,j\in\ZZ_{>0}$. Our assumptions that (i) $\alpha+\beta$ is simple and (ii) $\beta$ is a long root imply that any such $i\alpha+j\beta$ is positive (if it is a root), and the only simple root of this form is $\alpha+\beta$. Therefore, $[u_{\alpha}(c_\alpha),u_{\beta}(c_\beta)]\in U$, and its image in $U/[U,U]$ is $u_{\alpha_i}(\epsilon c_\alpha c_\beta)$, where $\epsilon=\pm1$.

By \eqref{comm}, $\Ad(u)u_\beta(c_\beta)\in U\barB$ if and only if $u_\beta(c_\beta)u^{\alpha,-1}\in U\barB$. Moreover, when $\Ad(u)u_{\beta}(c_\beta)\in U\barB$, its image in $U/[U,U]$ is $u_{\alpha_i}(\epsilon c_\alpha c_\beta)$ times another element which only depends on $c_\beta$ and $u^{\alpha}$.

Under the decomposition \eqref{qmstratum}, $\Ad(U)U^\times_{\beta,-1}\cong U_\alpha\times\bbA$. Let $\pr_\bbA:\Ad(U)U^\times_{\beta,-1}\to\bbA$ be the projection. By the above discussion, $\ev^*J$ restricted to the fibers of $\pr_\bbA$ are isomorphic to Artin-Schreier sheaves on $U_\alpha$. Therefore $\pr_{\bbA,!}\ev^*J=0$, hence
\begin{equation*}
H^*_c(\Ad(U)U^\times_{\beta,-1},\ev^*J)=H^*_c(\bbA,\pr_{\bbA,!}\ev^*J)=0.
\end{equation*}

This proves the Claim, and completes the proof of Theorem \ref{th:Eulerchar} (1).


\subsection{The adjoint Schubert variety}\label{adjSch}
In this subsection, let $G$ be non-simply-laced. We always assume that char$(k)$ is a good prime for $G$. So char$(k)>2$ when $G$ is of type $B_n,C_n$, and char$(k)>3$ when $G$ is of type $F_4,G_2$.

Let $\gamma$ be the short dominant root of $G$, and $\gav$ be the corresponding long coroot. We define $\Png$, $W_\gamma$, etc. in the same way as $\Pnth$, $W_\theta$, etc. According to the possible values of $\jiao{\alpha,\gav}$ for roots $\alpha\in\Phi$, we have two cases:
\begin{itemize}
\item Type ($B_n$,$C_n$,$F_4$): $|\jiao{\alpha,\gav}|=0,1,2$;
\item Type ($G_2$): $|\jiao{\alpha,\gav}|=0,1,2,3$.
\end{itemize}
Let $\Phi_n^\gamma:=\{\alpha\in\Phi|\jiao{\alpha,\gav}=n\}$.

\begin{lemme}\label{l:adj} Suppose $G$ is of type $B_n, C_n$ or $F_4$, then $\Gr^\triv_{\gav}$ consists of a single $G$-orbit. More precisely, let
\begin{equation*}
V_{-\gamma}=\prod_{\jiao{\alpha,\gav}=2}U_{-\alpha,-1},
\end{equation*}
We identify $\Vng$ with its Lie algebra. The adjoint action of $\Png$ on $\Vng$ (which factors through its Levi factor $L_\gamma$) stabilizes a unique quadric $\Qng$ defined by a nondegenerate quadratic form $\qng$ on $\Vng$, so that the action of $L_{\gamma}$ factors through $\GO(\Vng,\qng)$ \footnote{$\GO(\Vng,\qng)$ consists of invertible linear automorphisms of $\Vng$ preserving $\qng$ up to a scalar.}, with a dense orbits $\Vng-\Qng$.

The natural morphism
\begin{equation}\label{Gradj}
G\twtimes{\Png}(\Vng-\Qng)\ni(g,v)\mapsto\Ad(g)v\in G[\tau^{-1}]_1
\end{equation}
gives an isomorphism $G\twtimes{\Png}(\Vng-\Qng)\cong\Gr^\triv_{\gav}$.
\end{lemme}
\begin{proof}
In the case $G$ is of type $B$, $C$ or $F_4$, the product in \eqref{GPU} becomes $\Lambda=\prod_{\jiao{\alpha,\gav}=2}U_{\alpha,1}$ (a commutative unipotent group), which can be identified with its Lie algebra. Since $\Gr^\triv_{\gav}$ is stable under $G$-conjugation, it takes the form $G\twtimes{\Png}\Lambda^\triv$ for some $\Png$-stable open subset $\Lambda^\triv\subset\Lambda$. The action of $\Png$ on $\Lambda$ factors through the Levi quotient $L_\gamma$.

The vector space $\Lambda$ carries a bilinear form
\begin{equation}\label{pairing}
(x,y)_\Lambda:=(x,\Ad(w_\gamma)y)_{\frg}
\end{equation}
where $w_\gamma$ is the image of $\mat{0}{1}{-1}{0}$ under the homomorphism $\SL_2\to G$ corresponding to the root $\gamma$, and $(\cdot,\cdot)_\frg$ is an $\Ad(G)$-invariant nondegenerate symmetric bilinear pairing on $\frg$ (which exists when char$(k)$ is good, see \cite[\S1.16]{Ca}). It is easy to check that $(\cdot,\cdot)_\Lambda$ is a nondegenerate symmetric bilinear pairing, and the $L_\gamma$-action on $\Lambda$ preserves this pairing up to scalar. Let $q_\Lambda$ be the quadratic form associated to $(\cdot,\cdot)_\Lambda$.

By a quick case-by-case analysis, one checks that the action map $L_\gamma\to\GO(\Lambda,q_\Lambda)$ is surjective. Therefore, $\Lambda$ contains a unique $L_\gamma$-stable irreducible divisor $Q_\Lambda=\{q_\Lambda=0\}$ whose complement is a single $L_\gamma$-orbit.

On the other hand, the complement $\Gr_{\leq\gav}-\Gr^\triv_{\leq\gav}$ is an ample divisor representing the class of the determinant line bundle, hence has codimension 1. Since $\Gr_{\leq\gav}-\Gr_{\gav}$ has codimension at least 2, $\Gr_{\gav}-\Gr^\triv_{\gav}$ also has codimension 1 in $\Gr_{\gav}$. Hence $\Lambda-\Lambda^\triv$ also has codimension 1 in $\Lambda$, therefore must be the irreducible divisor $Q_\Lambda$. This implies $\Lambda^\triv=\Lambda-Q_\Lambda$, which is a single $\Png$-orbit, therefore $\Gr^\triv_{\gav}$ is a single $G$-orbit.

By the same $\SL_2$-calculation as in \eqref{sl2}, we have
\begin{equation*}
u_\gamma(\tau)\tau^\gav=u_{-\gamma}(\tau^{-1})\in\Gr^\triv.
\end{equation*}
Therefore, $\Lambda^\triv\tau^\gav=\Ad(\Png)u_\gamma(\tau)\tau^\gav=\Ad(\Png)u_{-\gamma}(\tau^{-1})\subset\Gr^\triv$. By a similar argument as above, $\Ad(\Png)u_{-\gamma}(\tau^{-1})$ is the open subset of $\Vng=\prod_{\jiao{\alpha,\gav}=2}U_{-\alpha,-1}$ defined by the complement of a quadric $\Qng$. Therefore
\begin{equation*}
\Gr^\triv_{\gav}=G\twtimes{\Png}(\Lambda-Q_\Lambda)\tau^\gav=G\twtimes{\Png}(\Vng-\Qng)\subset G[\tau^{-1}]_1.
\end{equation*}
\end{proof}

\begin{lemme}\label{l:adjres}
The morphism \eqref{Gradj} extends to a resolution:
\begin{equation*}
\nu:G\twtimes{\Png}\Vng\to\Gr^\triv_{\leq\gav}
\end{equation*}
which is an isomorphism over $\Gr^\triv_\gav$ by Lemma \ref{l:adj}.
\begin{enumerate}
 \item The fiber $\nu^{-1}(\star)\cong G/\Png$;
\item The fibers over $\Gr^\triv_\thv$ are isomorphic to $ L_\theta/L_\theta\cap\Png$. 
\end{enumerate}
\end{lemme}
\begin{proof}
Since $G\twtimes{\Png}(\Vng-\Qng)$ is dense in $G\twtimes{\Png}\Vng$, and $\Gr^{\triv}_{\leq\gav}$ is the closure of $\Gr^\triv_\gav$ in $\Gr^\triv$, the morphism \eqref{Gradj} extends to $\nu$. By Lemma \ref{l:adj}, there are three $\Ad(\Png)$ (or $\Ad(L_\gamma)$)-orbits on $\Vng$: $\Vng-\Qng$, $\Qng^\times:=\Qng-\{0\}$ and $\{0\}$, which give three $G$-orbits of $G\twtimes{\Png}\Vng$.

The orbit $G\twtimes{\Png}\{0\}$ maps to $\Gr^\triv_0=\star$, which proves (1).

The orbit  $G\twtimes{\Png}\Qng^\times$ must then map to $\Gr^\triv_\thv$. The $G$-stabilizer of $u_{-\theta}(\tau^{-1})$ in $G\twtimes{\Png}\Qng^\times$ and in $\Gr^\triv_\thv=G\twtimes{\Pnth}U^\times_{-\theta,-1}$ are $\Png\cap\Pnth^1$ and $\Pnth^1$ respectively, where $\Png^1=\ker(\theta:\Pnth\to\GG_m)$. The fiber $\nu^{-1}(u_{-\theta}(\tau^{-1}))$ thus equals $\Pnth^1/\Pnth^1\cap\Png$. It is easy to check that the inclusions $L_\theta\hookrightarrow\Pnth$ and $\Pnth^1\hookrightarrow\Pnth$ induce isomorphisms
\begin{equation*}
L_\theta/L_\theta\cap\Png\isom\Pnth/\Pnth\cap\Png\xleftarrow{\sim}\Pnth^1/\Pnth^1\cap\Png.
\end{equation*}
Since $\Gr^\triv_\thv$ is a single $G$-orbit by Lemma \ref{l:qm}, all the fibers of $\nu$ over $\Gr^\triv_\thv$ are isomorphic to $\nu^{-1}(u_{-\theta}(\tau^{-1}))\cong L_\theta/L_\theta\cap\Png$.

\end{proof}

The Bruhat decomposition $G=\sqcup_{W/W_\gamma}Uw\Png$ gives a decomposition
\begin{equation*}
G\twtimes{\Png}\Vng=\bigsqcup_{W/W_\gamma}Uw\Png\twtimes{\Png}\Vng
\end{equation*}
The map $w\mapsto -w\gamma$ sets up a bijection between $W/W_{\gamma}$ and the set of short roots of $G$. For $\beta=-w\gamma$, let $V_\beta=\Ad(w)\Vng$. Then the above decomposition can be rewritten as
\begin{equation}\label{decompadj}
G\twtimes{\Png}\Vng=\bigsqcup_{\beta\textup{ short root}}\Ad(U)V_\beta.
\end{equation}
As in the quasi-minuscule case, we can further write each stratum as
\begin{eqnarray}\label{adjstratum}
\Ad(U)V_\beta &\cong & \prod_{\alpha\in\Phi^+,\jiao{\alpha,\beta^\vee}<0}U_\alpha\times\prod_{\jiao{\beta',\beta^\vee}=2}U_{\beta',-1}\\
\notag\Ad(\prod u_\alpha) \prod u_{\beta'}(c_{\beta'}\tau^{-1})&\leftrightarrow& (\prod u_\alpha, u_{\beta'}(c_{\beta'}\tau^{-1})).
\end{eqnarray}

\begin{remark}\label{r:Phibeta2}
(1) For each short root $\beta$, the set $\Phi^\beta_2$ (those roots which appear in the factors of $V_\beta$) is {\em totally ordered}  according to their heights. To see this, we only need to show that for different $\beta',\beta''\in\Phi^\beta_2$, $\jiao{\rho^\vee,\beta'-\beta''}\neq0$. In the proof of Lemma \ref{l:adj}, we remarked that $L_\gamma\to\GO(\Lambda,q_\Lambda)$ is surjective, which means that for any two different roots $\alpha',\alpha''\in\Phi^\gamma_2$ (i.e., they appear in the factors of $\Lambda$), the difference $\alpha'-\alpha''$ is a nonzero multiple of a root of $L_\gamma$. Since $\Phi^\beta_2=w\Phi^\gamma_2$ if $\beta=w\gamma$, the difference $\beta'-\beta''$ is also a nonzero multiple of a root of $L_\beta$, therefore $\jiao{\rho^\vee,\beta'-\beta''}\neq0$.

The set $\Phi^\beta_2$ carries an involution $\beta'\mapsto2\beta-\beta'$, with $\beta$ the only fixed point. This involution is order-reversing. 

(2) Again since $L_\gamma\twoheadrightarrow\GO(\Lambda,q_\Lambda)$ in the proof of Lemma \ref{l:adj}, all roots $\gamma'\in\Phi^\gamma_2\backslash\{\gamma\}$ are permuted by $W_\gamma$. Similar statement hold if  $\gamma$ is replaced by any short root $\beta$. In particular, all roots in $\Phi^\beta_2$ are long roots except $\beta$ itself.
\end{remark}

\begin{lemme}\label{l:adjf0}
For a short root $\beta$, $\prod u_{\beta'}(c_{\beta'}\tau^{-1})\in V_\beta$ (the product over $\Phi^\beta_2$), and $u\in U$, we have
\begin{equation*}
a_{-\theta}(\Ad(u)(\prod u_{\beta'}(c_{\beta'}\tau^{-1})))=\begin{cases} c_{-\theta}\bbx_{-\theta} &\textup{ if }-\theta\in\Phi^\beta_2\\
0 & \textup{otherwise.}\end{cases}
\end{equation*}
\end{lemme}
The proof is similar to that of Lemma \ref{l:qmf01}.

\subsection{Proof of Theorem \ref{th:Eulerchar}(2)}\label{ss:adjproof} Here we assume $G$ is of type $B$,$C$ or $F_4$. For $G$ of type $A$,$D$,$E$, the statement (2) is identical to (1) in Theorem \ref{th:Eulerchar}; the proof for $G=G_2$ will be given in \S\ref{a:g2}.

Applying Corollary \ref{c:chiKlV} to $V=\dualg$, which has two nonzero dominant weights $\thv$ and $\gav$, each with multiplicity one, we conclude
\begin{equation*}
\chi_c(\pline,\Kl^\Ad_\dualG(\phi,\chi))=\chi_c(\Gr^{\triv,a=0}_\thv,\ev^*J)+\chi_c(\Gr^{\triv,a=0}_\gav,\ev^*J).
\end{equation*}
We already know from (1) that $\chi_c(\Gr^{\triv,a=0}_\thv,\ev^*J)=r_\ell(G)$, we only need to show that $\chi_c(\Gr^{\triv,a=0}_\gav,\ev^*J)=r_s(G)$.

For any subset $S\subset\Phi$, let
\begin{equation*}
V^{S}_\beta=\{\prod u_{\beta''}(c_{\beta''}\tau^{-1})\in V_\beta|c_{\beta''}=0\textup{ for }\beta''\notin S\}.
\end{equation*}
Using this definition and the standard partial ordering on $\Phi$ (which restricts to the total ordering on $\Phi^\beta_2$, by Remark \ref{r:Phibeta2}(1)), the meanings of $V^{>-\theta}_\beta$, $V^{>0}_\beta$, etc.  are obvious.

Using the Lemma \ref{l:adjres}, Lemma \ref{l:adjf0} and the decomposition \eqref{decompadj}, the resolution $\nu$ restricted to $\Gr^{\triv,a=0}_{\leq\gav}$ reads:
\begin{equation*}
\nu^{a=0}:\bigsqcup_{\beta\textup{ short root}}\Ad(U)(V^{>-\theta}_\beta)\to\Gr^{\triv,a=0}_{\leq\gav}.
\end{equation*}
By Lemma \ref{l:adjres}, the Euler characteristic of $\ev^*J$ on the target of $\nu^{a=0}$ has contributions from $\star$, $\Gr^{\triv,a=0}_\thv$ and $\Gr^{\triv,a=0}_\gav$:
\begin{eqnarray}
\label{3strata}&&\sum_{\beta\textup{ short root}}\chi_c(\Ad(U)V^{>-\theta}_\beta,\ev^*J)\\
\notag&=&\chi_c(G/\Png)+\chi_c(\Gr^{\triv,a=0}_\thv,\ev^*J)\chi_c(L_\theta/L_\theta\cap\Png)+\chi_c(\Gr^{\triv,a=0}_\gav,\ev^*J).
\end{eqnarray}
The Bruhat-decomposition implies:
\begin{eqnarray}
\label{0fiber}\chi_c(G/\Png)=\#(W/W_\gamma)=\#\{\textup{short roots in }\Phi\};\\
\label{thvfiber}\chi_c(L_\theta/L_\theta\cap\Png)=\#(W_\theta/W_{\theta}\cap W_\gamma).
\end{eqnarray}

\begin{claim} For a short root $\beta$
\begin{equation*}
\chi_c(\Ad(U)(V^{>-\theta}_\beta),\ev^*J)=\begin{cases}
0 & \Phi^\beta_2\textup{ contains a simple root}\\
1 & \textup{otherwise.}\end{cases}
\end{equation*}
\end{claim}

Admitting the Claim first, we finish the proof. Combining \eqref{3strata},\eqref{0fiber}, \eqref{thvfiber} and the Claim, we get
\begin{eqnarray}\label{chinumbers}
&&-\chi_c(\Gr^{\triv,a=0}_\gav,\ev^*J)\\
\notag&=&\#\{\textup{short roots}\}-r_\ell(G)\#(W_\theta/W_\theta\cap W_\gamma)\\
\notag&&-\#\{\textup{short roots  }\beta|\Phi^\beta_2\textup{ does not contain simple root}\}\\
\notag&=&\#\{\textup{short roots  }\beta|\Phi^\beta_2\textup{ contains a simple root}\}-r_\ell(G)\#(W_\theta/W_\theta\cap W_\gamma).
\end{eqnarray}
By Remark \ref{r:Phibeta2}(1), $\Phi^\beta_2$ is totally ordered by heights, hence contains at most one simple root. Therefore
\begin{equation*}
\#\{\textup{short roots  }\beta|\Phi^\beta_2\textup{ contains a simple root}\}=N_s+N_\ell
\end{equation*}
where $N_s$ (resp. $N_\ell$) is the number of short roots $\beta$ such that $\Phi^\beta_2$ contains a short (resp. long) simple root.
\begin{itemize}
\item If $\Phi^\beta_2$ contains a short simple root, this simple root must be $\beta$ since $\beta$ is the only short root in $\Phi^\beta_2$ by Remark \ref{r:Phibeta2}(2). Therefore $N_s=r_s(G)$.
\item Any simple long root is in the $W$-orbit of $\theta$, therefore $N_\ell=r_\ell\cdot N_\theta$ where $N_\theta=\#\{\textup{short roots  }\beta|\theta\in\Phi^\beta_2\}$. Such short roots $\beta$ are in the $W_\theta$-orbit (this follows by applying Remark \ref{r:Phibeta2}(2) to the dual root system), and $\gamma$ is one of them, hence $N_\theta=\#(W_\theta/W_\theta\cap W_\gamma)$. Hence $N_\ell=r_\ell(G)\#(W_\theta/W_\theta\cap W_\gamma)$.
\end{itemize}
Combining these calculations and \eqref{chinumbers}, we conclude that $-\chi_c(\Gr^{\triv,a=0}_\gav,\ev^*J)=r_s(G)$, hence proving Theorem \ref{th:Eulerchar}(2).

It remains to prove the Claim.

\noindent{\bf Case I.} $\Phi^\beta_2$ contains a simple root $\alpha_i$. For $u\in U$ and $v=u_{\alpha_i}(c_{\alpha_i}\tau^{-1})v^{\alpha_i}\in V_\beta$, where $v^{\alpha_i}\in\prod_{\beta'\neq\alpha_i}U_{\beta',1}$, we have
\begin{equation*}
\ev(\Ad(u)v)=\Ad(u)u_{\alpha_i}(c_{\alpha_i})\cdot\Ad(u)\ev(v^{\alpha_i}).
\end{equation*}
This means that $\Ad(u)v\in U\barB$ if and only if $\Ad(u)\ev(v^{\alpha_i})\in U\barB$, and in case this happens, its image in $U/[U,U]$ is $u_{\alpha_i}(c_{\alpha_i})$ times the image of $\Ad(u)\ev(v^{\alpha_i})$ in $U/[U,U]$. Using the decomposition \eqref{adjstratum}, we can write $\Ad(U)V^{>-\theta}_\beta=U_{\alpha_i,-1}\times\bbA$ for some affine space $\bbA$ (the product of other factors), and the sheaf $\ev^*J$ is an exterior product $\ev^*\AS_{\phi_i}\boxtimes L$ for some local system $L$ on $\bbA$. By K\"unneth formula,
\begin{equation*}
H^*_c(\Ad(U)V^{>-\theta}_\beta,\ev^*J)=H^*_c(U_{-\alpha_i,-1},\ev^*\AS_{\phi_i})\otimes H^*_c(\bbA,L)=0.
\end{equation*}

\noindent{\bf Case II.} $\Phi^\beta_2$ does not contain any simple root.

We stratify the vector space $V^{>-\theta}_\beta$ into
\begin{equation*}
V^{>-\theta}_\beta=V^{>0}_\beta\bigsqcup\left(\bigsqcup_{\beta'\in\Phi^\beta_2,-\theta<\beta'<0}(V^{\geq\beta'}_\beta-V^{>\beta'}_\beta)\right).
\end{equation*}

First, we show that $\chi_c(\Ad(U)V^{>0}_\beta,\ev^*J)=1$. In fact, since $\Phi^\beta_2\cap\Phi^+$ contains no simple root,  $\ev(\Ad(U)V^{>0}_\beta)\subset[U,U]$, hence $\ev^*J$ is the constans sheaf on $\Ad(U)V^{>0}_\beta$. Since $\Ad(U)V^{>0}_\beta$ is an affine space by \eqref{adjstratum}, we get the conclusion.

Second, we prove that $\chi_c(\Ad(U)(V^{\geq\beta'}_\beta-V^{>\beta'}_\beta),\ev^*J)=0$ for each $\beta'\in\Phi^\beta_2,-\theta<\beta'<0$. 

Since $\beta'\neq-\theta$, $\beta'-\alpha_i$ is still a root for some simple $\alpha_i$. Then $\alpha=\alpha_i-\beta'$ is a positive root. By assumption, $\alpha_i\notin\Phi^\beta_2$, therefore $\jiao{\alpha,\beta^\vee}=\jiao{\alpha_i,\beta^\vee}-2<0$, i.e., $\alpha$ appears in the first product in \eqref{adjstratum}.

\begin{lemme}\label{l:largerbeta}
Let $a,b\in\ZZ_{>0}$ and $\beta'\leq\beta_1,\cdots,\beta_b$ be roots in $\Phi^\beta_2$ (not necessarily distinct),  then $a\alpha+\sum_{i=1}^b\beta_i\neq0$. If  $a\alpha+\sum_{i=1}^b\beta_i$ is a root, then one of the following situations happens
\begin{enumerate}
\item  $a\alpha+\sum_{i=1}^b\beta_i$ is a negative root. Then $a=b-1$, $a\alpha+\sum_{i=1}^{b}\beta_i\in\Phi^\beta_2$ and is larger than any of the $\beta_i$'s (in the total order of $\Phi^\beta_2$);
\item $a\alpha+\sum_{i=1}^b\beta_i$ is positive but not simple;
\item $a=b=1, \beta_1=\beta'$ and $a\alpha+\sum_{i=1}^b\beta_i=\alpha_i$.
\end{enumerate}
\end{lemme}
\begin{proof}
If $a\geq b$. Since $\alpha+\beta'>0$ and $\beta_i\geq\beta'$,  $\alpha+\beta_i$ must have positive height. Therefore, $$a\alpha+\sum_{i=1}^b\beta_i=(a-b)\alpha+\sum_{i=1}^b(\alpha+\beta_i)>0.$$ Since each term on RHS of the above sum has positive height, it is a simple root if there is only one summand, which must be the case (3). 

If $a<b$, then $\jiao{a\alpha+\sum_{i=1}^b\beta_i,\beta^\vee}\geq 2b-2a\geq2$. Therefore, if $a\alpha+\sum_{i=1}^b\beta_i$ is a root, we must have $a=b-1$ and $a\alpha+\sum_{i=1}^b\beta_i\in\Phi^\beta_2$. Since $\Phi^\beta_2$ does not contain simple roots, we get either case (1) or case (2). 

In any case, we have $a\alpha+\sum_{i=1}^b\beta_i\neq0$.
\end{proof}

For $u=\prod u_{\alpha'}(c_{\alpha'})$ (the product over $\alpha'\in\Phi^+,\jiao{\alpha,\beta^\vee}<0$), let $u^\alpha=\prod_{\alpha'\neq\alpha} u_{\alpha'}(c_{\alpha'})$. Similarly, for  $v\in V^{\geq\beta'}_\beta-V^{>\beta'}_\beta$, write $v=u_{\beta'}(c_{\beta'}\tau^{-1})v^{\beta'}\prod u_{\beta''}(c_{\beta''}\tau^{-1})$ (product over $\beta'<\beta''\in\Phi^\beta_2$, and $c_{\beta'}$ invertible). Then
\begin{equation}\label{commadj}
\ev(\Ad(u)v)=u^\alpha[u_\alpha(c_\alpha),u_{\beta'}(c_{\beta'})]u_{\beta'}(c_{\beta'})\left(\prod_{\beta''>\beta'}[u_{\alpha}(c_\alpha),u_{\beta''}(c_{\beta''})]u_{\beta''}(c_{\beta''})\right)u^{\alpha,-1}.
\end{equation}

To calculate \eqref{commadj}, we use Chevalley commutator relation to write each $[u_\alpha(c_\alpha),u_{\beta''}(c_{\beta''})]$ into products of positive and negative root factors, and try to pass the negative root factors (which necessarily appear in $\Phi^\beta_2$ by Lemma \ref{l:largerbeta}(1)) to the right. By Chevalley's commutator relation, each time we will produce new factors of the form $a\alpha+\beta_1+\cdots\beta_b$ for $\beta_i\in\Phi^\beta_2, a,b>0$. We keep the positive factors and pass the negative factors further to the right. In only thing we need to make sure in this process is when we do commutators $[U_{\alpha'}, U_{\gamma'}]$, we always have that $\alpha'$ and $\gamma'$ are linearly independent so that Chevalley's commutator relation is applicable. In fact, in the process, we only encounter the case where $\alpha'\in\Phi^+,\gamma'\in\Phi^-$ and $\alpha'+\gamma'$ has the form $a\alpha+\beta_1
+\cdots+\beta_b$ for $\beta'\leq\beta_i\in\Phi^\beta_2$ ($a,b\in\ZZ_{>0}$).  The only possibility for $\alpha'$ and $\gamma'$ to be linearly dependent is $\alpha'+\gamma'=0$, which was eliminated by Lemma \ref{l:largerbeta}.

In the end of the process, we get
\begin{equation}\label{final}
\ev(\Ad(u)v)=u^{\alpha}u_{\alpha_i}(\epsilon c_\alpha c_{\beta'}) u^+\left(\prod_{\beta''\geq\beta',\beta''\in\Phi^\beta_2} u_{\beta''}(c_{\beta''}+\tilc_{\beta''})\right) u^{\alpha,-1}.
\end{equation}
The term $u_{\alpha_i}(\epsilon c_\alpha c_{\beta'})$ (where $\epsilon=\pm1$) comes from $[u_\alpha(c_\alpha),u_{\beta'}(c_{\beta'})]$. The term $u^+$ is the product all the other {\em positive} factors in $U_{a\alpha+\beta_1+\cdots\beta_b}$ (i.e., $a\alpha+\beta_1+\cdots\beta_b\in\Phi^+$). By Lemma \ref{l:largerbeta}(2)(3), these $a\alpha+\beta_1+\cdots+\beta_b$ are never simple, therefore $u^+\in[U,U]$. Finally, the extra coefficient $\tilc_{\beta''}$ comes from the negative factors in $U_{a\alpha+\beta_1+\cdots\beta_b}$, which is a polynomial functions in $c_\alpha$ and $c_{\beta'''}$ for $\beta'''\in\Phi^\beta_2$. By Lemma \ref{l:largerbeta}(1), $c_{\beta''}$ only involves those $c_{\beta'''}$ such that $\beta'''<\beta''$.

Therefore we can make a change of variables
\begin{eqnarray}\label{changevar}
\Ad(U)(V^{\geq\beta'}_\beta-V^{>\beta'}_\beta)&\cong& U_{\alpha}\times U^\times_{\beta',-1}\times\prod_{\jiao{\alpha',\beta^\vee}<0,\alpha'\neq\alpha}U_{\alpha'}\times\prod_{\beta''\in\Phi^\beta_2,\beta''>\beta'}U_{\beta'',-1}\\
\notag\Ad(u)(v)&\leftrightarrow&(u_\alpha(c_\alpha),u_{\beta'}(c_{\beta'}\tau^{-1}),u^{\alpha},\prod u_{\beta''}(c_{\beta''}+\tilc_{\beta''})).
\end{eqnarray}
Let $\bbA$ be the product of the last three terms  in \eqref{changevar}. Let $\pr_\bbA:\Ad(U)(V^{\geq\beta'}_\beta-V^{>\beta'}_\beta)\to\bbA$ be the projection. In view of \eqref{final}, the restriction of $\ev^*J$ on the fibers of $\pr_\bbA$ are isomorphic to Artin-Schreier sheaves on $U_\alpha$ , therefore $\pr_{\bbA,!}\ev^*J=0$, hence
\begin{equation*}
H^*_c(\Ad(U)U^\times_{\beta,-1},\ev^*J)=H^*_c(\bbA,\pr_{\bbA,!}\ev^*J)=0.
\end{equation*}

This completes the proof of Case II of the Claim.


\section{Global monodromy}\label{s:mono}

This section is devoted to the proof of Theorem \ref{th:glob}. Let $G$ be split almost simple over $k$.

\subsection{Dependence on the additive character}


Recall from Remark \ref{RemarkHeckeDiagram}\eqref{sym} that $T\rtimes\Aut^\dagger(G)\times\grot$ acts on $\Bun_{G(0,2)}$ and the Hecke correspondence \eqref{HeckeOp}.

The group  $T(k)\rtimes\Aut^\dagger(G)\times\grot(k)$ also acts on $I_\infty(1)/I_\infty(2)$, hence on the space of generic additive characters. Let $S_\phi$ be the stabilizer of $\phi$ under the action of $T\rtimes\Aut^\dagger(G)\times\grot$. This is a finite group scheme over $k$.

When $G$ is of adjoint type, for each $\sigma\in\Aut^\dagger(G)$, there is a unique $(t,s)\in(T\times\grot)(k)$ such that $(t,\sigma,s)$ fixes $\phi$. Therefore, in this case, the projection $S_\phi\to\Aut^\dagger(G)\isom\Out(G)$ is an isomorphism (as discrete groups over $k$). In general, $S_\phi\to\Aut^\dagger(G)$ is a $ZG$-torsor, but may not be surjective on $k$-points.

The following lemma follows immediately from the definition of the geometric Hecke operators. 
\begin{lemme}\label{l:Klequiv}
The tensor functor defining $\Kl_\dualG(\phi)$
\begin{equation*}
\Hk_\phi:\cS\ni\IC_V\mapsto\Hk_V(A_\phi)|_{\star\times\pline}\in\Loc(\pline)
\end{equation*}
carries a natural $S_\phi$-equivariant structure. Here $S_\phi(k)$ acts on $\cS$ via its image in $\Aut^\dagger(G)$ and 
$S_\phi(k)$ acts on $\Loc(\pline)$ via its action on $\pline$ through $S_\phi\to\grot$.
\end{lemme}
In Lemma \ref{l:outSatakeH}(1) we will check that under the equivalence $\Rep(\dualG) \cong \cS$ the $\Aut^\dagger(G)$-action on $\cS$ coincides with the action of $\Aut^\dagger(G)\cong \Aut^\dagger(\dualG)$ with respect to the chosen pinning of the dual group. 

Let $S^1_\phi(k)=\ker(S_\phi(k)\to\grot(k))$ and let $\Srot_\phi(k)$ be the image of $S_\phi(k)$ in $\grot(k)$, which is a finite cyclic group. We would like to use the above equivariance to conclude that the monodromy representation of $\Kl_\dualG(\phi)$ can be chosen to take values in $\dualG^{S_\phi^1(k)}$. However, since this representation is only  defined up to inner automorphisms, this requires an extra argument. 

Recall that we have chosen a geometric generic point $\generic$ over $\Spec K_0$, the formal punctured disc at $0$. We defined $\Kl_\dualG(\phi,1)$ as a functor from the Satake category $\cS\cong \Rep(\dualG)$ to $\Loc(\pline)$, so that  restriction to $\generic$ is a tensor functor
\begin{equation}
\omega_\phi:\cS\xrightarrow{\Kl_\dualG(\phi)}\Loc(\pline)\xrightarrow{j^*_\generic}\Vect,
\end{equation}
As this is a fiber functor of $\cS$, we can define $\dG=\uAut^\otimes(\omega_\phi)$. Of course $\dG$ is isomorphic to $\dualG$ but not canonically so.
The group $\dG$ has the advantage, that by Tannakian formalism, we get a canonical homomorphism $\uAut^\otimes(j^*_\generic)\to\uAut^\otimes(\omega_\phi)$, where $\uAut^\otimes(j^*_\generic)$ is the pro-algebraic envelope of the fundamental group $\pi_1(\pline,\generic)$. Therefore we get a homomorphism
\begin{equation}\label{globalrep}
\varphi:\pi_1(\pline,\generic)\to\dG(\Ql).
\end{equation}


By Lemma \ref{l:Klequiv} we get a homomorphism
\begin{equation}\label{S1actondG}
S^1_\phi(k)\to\Aut^\otimes(\cS,\Hk_\phi)\xrightarrow{(j^*_\generic)_*}\Aut^\otimes(\cS,\omega_\phi)=\Aut(\dG).
\end{equation}
In other words, we have an action of $S^1_\phi$ on $\dG$. We will prove in Lemma \ref{l:outSatake} that $\dG$ also carries a natural pinning $\ddagger$ and that the above action preserves this pinning. The pinnings $\ddagger$ and $\dagger$ define a canonical isomorphism $\canon\colon\dG\cong \dualG$. Using this identification we obtain:
\begin{corollary}\label{c:outfix}
The monodromy representation $\varphi$ extends to a homomorphism between exact sequences:
\begin{equation}\label{extmonorep}
\xymatrix{\pi_1(\pline,\generic)\ar[r]^{[\#\Srot_\phi(k)]}\ar[d]^{\varphi} & \pi_1(\pline,\generic)\ar[r]\ar[d]^{\widetilde{\varphi}} & \Srot_\phi(k)\ar@{=}[d]\\
\dualG^{S^1_\phi(k)}\ar[r] & \dualG^{S^1_\phi(k)}\rtimes\Srot_\phi(k)\ar[r] & \Srot_\phi(k)}
\end{equation}
\end{corollary}
\begin{proof}
By \eqref{S1actondG}, the $S^1_\phi(k)$-action on $\dG$ factors through $S^1_\phi(k)\to\Aut^\otimes(\cS,\Hk_\phi)$. By Lemma \ref{l:inv}, the monodromy representation $\varphi$ thus factors through $\pi_1(\pline,\generic)\to\dG^{S^1_\phi}\subset\dG$, i.e.,
\begin{equation*}
\Hk_\phi:\cS\cong\Rep(\dG)\xrightarrow{\Res}\Rep(\dG^{S^1_\phi(k)})\xrightarrow{\kappa}\Loc(\pline)
\end{equation*}
for some tensor functor $\kappa$. 

We identified $\dualG\cong\dG$ using the pinned isomorphism ``$\canon$'', so we can view the functor $\kappa$ as $\Rep(\dualG^{S^1_\phi(k)})\to\Loc(\pline)$. Then $\Srot_\phi(k)=S_\phi(k)/S^1_\phi(k)$ acts on $\dualG^{S^1_\phi(k)}$ via $S_\phi(k)\to\Out(G)\xrightarrow{\iota}\Aut^\dagger(\dualG)$, hence acting on the source and target of $\kappa$. Lemma \ref{l:Klequiv} implies that $\kappa$ carries a natural $\Srot_\phi(k)$-equivariant structure. Taking $\Srot_\phi(k)$-invariants of both tensor categories, we get a functor $\widetilde{\kappa}$:
\begin{eqnarray*}
\xymatrix{\Rep(\dualG^{S^1_\phi(k)}\rtimes\Srot_\phi(k))\ar[d]^{\wr}\ar[r]^{\widetilde{\kappa}} & \Loc(\pline/\Srot_\phi(k))\ar[d]^{\wr}\\ 
\Rep(\dualG^{S^1_\phi(k)})^{\Srot_\phi(k)}\ar[r]^{(\kappa)^{\Srot_\phi(k)}} &\Loc(\pline)^{\Srot_\phi(k)}}
\end{eqnarray*}
Since the quotient map $\pline\to\pline/\Srot_\phi(k)$ can be identified with the $\#\Srot_\phi(k)$-th power map of $\pline=\Gm$, we arrive at the diagram \eqref{extmonorep}.
\end{proof}

\begin{remark}
Corollary \ref{c:outfix} remains true if $k$ is replaced by an extension $k'$ and $\pline$ is replaced by $\pline\otimes_kk'$. In particular, to get information about the geometric monodromy, we take for $k'=\bark$.

In the case $\Out(\dualG)$ is nontrivial, we have
\begin{table}[h]
\caption{Outer automorphisms and stabilizers of $\phi$}
\begin{tabular}{l|l|l|l|l}
$\dualG$ & $\Aut^\dagger(\dualG)\isom\Out(\dualG)$ & $\dualG^{\Aut^\dagger(\dualG)}$ & $S^1_\phi(\bark)$ & $\Srot_\phi(\bark)$ \\\hline
 $A_{2n-1}$ $(n\geq2)$& $\ZZ/2$ & $C_n$ & $\Out(\dualG)$ & 1\\ 
$A_{2n}$ & $\ZZ/2$ & $B_{n}$ & 1 & $\ZZ/2$\\ 
$D_4$ & $S_3$ & $G_2$ & $\Out(\dualG)$ & 1\\ 
$D_{n}$ $(n\geq5)$ & $\ZZ/2$ & $B_{n-1}$ & $\Out(\dualG)$ & 1\\ 
$E_6$ & $\ZZ/2$ & $F_4$ & $\Out(\dualG)$ & 1
\end{tabular}
\label{t:outer}
\end{table}
\end{remark}

\subsection{Zariski closure of global monodromy}\label{ss:glob}
Let $\dualG_\geom\subset\dualG$ be the Zariski closure of the image of the geometric monodromy representation $\varphi^\geom:\pi_1(\geompline,\generic)\to\dG\cong\dualG$. We first show that $\dualG_\geom$ is not too small.

\begin{proposition}[B. Gross]\label{p:notPGL2}
If $\textup{char}(k)>2$, and the rank of $G$ is at least 2, the $\dualG_\geom$ is {\em not} contained in any principal $\PGL_2\subset\dualG$.
\end{proposition}
\begin{proof}
Suppose instead $\dualG_\geom\subset\PGL_2\subset\dualG$, where $\PGL_2$ contains a principal unipotent element (image of $\inert^t_0$) of $\dualG$.  The image of the wild inertia $\iinf^+$ must be nontrivial because $\Kl^\thv_\dualG(\phi)$ has nonzero Swan conductor at $\infty$. 

Since $p=\textup{char}(k)>2$, $\varphi(\iinf^+)$ lies in a maximal torus $\GG_m\subset\PGL_2$ and contains $\mu_p\subset\GG_m$. Since  $\varphi(\iinf)$ normalizes $\varphi(\iinf^+)$,  it must be contained in the normalizer $N(\GG_m)\subset\PGL_2$ of the torus $\GG_m$.

For any irreducible representation $S^{2\ell}=\Sym^{2\ell}(\St)$ of $\PGL_2$ (where $\St$ is the 2-dimensional representation of $\SL_2$), every pair of weight spaces $S^{2\ell}(n)\oplus S^{2\ell}(-n)$ ($0\leq n\leq\ell $) is stable under $N(\GG_m)$, hence under $\iinf$. If the weight $n$ does not divide $p$, the Swan conductor of $S^{2\ell}(n)\oplus S^{2\ell}(-n)$ is at least 1. Therefore
\begin{equation*}
\Swan_{\iinf}(S^{2\ell})\geq\ell-[\ell/p].
\end{equation*}

Consider the action of $\iinf$ on the quasi-minuscule representation $V_\thv$ of $\dualG$. By Lemma \ref{l:samenumber}, $V_\thv$ decomposes into $r_s(\dualG)$ irreducible representations of the principal $\PGL_2$: $V_\thv=S^{2\ell_1}\oplus\cdots\oplus S^{2\ell_{r_s(\dualG)}}$.  Therefore
\begin{equation*}
\Swan_{\iinf}(V_\thv)\geq\sum_{i=1}^{r_s(\dualG)}(\ell_i-[\ell_i/p])\geq(1-\frac{1}{p})\frac{hr_s(\dualG)}{2}.
\end{equation*}
Here we used $\sum_i\ell_i=hr_s(\dualG)/2$ (see the proof of Lemma \ref{l:samenumber}). So as long as $h>3$ or $p>3$, we get contradiction. Even when $h=3$ and $p=3$, then $\dualG=\PGL_3$ with $V_\thv=S^4\oplus S^2$, we still have $\Swan(V_\thv)\geq 2+1=3>r_s(\dualG)$, contradiction!
\end{proof}


\begin{proof}[Proof of Theorem \ref{th:glob}]
Since $\Kl^V_\dualG(\phi)$ is a pure of weight 0 for every $V\in\Rep(\dualG)$, hence geometrically semisimple. By Deligne \cite[Corollaire 1.3.9]{WeilII}, the neutral component $\dualG_\geom^\circ$ of $\dualG_\geom$ is a semisimple group.

\noindent{\bf Step I}. Assume $G$ is not of type $A_1$, $A_{2n}$ or $B_3$. We first determine $\dualG_\geom^\circ$, or equivalently, its Lie algebra $\dualg_\geom$.

On the one hand, by Theorem \ref{Heckeeigensheaf}(2), $\dualG_\geom^\circ$ contains a principal unipotent element, hence containing a principal $\PGL_2$. Proposition \ref{p:notPGL2} says that $\dualG_\geom$ cannot be equal to $\PGL_2$. Since a principal $\PGL_2$ is its own normalizer in $\dualG$, we conclude that $\dualG_\geom^\circ$ cannot be equal to $\PGL_2$, i.e., $\frsl_2\subsetneqq\dualg_\geom$.

On the other hand, since $G$ is not of type $A_{2n}$, $S^1_\phi(\bark)=\Aut^\dagger(\dualG)\isom\Out(\dualG)$. 
In this case, Corollary \ref{c:outfix} implies $\dualG_\geom\subset\dualG^{\Aut^\dagger(\dualG)}$.

Dynkin classified all Lie subalgebras $\dualg_\geom\subset\dualg$ which contain a principal $\frsl_2$. If $\dualg\neq\so_7$, then either $\dualg_\geom=\frsl_2$ or $\dualg_\geom$ is the fixed point algebra of some pinned automorphism of $\dualG$. In our case, if $G$ is not of type $A_{2n}$, we can already  conclude that $\dualg_\geom=\dualg^{\Aut^\dagger(\dualG)}$.

\noindent{\bf Step II}. Suppose $\dualG$ is of type $B_3$ and $\textup{char}(k)>3$, we claim that $\dualG_\geom^\circ=G_2$.

For this it suffices to show that $\Kl^{\wedge^3V_7}_{\SO_7}(\phi)$ contains has a global section over $\geompline$, where $V_7$ is the 7-dimensional standard representation of $\SO_7$. Suppose the contrary, then the long exact sequence \eqref{long1}--\eqref{long3} with $\dualg$ replaced by $\wedge^3V_7$ would imply
\begin{equation}\label{ineq}
\dim(\wedge^3V_7)^{\inert_0}+\dim(\wedge^3V_7)^{\iinf}\leq\dim H^1_c(\pline,\Kl^{\wedge^3V_7}_{\SO_7}(\phi))=\Swan_{\iinf}(\wedge^3V_7).
\end{equation}
We use the standard basis $\{e_{-3},\cdots,e_0,\cdots,e_3\}$ for $V_7$ (the quadratic form is $e_0^2+e_1e_{-1}+e_2e_{-2}+e_3e_{-3}$). The 0-weight space of $\wedge^3V_7$ under the principal $\frsl_2$ is spanned by $\{e_1\wedge e_{-2}\wedge e_{-3},e_{-1}\wedge e_2\wedge e_3, e_0\wedge e_i\wedge e_{-i}, i=1,2,3\}$, hence has dimension 5. Since $\inert_0$ acts on $V_7$ as a principal unipotent element by Theorem \ref{Heckeeigensheaf}(3), one concludes $\dim(\wedge^3V_7)^{\inert_0}=5$. 

Since $\textup{char}(k)>3$, it does not divide the Coxeter number of $\dualG$, therefore Corollary \ref{c:simple-wild} is applicable. Since the breaks of the $\iinf^+$-action on $V_7$ are $1/h=1/6$, the breaks of $\iinf^+$ on the nonzero weight spaces of $\wedge^3V_7$ are $\leq1/6$. Therefore $\Swan_{\iinf}(\wedge^3V_7)\leq[\frac{1}{6}\dim\wedge^3V_7]=[35/6]=5$. Moreover, by the description of $\varphi(\iinf^+)$ in Corollary \ref{c:simple-wild}, $\iinf^+$ acts trivially on $\Span\{e_1\wedge e_2\wedge e_3,e_{-1}\wedge e_{-2}\wedge e_{-3}\}$,  and the Coxeter permutation $e_1\to e_2\to e_3\to e_{-1}\to e_{-2}\to e_{-3}\to e_1$ permutes $e_1\wedge e_2\wedge e_3$ and $e_{-1}\wedge e_{-2}\wedge e_{-3}$. We get $\dim(\wedge^3V_7)^{\iinf}\geq1$. Thus
\begin{equation*}
\dim(\wedge^3V_7)^{\inert_0}+\dim(\wedge^3V_7)^{\iinf}\geq 5+1>5\geq\Swan_{\iinf}(\wedge^3V_7),
\end{equation*}
 which contradicts \eqref{ineq}. This proves that $\dualG_\geom=G_2$ in the case $\dualG=\SO_7$.

\noindent{\bf Step III}. When $\dualG$ is of type $A_{2n}$. In this case, $\Kl_\dualG(\phi)$ comes from the classical Kloosterman sheaf $\Kl_n(\phi)$, whose global monodromy is treated by Katz in \cite{KatzKloosterman}, see \eqref{Klnimage}.

\noindent{\bf Step IV}. It remains to prove that $\dualG_\geom$ is connected.

The case of $A_{n}$ is treated by Katz in \lc. In the case $\dualG$ is of type $E_6$, one checks that $\dualG^{\Aut^\dagger(\dualG)}$ is already connected. 

In general, we have $\dualG_\geom\subset N_\dualG(\dualG_\geom^\circ)=\dualG_\geom^\circ Z\dualG$. Consider the surjective homomorphism
\begin{equation*}
\overline{\varphi^\geom}:\pi_1^\geom(\pline,\generic)\xrightarrow{\varphi^\geom}\dualG_\geom\to\pi_0(\dualG_\geom)=Z\dualG/Z\dualG\cap\dualG_\geom^\circ.
\end{equation*}
Assuming $\dualG$ is not of type $A_n$ or $E_6$, then $\pi_0(\dualG)$ is a 2-group. Since char$(k)>2$, $\overline{\varphi^\geom}$ factors through the tame quotient. On the other hand, the tame generator in $\inert_0\cong\pi^\geom(\pline,\generic)$ must map to a unipotent element in $\dualG_\geom$, hence inside $\dualG_\geom^\circ$, therefore the map $\overline{\varphi^\geom}$ is trivial, i.e., $\dualG_\geom$ is connected.
\end{proof}

\section{Functoriality of Kloosterman sheaves--conjectures}\label{s:fonc}


In this section, we offer some conjectures for a further study on Kloosterman sheaves. In particular, according to a rigidity property that is known in the case of $\GL_n$ according to Katz and Gabber, the Kloosterman sheaves with the same geometric monodromy tabulated in Table \ref{t:globalmono} should be isomorphic after matching the additive characters $\phi$. We will also give a conjectural description of the local and global monodromy of $\Kl_{\LG}(\phi)$ for certain quasi-split groups $\cG$. 

\subsection{Rigidity}

\begin{conjecture}[Physical rigidity of Kloosterman sheaves]\label{conj:weak}
Suppose $L$ is a $\dualG$-local system on $\geompline$ which, as $\inert_0$ and $\iinf$-representations, has the same isomorphism types as $\Kl_\dualG(\phi,\chi)$. Then $L\cong\Kl_\dualG(\phi,\chi)$ over $\geompline$.
\end{conjecture}

Even stronger, we expect

\begin{conjecture}\label{conj:strong}
Suppose $L$ is a $\dualG$-local system on $\pline$ satisfying
\begin{itemize}
\item $L$ is tame at $\{0\}$, the semisimple part of the image of a topological generator of $\inert^t_0$ is conjugate to an element in $\dualT[q-1]$ which corresponds to a multiplicative character $\chi:T(k)\to\Qlt$. 
\item $\Swan_\infty(L^{\Ad})=r$, and $(L^{\Ad})^{\iinf}=0$.
\end{itemize}
Then there exists a generic linear function $\phi:I(1)/I(2)\to\GG_a$ such that $L\cong\Kl_\dualG(\phi,\chi)$ up to an unramified twist (given by $\Gal(\bark/k)\to Z\dualG$).
\end{conjecture}

Inspired by the above conjectures, and the calculation of the local and global monodromy of Kloosterman sheaves given in Theorem \ref{Heckeeigensheaf}(2), Corollary \ref{c:simple-wild} and Theorem \ref{th:glob}, we  conjecture that there should be functorial relationship between the Kloosterman sheaves. 

\begin{conjecture}\label{conj:Hasse-Davenport}
Let $G, G'$ be split almost simple groups over $k$ whose dual groups $\dualG \supset \dualG'$ appear in the same line of Table \ref{t:globalmono}. Then for every generic linear function $\phi$ of $G$, there exists a generic linear function $\phi'$ of $G'$ such that over $\geompline$, the Kloosterman sheaf $\Kl_{\dualG}(\phi)$ is the pushout of the Kloosterman sheaf $\Kl_{\dualG'}(\phi')$.
\end{conjecture}

\subsection{Quasi-split groups}\label{ss:conjqs}
Let $G$ be a split, almost simple and simply-connected group over $k$. Let $\cG$ be the quasi-split group scheme on $\PP^1$ whose restriction to $\pline$ is given by the twisting $\sigma:\mu_N\hookrightarrow\Aut^\dagger(G)$ as in \S\ref{ss:qs}. Recall that $N$ is assumed to be prime to $\textup{char}(k)$, and $N=2$ unless $G$ is of type $D_4$, in which case $N=3$. We abuse the notation to denote still by $\sigma$ the image of a generator of $\mu_N$ in $\Aut^\dagger(G)$. We identify $\Aut^\dagger(G)$ with $\Aut^\dagger(\dualG)$ using the isomorphism in Lemma \ref{l:outSatakeH}(2). We write $\jiao{\sigma}\subset\Aut^\dagger(G)\cong\Aut^\dagger(\dualG)$ for the subgroup generated by $\sigma$. The associated $L$-group can be taken as $\LG=\dualG\rtimes\jiao{\sigma}$. The Kloosterman sheaf $\Kl_\LG(\phi)$ constructed in Theorem \ref{Heckeeigensheaf} gives a monodromy representation
\begin{equation*}
\varphi:\pi_1(\pline,\generic)\to\LG=\dualG\rtimes\jiao{\sigma}.
\end{equation*}
The adjoint representation ${\rm Ad}:\dualG \to {\rm GL}(\dualg)$ can be extend to a homomorphism
$$ {\rm Ad}_\sigma: \LG\to {\rm GL}(\dualg),$$
so that the Kloosterman sheaf $\Kl_\LG(\phi)$ induces a local system $\Kl_\LG^{\Ad}(\phi)$.

We have the following predictions for the local monodromy of $\Kl_\LG(\phi)$. The tame monodromy at $\{0\}$ should be generated by the element $(u,\sigma)\in\LG=\dualG\rtimes\jiao{\sigma}$, where $u\in\dualG^\sigma$ is a principal unipotent element. The Swan conductor $\Swan_\infty(\Kl_{\LG}(\phi))=-\chi_c(\pline,\Kl_\LG^{\Ad}(\phi))$ should equal the rank of the neutral component of $\dualG^{\sigma}$. Moreover, we also expect that the analog of Corollary \ref{c:simple-wild} holds for $\Kl_\LG(\phi)$, with the Coxeter element replaced by the $\sigma$-twisted Coxeter element $\Cox_\sigma\in W\times\{\sigma\}\subset W\rtimes\jiao{\sigma}$ (see \cite{Springer} and \cite[\S 5]{Reeder}), and the Coxeter number replaced by the $\sigma$-twisted Coxeter number $h_\sigma$ (the order of the $\Cox_\sigma$).

As for the global monodromy, we expect that for  $\textup{char}(k)$ not too small, the global geometric monodromy representation $\varphi^\geom$ for $\Kl_\LG(\phi)$ has Zariski dense image.

We observe from Table \ref{t:globalmono} that for simply-laced split groups $G$ not of type $A_{2n}$, the Zariski closure $\dualG_\geom$ of the geometric monodromy of $\Kl_\dualG(\phi)$ is smaller than $\dualG$. Now pick a quasi-split form $\cG$ of $G$ built out of a nontrivial $\sigma\in\Aut^\dagger(G)$ of order $N$, then according to our expectation, the $\dualG$-local system $[N]^*\Kl_\LG(\phi)$ on $\tpline$ (the $N$-th Kummer cover of $\pline$, which is still isomorphic to $\pline$) should have Zariski dense geometric monodromy in $\dualG$. This compensates the smallness of $\dualG_\geom$ for simply-laced $G$.

When $G$ is of type $A_{2n}$, by Corollary\ref{c:outfix}, $\varphi$ extends to
\begin{equation*}
\pi_1(\pline,\generic)\to\dualG\rtimes\mu_2,
\end{equation*}
hence giving a $\dualG\rtimes\mu_2$-local system $\overline{\Kl}_\dualG(\phi)$ on $\pline/\mu_2$. On the other hand, let $\cG$ be the quasi-split unitary group given by the nontrivial $\sigma\in\Aut^\dagger(G)$. We have a Kloosterman sheaf $\Kl_\LG(\phi')$ on $\pline$. After identifying $\LG$ with $\dualG\rtimes\mu_2$ and matching $\phi$ with $\phi'$, we expect  that $\Kl_\LG(\phi)\cong\overline{\Kl}_\dualG(\phi')$ over $\geompline\cong(\pline/\mu_2)\otimes_k\bark$.

\section{Appendix: Proof of Proposition \ref{P1-uniformization} on the geometry of moduli spaces of $\cG$-bundles on $\bP^1$}\label{a:moduli}

In this appendix, we give a proof of Proposition \ref{P1-uniformization}.
\begin{proof}
The proof of the first two claims uses the same argument as in \cite[Proposition 3]{HainesRapoport}. We will first assume that the base field $k$ is algebraically closed.

A. Assume that $\cG=\cT$ is a torus. In this case (2) has been proved in \cite[Lemma 16]{Uniformization}.

Let us prove (1). Since $\cT$ is abelian, for any $\cT$-torsor $\cP$ we have $H^1(\bP^1,\cP\times^\cT \Lie(\cT)) = H^1(\bP^1, \Lie(\cT))$. We claim that this group vanishes. 
By construction we know that $\Lie(\cT) \subset \Lie([N_*](\bG_m^r\times \bP^1))$ is a direct summand. 
Also $\Lie([N]_*(\bG_m^r\times \bP^1))=[N]_*\cO_{\bP^1}^r$ and $H^1(\bP^1,[N]_*(\cO_{\bP^1}))=H^1(\bP^1,\cO_{\bP^1})=0$. Thus $H^1(\bP^1,\cP\times^\cT \Lie(\cT))=0$ for any $\cP$, so that there are no non-trivial deformations of $\cT$-bundles. Since $k$ is algebraically closed $H^1(k(\bP^1),\cT)=0$, i.e., all $\cT$-bundles are generically trivial. Thus the map $\oplus_{x\in \bP^1} \cT(K_x)/\cT(\cO_x) \to \Bun_{\cT}(k)$ is surjective. By construction, the isomorphism $\pi_0(\Bun_{\cT}) \cong \pi_1(\cT)_{\pi_1(\bG_m)}$ is induced from the Kottwitz homomorphism $\pi_0(\Gr_{\cT,x})\to \xcoch(\cT)_{\Gal(K_x^{sep}/K_x)}$. Thus for any $x$ the map $\cT(K_x)\to \Bun_{\cT}(k)$ is surjective. 
This proves (1) in the case of tori. Finally (1) and (2) imply (4) for tori.

B. For semisimple, simply connected groups $\cG$ statements (1) and (2) have been proved in \cite[ Theorem 4]{Uniformization}, again assuming $k$ to be algebraically closed. 

C. Suppose that the derived group $\cG_{\der}|\pline$ of $\cG|\pline$ is simply connected and denote $\cD|\pline:=\cG/\cG_{\der}|\pline$. Then Haines and Rapoport show that there is an induced exact sequence of group schemes over $\bP^1$:
$$ 1 \to \cG_{\der}\to \cG \to \cD \to 1.$$
Since any $\cD$-bundle is trivial over $\bP^1\backslash\{x\}$ and we know the same for all forms of $\cG_{der}$ as well, this implies (1) for $\cG$. Moreover, we know from (loc.cit.) that $\cG(k((s)))\to \cD(k((s)))$ is surjective. Thus $\pi_0(\Bun_\cG)\cong \pi_0(\Bun_\cD)$, which proves (2) for $\cG$.

D. For general groups $\cG$ we can choose a z-extension $$1 \to \cZ \to \cG^\prime \to \cG \to 1,$$ such that $\cZ$ is an induced torus and $\cG_{der}^\prime$ is simply connected. By Tsen's theorem $H^2(\bP^1,\cZ)=0$, thus any $\cG$ bundle can be lifted to a $\cG^\prime$ bundle. Since we know (1) for $\cG^\prime$ this proves (1) for $\cG$. To prove (2) we use that under our assumption $(\xcoch(\cZ))_{\pi_1(\pline)}$ is torsion free and therefore as in \cite{HainesRapoport} we have an exact sequence:
$$0 \to \xch(\check{\cZ}^{\pi_1(\pline)}) \to \xch(Z(\check{\cG}^{\prime})^{\pi_1(\pline)}) \to  \xch(Z(\check{\cG})^{\pi_1(\pline)}) \to 0.$$

To prove the third claim we follow the arguments of Ramanathan \cite{Ramanathan} or Faltings \cite[ Lemma 4]{Faltings}. 

We need to fix notations for dominant weights. First $\xch(\cT)^{\pi_1(\bG_m)}_{\bQ}=\xch(S)_{\bQ}$. Moreover, the relative roots $\Phi(\cG,S)$ span the subspace of characters of $S$ that are trivial on the center of $\cG$. 
We denote by $\xcoch(\cT)_{\pi_1(\bG_m),+}\subset \pi_0(\Bun_\cT)$ the the subset of $\gamma\in \xcoch(\cT)_{\pi_1(\bG_m)}$ such that for any positive, relative root $a\in \Phi(\cG,S)^+$ we have $a(\gamma)\geq 0$.

First, we want to prove that 
$$\cG(k((s))) = \coprod_{\lambda \in \xcoch(\cT)_{\pi_1(\bG_m),+}} \cG(k[t])\lambda(s)\cG(k[[s]]).$$
From (1) we conclude that every $\cG$-bundle is trivial outside any point. In particular any $\cG$ bundle admits a reduction to the Borel subgroup $\cB$. Let $\cE$ be a $\cG$-bundle and choose a reduction $\cE_\cB$ of $\cE$ to $\cB$. 

For any character $\alpha\colon\cT\to \bG_m$ we denote by $\cE_{\cB}(\alpha)$ the associated line bundle on $\bP^1$. Since $\xch(\cT)^{\pi_1(\bG_m)}_\bQ \cong \xch(S)_\bQ$ the degree $\deg(\cE_\cB(\alpha))\in \bQ$ is also defined for $\alpha \in \xch(S)$.

We claim that if for all positive, simple roots $a_i\in \Phi(\cG,S)^+$ we have 
 $\deg(\cE_{\cB}(a_i))\geq 0$ then the bundle $\cE_\cB$ admits a reduction to $\cT$. 

To show this, denote by $\cU$ the unipotent radical of $\cB$ and $\cE_{\cT}:=\cE_{\cB}/\cU$ the induced $\cT$-bundle. 
In order to show that $\cE_\cB$ is induced from $\cE_\cT$ we only need to show that $H^1(\bP^1,\cE_{\cT}\times^{\cT}\cU)=0$. The group $\cU$ has a filtration, such that the subquotients are given by root subgroups. 

Consider a positive, relative root $a\in \Phi(\cG,S)$. The root subgroup $\cU_a$ is a direct summand of $[N]_*(\oplus \cU_{\alpha^\prime})$ where the sum is over those roots $\alpha^\prime\in \Phi(G,T)$ that restrict to $a$ on $S$. Thus $\cU_a$ is a direct summand of a vector bundle $\cV$ satisfying $H^1(\bP^1,\cV)=0$. 


Similarly $\cE_{\cT}\times^{\cT} \cU_a$ is a direct summand of $\cE_{\cT}\times^{\cT} [N]_*(\oplus\cU_{\alpha^\prime})$. Since $[N]^*\cE_\cT$ is $\pi_1(\bG_m)$-invariant this implies that $H^1$ of this bundle is $0$ if $\deg(\cE_{\cB}(a))\geq 0$. Thus we also find $H^1(\bP^1,\cE_{\cT}\times^{\cT}\cU)=0$. 

If the reduction $\cE_\cB$ does not satisfy the condition $\deg(\cE_{\cB}(a_i))\geq 0$ for some simple root $a_i$ we want to modify the reduction $\cE_\cB$. Consider the parabolic subgroup $\cP_{a_i}\subset \cG$ generated by $\cB$ and $\cU_{-a_i}$. The root subgroups $\cU_{\pm a_i}$ define a subgroup $\cL$ of $\cP$ such that the simply connected cover of $\cL$ is either isomorphic to $[n]_*\SL_2$ or isomorphic to $[n]_*\SU_3$ for some $n$ dividing $N$ \cite[\S 4.1.4]{BruhatTits2}.  The semisimple quotient $\cP^{ss}/Z(\cP^{ss})$ is isomorphic to $\cL^{ad}$ and furthermore $\cP/\cB \cong \cL/(\cL\cap\cB)$. We claim that the result holds for $\cL^{ad}$ bundles, because these bundles can be described in terms of vector bundles. If $\cL^{ad}=[n]_*\PGL_2$, then the result follows from the result for vector bundles of rank $2$. The case of unitary groups is similar, we will explain it below in Lemma \ref{SU3Lemma}.

Thus we can find a new reduction $\cE^\prime_\cB$ to $\cB$ such that $\cE_{\cP_{a_i}}$ is unchanged, but $\cE_\cB^\prime(a_i)\geq 0$. This implies that for the fundamental weights $\epsilon_k$ (multiples of the determinant of the adjoint representation of the maximal parabolic subgroups $\cP^k$, generated by $\cB$ and all $\cU_{a_j}$ with $j\neq k$), the degree of $\cE_{\cB}(\epsilon_k)$ for $k\neq i$ is unchanged but the degree $\cE_{\cB}^\prime(\epsilon_i)$ is larger than the degree of $\cE_\cB(\epsilon_i)$. 

However $\cE_{\cB}^\prime(\epsilon_i)$ is a subbundle of $\wedge^{\dim(\Lie(\cP^i))}\cE(\Lie(\cG))$, so the degree of all of these line bundles is bounded, so the procedure must eventually produce a $\cB$-reduction satisfying $\deg(\cE_{\cB}(a_i))\geq 0$ for all simple roots. This proves (3).

Let us deduce the Birkhoff decomposition (4). In the case of constant groups, this is usually deduced from the  decomposition $G=B^-W_0 B$ and (3). For general $\cG$ the analog of this is provided by \cite[Theorem  4.6.33]{BruhatTits2}: Denote by $G_0:=\cG_0^{red}$ the reductive quotient of the fiber of $\cG$ over $0$ and $B_0\subset G_0$ the image of $\cB$. Denote by $U_0$ the unipotent radical of $B_0$, so that $G_0= U_0 W_0 B_0$. The quoted result says that the inverse image of $B_0$ in $\cG$ is $I(0)$. By construction of $\cG$, elements of $U_0(k)$ can be lifted to $\cU(\bP^1)$. Thus we have $\cG(k[[s]])=U_0 W_0 I(0)$. 

Similarly, by our construction of $\cG$ this decomposition implies $\cG(\bA^1)=I^-(0)W_0 U_0$. 

For dominant $t\in \cT(k((s)))$ and $b\in N_0$ we have $t^{-1}bt \in I(0)$. Thus using (3) we find $\cG(k((s)))=I^-(0) \cT(k((s))) \cG(k[[s]])$. Now we want to argue in the same way, decomposing $\cG(k[[s]])$. For any $t$ we can choose $w\in W_0$ such that $wtw^{-1}$ is dominant. Then choose $B_0^\prime \subset G_0$ as $wB_0w^{-1}$ and write $G_0=B_0^\prime W_0 B_0$. Then we can use the same argument as before to deduce (4) (still assuming $k$ to be algebraically closed).

Let us deduce the case that $k$ is a finite field. First assume that $\cG$ splits over the totally ramified covering $[N]\colon\bP^1\to \bP^1$. The embedding $\widetilde{W}\to \Gr_{\cG,x}$ is then defined over $k$, so that all geometric points of $\Bun_{\cG}$ are defined over $k$. Moreover we claim that the automorphism group of any $\cG(0,0)$-bundle is connected. Once we show this, (1),(3) and (4) follow over $k$ by Lang's theorem.  

We have $H^0(\bP^1,\cT)\incl{} H^0(\bP^1,[N]_*(\bG_m^r))=\bG_m^r$. Therefore, $H^0(\bP^1,\cT)\cong \cT_0^{red}$, which is a connected group. Also for any $w\in \widetilde{W}$, the automorphism group of the corresponding $\cG(0,0)$-bundle over $\overline{k}$ is $I^-(0) \cap w I(0) w^{-1}$. This group admits $H^0(\bP^1,\cT)$ as a quotient and the kernel is a product of root subgroups for affine roots, which are connected as well.

The general case follows form this by Galois-descent, using \cite[Remark 9]{HainesRapoport} that the Iwahori-Weyl-group can be computed as the Galois-invariants in the Iwahori-Weyl-group over the separable closure of $k$.
\end{proof}
In the above proof we used the  following special case. Denote by $SU_3$ the quasi-split unitary group for the covering $[2]\colon\bP^1\to \bP^1$. This can be described as the special unitary group for the hermitian form $h(x_1,x_2,x_3)=x_1x_3^\sigma + x_2x_2^\sigma + x_3x_1^\sigma$. Denote by $PSU_3$ the corresponding adjoint group.
\begin{lemme}\label{SU3Lemma}
Any $\PSU_3$ bundle $\cP$ has a reduction $\cP_\cB$ to $\cB$ such that for the positive root $\alpha$ we have $\deg(\cP_\cB(\alpha))\geq 0$. 
\end{lemme}
\begin{proof}
Define $GU_3$ to be the group obtained from $SU_3$ by extending the center of $SU_3$ to $[2]_*\bG_m$, so that thee is an exact sequence $1 \to [2]_*\bG_m \to GU_3 \to \PSU_3 \to 1$. Again, every $PSU_3$ bundle is induced from a $GU_3$-bundle. Such a bundle can be viewed as a rank $3$ vector bundle $\cE$ on the covering $\bP^1\map{[2]} \bP^1$ with a hermitian form with values in a line bundle of the form $[2]^*\cL$. In this case, to give a reduction to $\cB$ it is sufficient to give an isotropic line sub-bundle $\cE_1 \to \cE$, as this defines a flag $\cE_1 \subset \cE_1^\perp \subset \cE$. If $\cE$ is not semistable, then the canonical subbundle of $\cE$ defines an isotropic subbundle of positive degree, so in this case a reduction exists. If $\cE$ is semistable, then the hermitian form defines a global isomorphism $\cE \map{\cong} \sigma^*\cE^\vee \tensor [2]^*\cL$. But such an isomorphism must be constant, so that we can find an isotropic subbundle of degree $\frac{\deg(\cE)}{3}$.
\end{proof}

\section{Appendix: Geometric Satake equivalence and pinnings of the Langlands dual group}\label{a:Sat}

We first need some general properties of tensor functors. Let $\cC$ be a rigid tensor category with a fiber functor $\omega:\cC\to\Vect_F$ where $F$ is a field. Let $H=\uAut^{\otimes}(\omega)$ be the algebraic group over $F$ determined by $(\cC,\omega)$.

For any $F$-algebra $R$, let $\omega_R:\cC\xrightarrow{\omega}\Vect_F\xrightarrow{\otimes_FR}\Mod_R$. Let $\Aut^\otimes(\cC,\omega_R)$ be isomorphism classes of pairs $(\sigma,\alpha)$ where $\sigma:\cC\isom\cC$ is a tensor auto-equivalence, and $\alpha:\omega_R\circ\sigma\Rightarrow\omega_R$ is a natural isomorphism of functors.  Then $\Aut^\otimes(\cC,\omega_R)$ has a natural group structure. Denote by $\uAut^\otimes(\cC,\omega)$ the functor $R\mapsto\Aut^\otimes(\cC,\omega_R)$, which defines a fppf sheaf of groups. On the other hand, let $\uAut(H)$ be the fppf sheaf of automorphisms of the pro-algebraic group $H$ over $F$.

\begin{lemme}\label{l:autgroup}
There is a natural isomorphism of fppf sheaves of groups $\uAut^\otimes(\cC,\omega)\isom\uAut(H)$. In particular, we have a natural isomorphism of groups $\Aut^\otimes(\cC,\omega)\isom\Aut(H)$, which induces an isomorphism of groups $[\Aut^\otimes(\cC)]\isom\Out(H)$, where $[\Aut^\otimes(\cC)]$ is the set of isomorphism classes of tensor auto-equivalences of $\cC$.
\end{lemme}
On the level of $F$-points, a pair $(\sigma,\alpha)\in\Aut^\otimes(\cC,\omega)$ gives the following automorphism of $H=\uAut^\otimes(\omega)$: it sends $h:\omega\Rightarrow\omega$ to the natural transformation 
\begin{eqnarray*}
\omega\stackrel{\alpha}{\Rightarrow}\omega\circ\sigma\stackrel{h\circ\id_\sigma}{\Longrightarrow}\omega\circ\sigma\stackrel{\alpha^{-1}}{\Rightarrow}\omega.
\end{eqnarray*}

More generally, suppose we are given a tensor functor $\Phi:\cC\to\cC'$ into another rigid tensor category $\cC'$, we can similarly define a sheaf of groups $\uAut^\otimes(\cC,\Phi)$.

Let $\omega':\cC'\to\Vect$ be a fiber functor and $\omega=\omega'\circ\Phi$. Then there is a natural homomorphism $\omega'_*:\uAut^\otimes(\cC,\Phi)\to\uAut^\otimes(\cC,\omega)=\uAut(H)$ by sending $(\sigma,\alpha)\in\Aut^\otimes(\cC,\Phi)$ to $(\sigma,\id_{\omega'}\circ\alpha)\in\Aut^\otimes(\cC,\omega)$. In other words, $\uAut^\otimes(\cC,\Phi)$ acts on the pro-algebraic group $H$.

On the other hand, we have natural homomorphism of pro-algebraic groups $\Phi^*:H'=\uAut^\otimes(\omega')\to H=\uAut^\otimes(\omega)$.

\begin{lemme}\label{l:inv}
The homomorphism  $\Phi^*:H'\to H$ factors through $H'\to H^{\uAut^\otimes(\cC,\Phi)}\subset H$.
\end{lemme}

Now we consider the normalized semisimple Satake category $\cS$ in \S\ref{ss:Hecke}. Following \cite{MV} and \cite{Ginz}, we use the global section functor $h$ to define $\dualG=\uAut^\otimes(h)$ and get the geometric Satake equivalence $\cS\cong\Rep(\dualG)$.

\begin{lemme}\label{l:outSatakeH}
\begin{enumerate}
\item []
\item There is a natural homomorphism $\Aut(G)\to\Aut^\otimes(\cS, h)\cong\Aut(\dualG)$ which factors through $\widetilde{\iota}:\Out(G)\to\Aut(\dualG)$;
\item There is a natural pinning $\dagger=(\dualB,\dualT,\{\bbx_{\alpha^\vee_i}\})$ of $\dualG$ preserved by the $\Out(G)$-action via $\widetilde{\iota}$. Let $\Aut^\dagger(\dualG)$ be the automorphism group of $\dualG$ fixing this pinning. Then $\widetilde{\iota}$ induces an isomorphism $\iota:\Out(G)\isom\Aut^\dagger(\dualG)$.
\end{enumerate}
\end{lemme}
\begin{proof}
(1) The $\Aut(G)$-action on $(\cS, h)$ is induced from its action on $\Gr_G$. Since objects in $\cS$ carry $G$-equivariant structures under the conjugation action of $G$ on $\Gr_G$, inner automorphisms of $G$ acts trivially on $(\cS, h)$, i.e., the $\Aut(G)$-action on $(\cS, h)$ factors through $\Out(G)$.

(2) We first need to exhibit a pinning of $\dualG$ which is preserved by the $\Aut(G)$-action. For this, we need to give a maximal torus $\dualT\subset\dualG$, a cocharacter $2\rho\in\xch(T)=\xcoch(\dualT)$ (half the sum of positive coroots in the pinning), and a principal nilpotent element $e\in\dualg$ (the sum of simple root vectors). Equivalently, we need to 
\begin{itemize}
\item[(P1)] Factor the fiber functor $h$ into a tensor functor
\begin{equation*}
h=\bigoplus_{\mu\in\xcoch(T)} h^\mu:\cS\to\Vect^{\xcoch(T)}\xrightarrow{\textup{forget}}\Vect.
\end{equation*}
\item[(P2)] Find a a tensor derivation $e:h\to h$ (i.e., $e(K_1\ast K_2)=e(K_1)\otimes\id_{K_2}+\id_{K_1}\otimes e(K_2)$) which sends $h^\mu$ to $\oplus_ih^{\mu+\alpha_i}$, such that for each simple root $\alpha_i$, the component $h^\mu\to h^{\mu+\alpha_i}$ is nonzero as a functor.
\end{itemize}
The factorization $h^\mu$ is given by Mirkovi\'c-Vilonen's ``weight functors'' \cite[Theorem  3.5, 3.6]{MV}. They also proved that $h=\bigoplus_ih^i$ (where $h^i$ is the sum of $h^\mu$ with $\jiao{2\rho,\mu}=i$) coincides with the cohomological grading of $h=H^*(\Gr_G,-)$. The tensor derivation $e$ is given by the cup product with $c_1(\cL_{\det})\in H^2(\Gr_G,\Ql)$, where $\cL_{\det}$ is the determinant line bundle on $\Gr_G$.

The action of $\Aut^\dagger(G)\subset\Aut(G)$ on $\cS$ permutes the weight functors $h^{\mu}$ in the same way as it permutes $\mu\in\xcoch(T)$, preserves the $\jiao{2\rho,\mu}$ (hence preserves $2\rho\in\xch(T)=\xcoch(\dualT)$), and commutes with $e=c_1(\cL_{\det})$. Therefore,  the $\Aut^\dagger(G)\isom\Out(G)$-action on $\cS$ preserves the above pinning.

An element $\sigma\in\Aut^\dagger(G)$ induces a dual automorphism $\dsigma$ of the Dynkin diagram of $\dualG$. Since $\sigma^*\IC_\mu\cong\IC_{\sigma^{-1}(\mu)}$, the self-equivalence $\sigma^*$ of $\cS\cong\Rep(\dualG)$ is isomorphic to the self-equivalence of $\Rep(\dualG)$ induced by the pinned automorphisms of $\dualG$ given by the dual automorphism $\dsigma^{-1}$ on the Dynkin diagram of $\dualG$. This proves
\begin{equation*}
\Out(G)\to[\Aut^\otimes(\cS)]\cong[\Aut^\otimes(\Rep(\dualG))]\cong\Out(\dualG)
\end{equation*}
is an isomorphism (the last isomorphism follows from Lemma \ref{l:autgroup}). Hence $\iota:\Out(G)\to\Aut^\dagger(\dualG)$ is also an isomorphism.
\end{proof}

For our purpose in \S\ref{s:mono}, we shall also need a different fiber functor $\omega_\phi$ defined in \eqref{omphi}. Recall that $S_\phi$ is the stabilizer of $\phi$ under the action of $T\rtimes\Aut^\dagger(G)\times\grot$, and $S^1_\phi(k)=\ker(S_\phi(k)\to\grot(k))$. Recall from \eqref{S1actondG} that we have an action of $S^1_\phi(k)$ on $\dG$.
 
\begin{lemme}\label{l:outSatake}
There is a natural pinning $\ddagger=(\dB,\dT,\{\leftexp{\phi}{\bbx}_{\alpha^\vee_i}\})$ of $\dG$ which is preserved by the $S^1_\phi(k)$-action.
\end{lemme}
\begin{proof}
Using \eqref{conv}, we can rewrite the fiber functor $\omega_\phi$ as
\begin{equation*}
\cS\xrightarrow{\Psi}\Perv_{I_0}(\Fl_G)\xrightarrow{V_\phi}\Vect
\end{equation*}
where $\Fl_G=G((t))/I_0$ is the affine flag variety at $\{0\}$, $\Psi:\cS\to\Perv(I_0\backslash G((t))/I_0)=\Perv_{I_0}(\Fl)$ is the nearby cycles functor of Gaitsgory \cite{Gaitsgory}, and $V_\phi(K):=R\Gamma_c(\Fl_G, K\otimes\pr^*_1A_\phi)$ as in \eqref{conv}. To exhibit a pinning of $\dG$, we need to find analogs of (P1) and (P2) as in the proof of Lemma \ref{l:outSatakeH}(2).

According to Arkhipov-Bezrukavnikov \cite[Theorem 4]{ArkhipovBezrukavnikov}, each object $\Psi(K)$ admits a $\xcoch(T)$-filtration with Wakimoto sheaves as associated graded pieces. More precisely, they constructed a functor
\begin{equation}\label{Waki}
\bigoplus_{\mu\in\xcoch(T)}W^\mu\circ\Psi(-):\cS\to\Vect^{\xcoch(T)}.
\end{equation}
Here for $\Psi(K)\in\Perv_{I_0}(\Fl_G)$ with Wakimoto filtration $\Psi(K)_{\leq\mu}$, we write $\gr_\mu \Psi(K)=\Psi(K)_{\leq\mu}/\Psi(K)_{<\mu}=J_\mu\otimes W^\mu(K)$, where $J_\mu$ is the Wakimoto sheaf \cite[\S3.2]{ArkhipovBezrukavnikov} and $W^\mu(K)$ is a vector space with Frobenius action. In \cite[Theorem 6]{ArkhipovBezrukavnikov}, it is proved that \eqref{Waki} is tensor.

For each $\mu$, let $j_\mu:\Fl_\mu=I_0t^\mu I_0/I_0\hookrightarrow\Fl_G$ be the inclusion. We fix the Frobenius structure of $J_\mu$ in the following way: for $\mu$ regular dominant, let $J_\mu=j_{\mu,*}\Ql[\jiao{2\rho,\mu}](\jiao{2\rho,\mu})$; for $\mu$ regular anti-dominant, let $J_\mu=j_{\mu,!}\Ql[\jiao{2\rho,\mu}]$; for general $\mu=\mu_1+\mu_2$ where $\mu_1$ is regular dominant and $\mu_2$ is regular anti-dominant, let $J_\mu=J_{\mu_1}\conv{I_0}J_{\mu_2}$ (where $\conv{I_0}$ is the convolution on $\Fl_G$). It follows from \cite[Lemma 8, Corollary 1]{ArkhipovBezrukavnikov} that $J_\mu$ is well-defined. This normalization makes sure that in the composition series of $J_\mu$, $\delta=\IC_1$ appears exactly once, with multiplicity space $\Ql$ as a trivial Frobenius module (see \cite[Lemma 3(a)]{ArkhipovBezrukavnikov} with obvious adjustment to the mixed setting). Since all $\IC_\tilw$ are killed by $V_\phi$ except $\tilw=1$, we conclude that
\begin{equation}\label{Wakipure}
V_\phi(J_\mu)=\Ql \textup{ as a trivial Frobenius module for all }\mu\in\xcoch(T).
\end{equation}

\begin{claim}
For $K\in\cS$, $W^\mu(K)$ is pure of weight $\jiao{2\rho,\mu}$. In fact, we have a natural isomorphism of funtors $h^\mu\cong W^\mu$ ($h^{\mu}$ is the weight functor in \cite[Theorem  3.5, 3.6]{MV}, which was used in the proof of Lemma \ref{l:outSatakeH}).
\end{claim}
\begin{proof}
We first recall the definition of the weight functors in \lc. For every $\mu\in\xcoch(T)$, let $\frS_\mu\subset\Gr_G$ be the $U((t))$-orbit containing $t^\mu$. The weight function defined in \lc is $h^{\mu}(K)=H^*_c(\frS_\mu,K)$, which is concentrated in degree $\jiao{2\rho,\mu}$.

Let $\pi:\Fl_G\to\Gr_G$ be the projection, then $\pi^{-1}(\frS_\mu)=\sqcup_{w\in W}\tfrS_{\mu w}$, where $\tfrS_{\mu w}\subset\Fl_G$ is the $U((t))$-orbit containing $t^\mu w$. We have natural isomorphisms (for $K\in\cS$)
\begin{equation*}
h^\mu(K)=H^*_c(\frS_\mu,\pi_!\Psi(K))\cong H^*_c(\pi^{-1}(\frS_\mu),\Psi(K))\isom H^*_c(\tfrS_\mu,\Psi(K))\isom W^\mu(K).
\end{equation*}
Here the first equality follows from $\pi_!\Psi(K)=K$, and the last two isomorphisms follow from \cite[Theorem 4(2)]{ArkhipovBezrukavnikov} (with extra care about the Frobenius structure).

Since $h^i(K)=H^i(\Gr_G,K)$ is pure of weight $i$, and $h^\mu(K)$ is a direct summand of $h^{\jiao{2\rho,\mu}}$, $h^\mu(K)$ is pure of weight $\jiao{2\rho,\mu}$. Hence $W^\mu(K)$ is also pure of weight $\jiao{2\rho,\mu}$.
\end{proof}

Now we construct a natural isomorphism $\bigoplus_\mu W^\mu\cong\omega_\phi$. For each $K\in\cS$, the Wakimoto filtration on $\Psi(K)$ gives a spectral sequence calculating $V_\phi(\Psi(K))$ with $E_1$ page $V_\phi(\gr_\mu\Psi(K))$. By \eqref{Wakipure},
\begin{equation}\label{Vpure}
V_\phi(\gr_\mu\Psi(K))=V_\phi(J_\mu)\otimes W^\mu(K)=W^\mu(K)
\end{equation}
is concentrated in degree 0, the spectral sequence degenerates at $E_1$. The limit of the spectral sequence gives a Wakimoto filtration on the  Frobenius module $\omega_\phi(K)=V_\phi(\Psi(K))$, which we denote by $w_{\leq\mu}$. By the Claim and \eqref{Vpure}, this filtration refines the weight filtration $w_{\leq i}$ on $\omega_\phi(K)$:
\begin{equation*}
\gr^w_{i}\omega_\phi(K)=\bigoplus_{\jiao{2\rho,\mu}=i}\gr^w_{\mu}\omega_\phi(K).
\end{equation*}
Since $\omega_\phi(K)$ is a Frobenius module in $\Ql$-vector spaces, the weight filtration splits canonically. Therefore, the Wakimoto filtration $w_{\leq\mu}\omega_\phi(K)$ also splits canonically. This gives a canonical isomorphism $\bigoplus_\mu W^\mu\cong\omega_\phi$.

Finally, the principal nilpotent element is given by the logarithm of the monodromy action (of a topological generator of the tame inertia group at 0) on the nearby cycles (\cite[Theorem 2]{Gaitsgory}) $M_K:\Psi(K)\to\Psi(K)(-1)$.

The action of $S^1_\phi(k)$ commutes with the nearby cycle functor, therefore commutes with the monodromy $M_K$. It permutes the Wakimoto sheaves, hence permutes the functors $W^\mu$ through the action of $S^1_\phi(k)\to\Aut^\dagger(G)$ on $\xcoch(T)$, and it preserves $\jiao{2\rho,\mu}$. Therefore the $S^1_\phi(k)$-action preserves the pinning $\ddagger$. 
\end{proof}


\section{Appendix: Quasi-minuscule combinatorics}\label{a:qmcom} We assume $G$ is almost simple of rank at least 2.

Let $\theta$ be the highest root (which is a long root) and $\thv$ the corresponding coroot (which is a short coroot). The set of roots $\Phi$ is partitioned into $\Phi^\theta_n=\{\alpha\in\Phi|\jiao{\alpha,\theta^\vee}=n\}$, $n=0,\pm1,\pm2$, with $\Phi^\theta_{\pm2}=\{\pm\theta\}$. The roots in $\Phi^\theta_1$ are coupled into pairs $(\alpha,\beta)$ with $\alpha+\beta=\theta$.



Let $V_{\thv}$ be the irreducible representation of $\dualG$ with highest weight $\thv$. The nonzero weights of $V_\thv$ are the short roots of $\dualG$, each with multiplicity 1. Let $\{e=\sum_{i}x_i,2\rho,f=\sum_iy_i\}\in\dualg$ be the principal $\frsl_2$-triple, where $x_i\in\dualg_{\alpha^\vee_i}, y_i\in\dualg_{-\alpha^\vee_i}$ are nonzero and $2\rho\in\dualt$ is the sum of positive coroots of $\dualG$.

\begin{lemme}\label{l:samenumber} The following numbers are the same:
\begin{enumerate}
\item $\dim V_\thv(0)$, where $V_\thv(0)$ is the zero weight space under the $2\rho$-action, and is also the zero weight space under the $\dualT$-action;
\item $\dim V^e_\thv$, where $V^e_\thv=\ker(\ad(e)|V_\thv)$;
\item the number of short simple roots of $\dualG$ (the {\em short rank} $r_s(\dualG)$);
\item $\#(W\thv)/h=\#\{\textup{short roots of }\dualG\}/h$ (Here $h$ is the Coxeter number of $W$).
\end{enumerate}
\end{lemme}
\begin{proof}
Under the principal $\frsl_2$ action, $V_\thv$ can be decomposed as a sum of irreducible representations of $\frsl_2$:
\begin{equation*}
V_\thv=\sum_{i=1}^{r_s(\dualG)}\Sym^{2\ell_i}(\St)
\end{equation*}
where $\St$ is the 2-dimensional representation of $\frsl_2$. Since the weights of the $2\rho$-action on $V_{\thv}$ are even, only even symmetric powers of $\St$ appear in $V_\thv$.  Since each $\Sym^{2\ell_i}(\St)$ contributes 1-dimension to both $V_\thv(0)$ and $V^e_\thv$, we have $r_s(\dualG)=\dim V_\thv(0)=\dim V^e_\thv$. This proves the equality of the numbers in (1) and (2).

Let $V_\thv(n)$ be the weight $n$-eigenspace of the $\rho$-action, then
\begin{equation*}
V_\thv(n)=\sum_{\alpha^\vee \textup{short},\jiao{\rho,\alpha^\vee}=n}V_\thv(\alpha^\vee).
\end{equation*}
In particular, $\dim V_\thv(1)$ is the number of short simple roots of $\dualG$. The map $e:V_\thv(0)\to V_\thv(1)$ is clearly surjective. It is also injective, because if $v\in V_\thv(0)$, $ev=0$ means $\dualg_{\alpha^\vee_i}v=0$ for all simple $\alpha^\vee_i$, i.e., $v$ is a highest weight, contradiction. This proves
\begin{equation}\label{V0V1}
\dim V_\thv(0)=\dim V_\thv(1)=r_s(\dualG).
\end{equation}

It remains prove (4) is the same as the rest. The argument is similar to that of \cite[\S 6.7]{Kostant}, where Kostant considered the adjoint representation instead of $V_\thv$. We only give a sketch. Let $\zeta$ be a primitive $h$-th root of unity and let $P=\rho(\zeta)\in\dualG$. Let $\gamma^\vee$ be the highest root of $\dualG$ and $z=e+\bbx_{-\gamma^\vee}$. Then $z$ is a regular semisimple element in $\dualg$, and $\Ad(P)z=\zeta z$. Let $\dualG_z$ be the centralizer of $z$ (which is a maximal torus), then $P\in N_{\dualG}(\dualG_z)$ and its image in the Weyl group is a Coxeter element. Choosing a basis $u_i$ for the highest weight line in $\Sym^{2\ell_i}(\St)\subset V_\thv$, then there exists a unique $v_i\in V_\thv(\ell_i-h)$ such that $u_i+v_i\in V_\thv^z$ (kernel of the $z$-action on $V_\thv$). Then $\{u_i+v_i\}_{1\leq i\leq r_s(\dualG)}$ form a basis of $V_\thv^z$, with eigenvalues $\zeta^{\ell_i}$ under the action of $\Ad(P)$. 

The representation $V_\thv$ of $\dualG$ is clearly self-dual. According to zero and nonzero weights under $\dualg_z$, we can write $V_\thv=V_\thv^z\oplus V'$ as $N_{\dualG}(\dualG_z)$-modules. Any self-duality $V_\thv\cong V_\thv^*$ of $\dualG$-modules necessarily restricts to a self-duality $V_\thv^z\cong (V_\thv^z)^*$ as $N_{\dualG}(\dualG_z)$-modules. Hence the eigenvalues of $\Ad(P)$ on $V^z_\thv$ are invariant under inversion: i.e., the multi-set $\{\zeta^{\ell_i}\}_{1\leq i\leq r_s(\dualG)}$ is invariant under inversion. This implies $\sum_{i}\ell_i=r_s(\dualG)h/2$. Hence $\dim V_\thv=\sum_i(2\ell_i+1)=(h+1)r_s(\dualG)$. Since the nonzero weight spaces of $V_\thv$ are indexed by short roots of $\dualG$, hence by the orbit $W\thv$, we get $\#(W\thv)=\dim V_\thv-\dim V_\thv(0)=hr_s(\dualG)$. This proves that (4) coincides with the rest of the numbers.

\end{proof}


\section{Appendix: The adjoint Schubert variety for $G_2$}\label{a:g2} In this section we assume char$(k)>3$. Recall (see tables in \cite{Ca}) that $G_2$ has four nonregular unipotent orbits: 
\begin{itemize}
\item the {\em subregular} orbit containing a generic element in the unipotent radical of $\Pnth$; it has dimension 10;
\item the orbit containing $U^\times_{-\gamma}$, which has dimension 8;
\item the orbit containing $U^\times_{-\theta}$, which has dimension 6;
\item the identity orbit, which has dimension 0.
\end{itemize}
The $G$-orbits of $\Gr^\triv_{\leq\gav}$ for $G_2$ turns out to be closely related to these unipotent orbits. More precisely,

\begin{lemme}\label{l:g2adj}
\begin{enumerate}
\item There are four $\Ad(G)$-orbits on $\Gr^\triv_{\leq\gav}$, 
\begin{equation*}
\Gr^\triv_{\leq\gav}=\Grsubreg\bigsqcup\Ad(G)U^\times_{-\gamma,-1}\bigsqcup\Gr^{\triv}_\thv\bigsqcup\{\star\}
\end{equation*}
of dimensions $10,8,6$ and $0$ respectively, which, under the evaluation map $\ev_{\tau=1}$, map onto the four nonregular unipotent orbits.
\item The morphism
\begin{eqnarray*}
\nu:G\twtimes{\Pnth}(\prod_{\jiao{\beta,-\thv}\geq1}U_{\beta,-1})&\to&\Gr^\triv_{\leq\gav}\\
(g,\prod u_{\beta}(c_\beta\tau^{-1}))&\mapsto&\Ad(g)(\prod u_{\beta}(c_\beta\tau^{-1}))
\end{eqnarray*}
is a resolution. Its fibers over the $G$-orbits are:
\begin{itemize}
\item $\nu$ is an isomorphism over $\Grsubreg$;
\item $\nu^{-1}(u_{-\gamma}(\tau^{-1}))\cong\PP^1$;
\item $\nu^{-1}(u_{-\theta}(\tau^{-1}))$ is a projective cone over $\PP^1$ (it contains a point, whose complement is a line bundle over $\PP^1$);
\item $\nu^{-1}(\star)\cong G/\Pnth$.
\end{itemize}
\item The morphism
\begin{equation*}
\nu':G\twtimes{\Png}(\prod_{\jiao{\beta,-\gav}\geq2}U_{\beta,-1})\to\Gr^\triv_{\leq\gav}
\end{equation*}
defined similarly as $\nu$, is a resolution of the closure of $\Ad(G)U^\times_{-\gamma,-1}$. Its fibers over the $G$-orbits are:
\begin{itemize}
\item $\nu$ is an isomorphism over $\Ad(G)U^\times_{-\gamma,-1}$;
\item $\nu^{-1}(u_{-\theta}(\tau^{-1}))\cong\PP^1$;
\item $\nu^{-1}(\star)\cong G/\Png$.
\end{itemize}
\end{enumerate}
\end{lemme}

\begin{proof}[Proof of Theorem \ref{th:Eulerchar}(2) for $G=G_2$.]
By the same reduction steps as in the case of other types, we reduce to showing
\begin{equation}\label{chig2}
\chi_c(\Gr^{\triv,a=0}_\gav,\ev^*J)=-r_s(G)=-1.
\end{equation}

\noindent{\bf Step I}. $\chi_c((\Ad(G)U^\times_{-\gamma,-1})^{a=0},\ev^*J)=-1$.

For any short root $\beta$, let $V_\beta=\prod_{\jiao{\beta',\beta^\vee}\geq2}U_{\beta',-1}$. Then the source of $\nu'$ has a Bruhat decomposition
\begin{equation*}
G\twtimes{\Png}\Vng=\bigsqcup_{\beta\textup{ short root}}\Ad(U)V_\beta.
\end{equation*}
The following Claim can be proved similarly as the claim in \S\ref{ss:adjproof}.
\begin{claim}
For a short root $\beta$,
\begin{equation*}
\chi_c(\Ad(U)V^{>-\theta}_\beta,\ev^*J)=\begin{cases}
0 & \Phi^\beta_{\geq2}\textup{ contains a simple root}\\
1 & \textup{otherwise.}\end{cases}
\end{equation*}
\end{claim}

Looking at the root system $G_2$, there are 3 short roots $\beta$ such that $\Phi^\beta_{\geq2}$ does not contain a simple root. Therefore $\chi_c((G\twtimes{\Png}\Vng)^{a=0},\nu'^*\ev^*J)=3$. On the other hand, by Lemma \ref{l:g2adj}(3), we have
\begin{eqnarray*}
3&=&\chi_c((G\twtimes{\Png}\Vng)^{a=0},\nu'^*\ev^*J)\\
&=&\chi_c(G/\Png)+\chi_c(\PP^1)\chi_c(\Gr^{\triv,a=0}_\thv,\ev^*J)+\chi_c((\Ad(G)U^\times_{-\gamma,-1})^{a=0},\ev^*J)
\end{eqnarray*}
Plugging in $\chi_c(G/\Png)=\#W/W_\gav=6$, $\chi_c(\PP^1)=2$ and $\chi_c(\Gr^{\triv,a=0}_\thv,\ev^*J)=-1$ from the proof of Theorem \ref{th:Eulerchar}(1), we conclude that  $\chi_c((\Ad(G)U^\times_{-\gamma,-1})^{a=0},\ev^*J)=-1$.

\noindent{\bf Step II}. $\chi_c(\Grsubreg^{a=0},\ev^*J)=0$.

For any long root $\alpha$, let $V_\alpha=\prod_{\jiao{\alpha',\alpha^\vee}\geq1}U_{\alpha',-1}$. Then the source of $\nu$ has a Bruhat decomposition
\begin{equation*}
G\twtimes{\Pnth}V_{-\theta}=\bigsqcup_{\alpha\textup{ long root}}\Ad(U)V_\alpha.
\end{equation*}
The following Claim can be proved similarly as the claim in \S\ref{ss:adjproof}.
\begin{claim}
For a long root $\alpha$,
\begin{equation*}
\chi_c(\Ad(U)V^{>-\theta}_\alpha,\ev^*J)=\begin{cases}
0 & \Phi^\alpha_{\geq1}\textup{ contains a simple root}\\
1 & \textup{otherwise.}\end{cases}
\end{equation*}
\end{claim}

Looking at the root system $G_2$, $\alpha=-\theta$ is the only long root for which $\Phi^\alpha_{\geq1}$ does not contain a simple root. Therefore $\chi_c((G\twtimes{\Pnth}V_{-\theta})^{a=0},\nu^*\ev^*J)=1$. On the other hand, by Lemma \ref{l:g2adj}(2), we have
\begin{eqnarray*}
1&=&\chi_c((G\twtimes{\Pnth}V_{-\theta})^{a=0},\nu*\ev^*J)\\
&=&\chi_c(G/\Pnth)+\chi_c(\PP^1)\chi_c((\Ad(G)U^\times_{-\gamma,-1})^{a=0},\ev^*J)\\
&&+\chi_c(\nu^{-1}(u_{-\theta}(\tau^{-1})))\chi_c(\Gr^{\triv,a=0}_\thv,\ev^*J)+\chi_c(\Grsubreg^{a=0},\ev^*J).
\end{eqnarray*}
Plugging in $\chi_c(G/\Pnth)=\#W/W_\theta=6$, $\chi_c(\PP^1)=2$, $\chi_c(\nu^{-1}(u_{-\theta}(\tau^{-1})))=3$ from Lemma \ref{l:g2adj}(2), $\chi_c((\Ad(G)U^\times_{-\gamma,-1})^{a=0},\ev^*J)=-1$ from Step I,  and $\chi_c(\Gr^{\triv,a=0}_\thv,\ev^*J)=-1$ from the proof of Theorem \ref{th:Eulerchar}(1), we conclude that  $\chi_c(\Grsubreg^{a=0},\ev^*J)=-1$.

Combining Step I and II, since $\Gr^{\triv,a=0}_\gav=\Grsubreg^{a=0}\sqcup(\Ad(G)U^\times_{-\gamma,-1})^{a=0}$, we get \eqref{chig2}. This proves Theorem \ref{th:Eulerchar}(2) in the case of $G_2$.
\end{proof}

\tableofcontents
\end{document}